\documentclass{amsart}%
\usepackage{amsfonts}%
\usepackage{amsmath}%
\usepackage{amssymb}%
\usepackage{graphicx}
\usepackage{amsthm}
\usepackage{setspace}
\usepackage{color}
\usepackage{newclude}
\usepackage[top = 1in, bottom = 1in, left = 1.2in, right =
1.2in]{geometry}
\providecommand{\U}[1]{\protect\rule{.1in}{.1in}}

\newtheorem{theorem}{Theorem}
\theoremstyle{plain}

\newtheorem{definition}{Definition}[section]

\newtheorem{lemma}{Lemma}[section]

\newtheorem{proposition}{Proposition}[section]
\newtheorem{remark}{Remark}[section]

\numberwithin{equation}{section}

\begin{document}
\title[]{On the structure of the singular set for the kinetic Fokker-Planck
equations in domains with boundaries.}


\author{Hyung Ju Hwang}
\address{Department of Mathematics, Pohang University of Science and
Technology, Pohang, GyungBuk 790-784, Republic of Korea}
\email{hjhwang@postech.ac.kr}
\thanks{}

\author{Juhi Jang}
\address{Department of Mathematics, University of Southern California,
Los Angeles, CA 90089, USA and Korea Institute for Advanced Study, Seoul, Korea}
\email{juhijang@usc.edu}
\curraddr{}
\thanks{}

\author{Juan J. L. Vel\'{a}zquez}
\address{Institute of Applied Mathematics, University of Bonn, Endenicher
Allee 60, 53115 Bonn, Germany}
\email{velazquez@iam.uni-bonn.de}
\curraddr{}
\thanks{}

\date{today}



\begin{abstract}
In this paper we compute asymptotics of solutions of the kinetic Fokker-Planck equation with inelastic boundary conditions
which indicate that the solutions are nonunique if $r < r_c$. The nonuniqueness is due to the fact that different solutions can interact in a different manner 
with a Dirac mass which appears at the singular point $(x,v)=(0,0)$. In particular, this nonuniqueness explains 
the different behaviours found in the physics literature
for numerical simulations of the stochastic differential equation associated to the kinetic Fokker-Planck
equation. The asymptotics obtained in this paper will be used in a companion paper \cite{HVJ2} to prove rigorously 
nonuniqueness of solutions for the kinetic Fokker-Planck equation with inelastic boundary conditions. 
\end{abstract}

\maketitle


%
%
%



\section{Introduction}




A general feature of several kinetic equations is the fact that their
solutions in a domain $\Omega $ with boundaries cannot be infinitely
differentiable at the points of the phase space $\left( x,v\right) $ for
which $x\in \partial \Omega $ and $v$ is tangent to $\partial \Omega ,$ even
if the initial data are arbitrarily smooth. This set of points is usually
denoted as the singular set. This property of the solutions of kinetic equations
was found by Guo in his study of the Vlasov-Poisson system (cf. \cite{G1}).

In the case of the one-species Vlasov-Poisson system in 3-dimensional
bounded domains with specular reflection boundary conditions, the convexity
of the boundary plays a crucial role in determining the regularity of
solutions (cf. \cite{Hw}, \cite{HV}). If the boundary of the domain is flat
or convex, then one can show that a trajectory on the singular set is
separated from regular one in the interior, so that classical solutions
exist globally in time and are $C^{1}$ in phase space variables except on
the singular set. 
However, if a domain is nonconvex, then even in the Vlasov equation without
interactions between the particles and with inflow boundary condition, a
classical $C^{1}$ solution fails to exist in general (cf. \cite{G0}). 

The lack of smoothness of solutions of the Vlasov-Poisson system is closely
related to the fact that the flow describing the evolution of the
characteristic curves is not a $C^{2}$\ function at the singular set,
something which has been observed in the dynamics of billiards, even if in
this case there are not gravitational or electrical fields (cf. \cite{CM},
\cite{Si2}).

In the case of the Boltzmann equation near Maxwellian, under inflow,
diffuse, specular reflection and bounce-back boundary conditions, it is
shown in \cite{GKTT} that if a domain is strictly convex, then the solution
is (weighted) $C^{1}$ away from the singular set. On the other hand, if
there is a concave boundary point, one can construct an initial condition
which induces a discontinuous solution in any given time interval (cf. \cite%
{Kim}).

A type of equation in which the regularizing effects can be expected to be
stronger than those mentioned above among the kinetic equations are the
Fokker-Planck equations. These equations contain second order operators in
some of the variables (the velocities), but due to the specific form of the
transport terms in the remaining variables (positions) the solutions of
these equations are smooth in all their variables. This hypoelliptic feature 
has been
extensively studied in the mathematics literature (cf. for instance \cite%
{HBook}, \cite{RS} as well as additional references in \cite{HVJ2}). 
From a PDE point of view, the solutions
are expected to be more regular due to the smoothing effect of the random
force. This was clarified by Bouchut in the case of the 3-dimensional whole
space problem for the Vlasov-Poisson-Fokker-Planck equations \cite{Bou}.

On the other hand, boundary value
problems for the Fokker Planck equation have been considered in many
physical situations, including diffusion controlled reactions and the
dynamics of semiflexible polymers (cf. \cite{B1}, \cite{B2}) and in general
in problems involving the dynamics of Brownian motion described by
Uhlenbeck-Ornstein processes in the presence of boundaries (cf. for instance
\cite{KB}, \cite{MW1}, \cite{MW2}). 
Despite physical importance, in the presence of the boundary, there had been not much
rigorous improvement for a few decades.

It is interesting to remark that the onset of singular behaviours at the
singular set takes place even for kinetic equations associated to
hypoelliptic operators. As indicated above, one of the best known
hypoelliptic operators is the one associated to the kinetic Fokker-Planck
equation, also known as Kolmogorov equation:%
\begin{equation}
\partial _{t}P+v\cdot \nabla _{x}P=\Delta _{v}P\ \ ,\ \ P=P\left(
x,v,t\right) \   \label{I0E1}
\end{equation}%
where $t\in \left( 0,T\right) \subset \mathbb{R}$ and $\left( x,v\right) \in
U\subset \mathbb{R}^{N}\times \mathbb{R}^{N}$ for some $N\geq 1$ and some
domain $U$. It is well known that any solution of (\ref{I0E1}) is $C^{\infty
}\left(  U  \times \left( 0,T\right)\right) $ in spite of the fact that the
second order operator $\partial _{vv}$ acts only in the variable $v$ (cf.
\cite{HBook}). However, it turns out that in domains $\Omega $ with
boundaries, singular sets in the sense of kinetic equations arise for some
classes of boundary conditions. More precisely, there are some classes of
boundary conditions, which arise naturally from the physical interpretation
of the equation (\ref{I0E1}) for which the hypoellipticity property can be
extended to a large fraction of the boundary $\partial U,$ but it fails
along some subsets of $\partial U.$ In those sets the solutions are just H%
\"{o}lder continuous. One example of this situation is the following.
Suppose that we restrict the analysis to the case $N=1$ and that we solve (%
\ref{I0E1})\ in the domain $U=\left( 0,1\right) \times \mathbb{R}$ . In
these cases the singular set reduces either to the two points $\left\{
\left( 0,0\right) ,\left( 1,0\right) \right\} $. Suppose that we study this
problem with the absorbing boundary condition:%
\begin{equation}
P\left( 0,v,t\right) =0,\text{ for }v>0,\ P\left( 1,v,t\right) =0,\ \text{
for }v<0,\ t>0.\text{\ }  \label{I0E2}
\end{equation}%
The problem (\ref{I0E1}), (\ref{I0E2}) has been studied in \cite{HVJ} where
it has been proved that the solutions are at most H\"{o}lder continuous in a
neighbourhood of the singular set $\left\{ \left( 0,0\right) ,\left(
1,0\right) \right\}.$ 

In this paper, we will consider the following problem: 
\begin{equation}
\partial _{t}P+v\partial _{x}P=\partial _{vv}P\ \ ,\ \ P=P\left(
x,v,t\right) ,\ \ \ x>0,\ \ v\in \mathbb{R}\ \ ,\ \ t>0  \label{S0E1}
\end{equation}%
with the so-called inelastic boundary condition:%
\begin{equation}
P\left( 0,-v,t\right) =r^{2}P\left( 0,rv,t\right) \ \ ,\ \ v>0\ \ ,\ \ t>0
\label{S0E2}
\end{equation}%
where $0<r<\infty .$ We are particularly interested in the case in which $%
r\leq 1.$ Notice that the physical meaning of (\ref{S0E2}) is that the
particles arriving to the wall $\left\{ x=0\right\}$ with a velocity $-v$
bounce back to the domain $\left\{ x>0\right\}$ with a new velocity $v.$
Notice that in the case $r=1$ the particles bounce elastically, but for $%
r<1, $ the collisions are inelastic and the particles lose a fraction of
their energy in the collisions. 
The problem (\ref{S0E1}), (\ref{S0E2}) will be solved 
with the initial condition:%
\begin{equation}
P\left( x,v,0\right) =P_{0}\left( x,v\right)  \label{S0E3}
\end{equation}%
where $P_{0}$ is a probability measure in $\left( x,v\right) \in \mathbb{R}%
^{+}\times \mathbb{R}$.

In this case the singular point is $\left( x,v\right) =\left( 0,0\right) .$
We are interested in the construction of measure valued solutions of the
problem (\ref{S0E1})-(\ref{S0E3}). It turns out that for some choices of $r$
we might have $\int_{\left\{ \left( 0,0\right) \right\} }P\left(
dxdv,t\right) >0$ for $t>0.$ Moreover, $P\left( x,v,t\right) $ is a $%
C^{\infty }$ function for $\left( x,v\right) \neq \left( 0,0\right) ,$ but
it would be only possible to obtain uniform estimates in some suitable H\"{o}%
lder norms if $\left( x,v\right) $ is close to the singular point. In fact, 
we will obtain information about the behaviour of the solutions near the
singular point for the adjoint problem of (\ref{S0E1}), (\ref{S0E2}) which
demonstrates that the solutions of this problem are at most H\"{o}lder continuous. 

It is possible to give the following physical motivation to the problem (\ref%
{S0E1})-(\ref{S0E3}). Suppose that we have a\ particle $X\left( t\right) $
which moves in the half-line $\left\{ x>0\right\}.$ We will assume that, as
long as $X\left( t\right) >0,$ the dynamics of the particle is similar to
that of a Brownian particle, but that, every time that the particle reaches $%
x=0,$ it bounces inelastically, that is, the absolute value of its velocity
is multiplied by a restitution coefficient $r,\ $with $0<r\leq 1.$ Formally 
we may state that 
the dynamics of the particle is given by the following set
of equations:%
\begin{eqnarray}
\frac{dX\left( t\right) }{dt} &=&V\left( t\right) \ ,\ \ \frac{dV\left(
t\right) }{dt}=\sqrt{2}\eta \left( t\right) \ \ ,\ \text{if}\ X\left(
t\right) >0\ ,\ \   \label{S1E1} \\
X\left( t_{\ast }^{-}\right) &=&0\ \ ,\ \ V\left( t_{\ast }^{-}\right) <0\ \
\Longrightarrow X\left( t_{\ast }^{+}\right) =0\ \ ,\ \ V\left( t_{\ast
}^{+}\right) =-rV\left( t_{\ast }^{-}\right)  \label{S1E2}
\end{eqnarray}%
where $\eta \left( t\right) $ is the white noise stochastic process, which
formally satisfies $\left\langle \eta \left( t\right) \right\rangle =0,$ $%
\left\langle \eta \left( t_{1}\right) \eta \left( t_{2}\right) \right\rangle
=\delta \left( t_{1}-t_{2}\right),$ where $\left\langle \cdot \right\rangle
$ denotes average. It is not a priori obvious if the process defined by
means of the equations (\ref{S1E1}), (\ref{S1E2}) can be given a precise
meaning. However, standard arguments of the Theory of Stochastic Processes
suggest that the probability of finding the particle $\left( X\left(
t\right) ,V\left( t\right) \right) $ in a region of the phase space $\left(
x,v\right) $ is given by the probability density $P\left( x,v,t\right).$

The first mathematical results for the problem (\ref{S1E1}), (\ref{S1E2})
were obtained by McKean (cf. \cite{McK}). More precisely, \cite{McK}
contains in particular the probability distribution of hitting times and
hitting points, defined by means of the solutions of (\ref{S1E1}) with
initial conditions $X\left( 0\right) =0,\ V\left( 0\right) =b>0.$ The
hitting times are defined as $\mathfrak{t}_{1}=\max \left\{ t>0:X\left(
s\right) >0\text{ for }0<s<t\right\} $ and the hitting points are defined by
means of $\mathfrak{h}_{1}=\left\vert V\left( \mathfrak{t}_{1}\right)
\right\vert .$ Notice that a rescaling argument allows us to reduce the
computation of the joint distribution of $\left( \mathfrak{t}_{1},\mathfrak{h%
}_{1}\right) $ for $b>0$ to the case $b=1.$ Suppose that we denote as $%
\left\{ t_{n}\right\} $ the sequence of consecutive times where the
solutions of (\ref{S1E1}), (\ref{S1E2}) with initial values $X\left(
0\right) =X_{0}\geq 0,\ V\left( 0\right) =V_{0}$ satisfy $X\left(
t_{n}\right) =0.$ Then, the results in \cite{McK} imply that $X\left(
t_{n}\right) \rightarrow 0$ as $n\rightarrow \infty ,$ $t_{n}\rightarrow
T<\infty $ with probability one if $r<r_{c},$ where
\begin{equation}
r_{c}=\exp \left( -\frac{\pi }{\sqrt{3}}\right)  \label{S1E3}
\end{equation}%
On the contrary, the results in \cite{McK} imply that that for $r\geq r_{c}$
we have $\left( X\left( t\right) ,\ V\left( t\right) \right) \neq 0$ for any
$t\geq 0.$ Details about these results can be found in \cite{J1}, \cite{J2}.

The results indicated above suggest that the solutions of (\ref{S0E1})-(\ref%
{S0E3}) should have a very different behaviour for $r<r_{c}$ and $r\geq
r_{c}.$ Indeed, in the first case the equations (\ref{S0E1})-(\ref{S0E3})
cannot be expected to define the evolution of $P$ without some additional
information about the behaviour of the solutions at the singular point $%
\left( x,v\right) =\left( 0,0\right) .$ On the contrary, the problem (\ref%
{S0E1})-(\ref{S0E3}) should be able to define a unique dynamics for $P$ if $%
r\geq r_{c}.$ The main goal of this paper as well as the companion paper \cite{HVJ2} 
is to show that we can have
different evolutions for the problem (\ref{S0E1})-(\ref{S0E3}) if $r<r_{c}.$
On the contrary, we will prove that for $r>r_{c}$ no additional information
is required in order to define the unique evolution of the problem.

The problem (\ref{S1E1}), (\ref{S1E2}) has been also considered in the
physics literature \cite{B1,BFG,BK,CSB,FBA,KB}. The paper \cite{CSB} arrives
to conclusions analogous to those which follow from McKean's results by
means of the analysis of some particular solutions of the Fokker-Planck
equation associated to (\ref{S1E1}), (\ref{S1E2}) as well as numerical
simulations. The authors of \cite{CSB} claimed that the solutions of (\ref%
{S1E1}), (\ref{S1E2}) undergo inelastic collapse for $r<r_{c}.$ More
precisely, it was claimed that for $r<r_{c}$ the Brownian particle with
inelastic collisions stops its motion in finite time with probability one,
while for $r>r_{c}$ the Brownian particle with inelastic collisions has
positive velocity for any time.

However, some papers in the physics literature raised doubts about the
results in \cite{CSB}. It was claimed in \cite{FBA}, \cite{A}, on the basis
of numerical simulations, that inelastic collapse for $r<r_{c}$ cannot be
observed.

There is a detailed discussion about these seemingly contradictory results
in \cite{BK}. The conclusion of this paper was that for $r<r_{c}$ the
particle $\left( X\left( t\right) ,V\left( t\right) \right) $ solutions of (%
\ref{S1E1}), (\ref{S1E2}) arrive to $\left( x,v\right) =\left( 0,0\right) $
in finite time. However the particle $\left( X\left( t\right) ,V\left(
t\right) \right) $ does not stay there for later times. The main argument
given in \cite{BK} to support this conclusion is a detailed analysis of a
stationary solution of the Fokker-Planck equation associated to the problem (%
\ref{S1E1}), (\ref{S1E2}). Such stationary solution exhibits a flux of
particles from the point $\left( x,v\right) =\left( 0,0\right) $ to the
region $\left\{ x>0\right\}.$ This flow of particles would be balanced by
the flux of particles arriving towards the origin as predicted in \cite{McK}
by McKean.

Suppose that we denote as $P\left( x,v,t\right) $ the probability density
which gives the probability of finding one particle solving (\ref{S1E1}), (%
\ref{S1E2}) in any region of the phase space $\left( x,v\right) .$ Inelastic
collapse would imply that $P\left( x,v,t\right) \rightarrow \delta \left(
x\right) \delta \left( v\right) $ as $t\rightarrow \infty .$ On the other
hand, the steady solution described in \cite{BK} is different from $\delta
\left( x\right) \delta \left( v\right) $ and this indicates that the
particles arriving to $\left( x,v\right) =\left( 0,0\right) $ do not
necessarily remain there.

We will obtain in this paper different asymptotic formulas for the behaviour of the solutions 
of (\ref{S0E1})-(\ref{S0E3}) near the singular point $(x,v)=(0,0)$. 
These formulas will be used in the paper \cite{HVJ2} to prove rigorously 
nonuniqueness of nonnegative solutions of  (\ref{S0E1})-(\ref{S0E3}).
More precisely, it is possible to obtain
different nonnegative solutions of this Fokker-Planck equation with a
different asymptotic behaviour as $\left( x,v\right) \rightarrow \left(
0,0\right) $. Actually, the solutions of (\ref{S0E1})-(\ref{S0E3}) that we
will construct in this paper will be measure-valued solutions and they will
differ in the amount of mass that they contain at the point $\left\{ \left(
0,0\right) \right\},$ i.e. $\int_{\left\{ \left( 0,0\right) \right\}
}P\left( dxdv,t\right).$

The intuitive explanation behind this nonuniqueness result is the following.
If $r<r_{c},$ the solutions of the problem (\ref{S1E1}), (\ref{S1E2}) reach
the point $\left( 0,0\right) $ in a finite time with probability one. After
reaching that point it is possible to give several dynamic laws for the
evolution of the particle $\left( X\left( t\right) ,V\left( t\right) \right)
.$ For instance, the particle could remain trapped at the origin, or it
could continue their movement, or it could remain trapped at $\left(
0,0\right) $ during some time and then restart its motion again. The
situation is essentially similar to the one of a Brownian particle moving in
a half-line $\left\{ x>0\right\} $ which reaches $x=0$ with probability one,
and there, it can be either absorbed, or to be reflected and continue its
motion back to the region $\left\{ x>0\right\}.$ All these possibilities
can be described in the case of a moving Brownian particle by means of
different boundary conditions at $x=0$ for the diffusion equation describing
the probability density associated to that process. In the case of the
Fokker-Planck equation (\ref{S0E1})-(\ref{S0E3}) the different evolution
laws assumed for $\left( X\left( t\right) ,V\left( t\right) \right) $ after
reaching the origin, will result in different boundary conditions for $%
P\left( x,v,t\right) $ as $\left( x,v\right) \rightarrow \left( 0,0\right).$

The nonuniqueness of the stochastic process associated to (\ref{S1E1}), (\ref%
{S1E2}) has been considered in the literature of stochastic processes. The
paper \cite{BE} concerns 
the construction of a stochastic process which
formally can be described by the stochastic differential equation $%
dX_{t}=V_{t}dt,\ dV_{t}=dW_{t},$ where $W_{t}$ is the Wiener process, and
where we impose in addition that $V_{t}=0$ if $X_{t}=0,$ $a.s.$ for $t\geq
0.$ Moreover, the process constructed in \cite{BE} satisfies that the
amount of time such that $X_{t}=0$ has zero measure $a.s.$ Notice that this
stochastic process can be considered as a solution of (\ref{S1E1}), (\ref%
{S1E2}) with $r=0,$ except for the fact that the solution is not requested
to remain in the half-plane $\left\{ X_{t}\geq 0\right\} .$ The surprising
feature of the stochastic process obtained in \cite{BE} is the fact that
"killing" the process at the times when $X_{t}$ reaches the line $\left\{
X_{t}=0\right\} $ does not force the particle to stay at this line for later
times. On the contrary, the particle is able to escape from the point $%
\left( X_{t},V_{t}\right) =\left( 0,0\right) $ at later times. The results
in \cite{BE} indicate that the solutions of (\ref{S1E1}), (\ref{S1E2}) can
be expected to be nonunique.

On the other hand, the problem (\ref{S1E1}), (\ref{S1E2}) has been
considered in \cite{J1}, \cite{J2}. It has been shown in \cite{J1} that for $%
r\geq r_{c}$ there is a unique entrance law in the half plane $\left\{ X>0,\
V\in \mathbb{R}\right\} $ under the assumption that $X_{0}=V_{0}=0.$ On the
other hand, if $r<r_{c}$ the bounces of the particle at the line $\left\{
X=0\right\} $ accumulate at the origin in a finite time. However, it is
proved in \cite{J2} that the particle does not remain necessarily at $\left(
X,V\right) =\left( 0,0\right) $ for later times, but on the contrary, it is
possible to define a "resurrected" process after the particle reaches the
origin. The measure of the times in which the solutions of this\
"resurrected " process remain in the line $\left\{ X=0\right\} $ is zero $%
a.s.$

In this paper we will obtain several results for the problem (\ref{S0E1})-(%
\ref{S0E3}) which can be viewed as the PDE reformulation of the previously
described results. More precisely, we will define a suitable concept of
measure valued solutions for (\ref{S0E1})-(\ref{S0E3}) if $r>r_{c}.$ On the
contrary it is possible to define different concepts of measure valued
solutions $P$ of (\ref{S0E1})-(\ref{S0E3}) if $r<r_{c}.$ The key point is
that for $r<r_{c}$ we must impose some boundary condition at the singular
point $\left( x,v\right) =\left( 0,0\right) $ in order to determine the
solution of (\ref{S0E1})-(\ref{S0E3}). Different boundary conditions can be
thought as to be associated to different physical meanings for the
corresponding evolution of the stochastic particle whose probability
distribution is given by $P.$ The boundary conditions obtained in this paper
will be denoted as trapping, nontrapping and partially trapping boundary
conditions at the point $\left( x,v\right) =\left( 0,0\right) $ and they can
be given respectively the meaning of a particle which after reaching the
point $\left( x,v\right) =\left( 0,0\right) $ remains trapped there for
later times, alternatively does not remain trapped at $\left( x,v\right)
=\left( 0,0\right) $ and continuous its motion in the half plane $\left\{
x>0\right\} $ or remains trapped during some characteristic time and
continues then its motion in the half-plane $\left\{ x>0\right\}.$ We will
not study in detail the case $r=r_{c}$ in this paper, because in that case
some logarithmic terms appear in the asymptotics of the solutions near the
singular point, something that makes the analysis rather cumbersome. 

The different solutions $P$ obtained will be nonnegative, classical
solutions of (\ref{S0E1})-(\ref{S0E3}) for $\left( x,v\right) \neq \left(
0,0\right) $ satisfying $\int_{\left\{ \left( x,v\right) \neq \left(
0,0\right) \right\} }P\left( x,v,t\right) dxdv<\infty ,$ but they differ in
their asymptotics as $\left( x,v\right) \rightarrow \left( 0,0\right).$ 
These different solutions will satisfy different definitions of
distributional solutions of (\ref{S0E1})-(\ref{S0E3}) which will
characterize the dynamics of probability densities $P\left( x,v,t\right) $
for which either the Brownian particles arriving to $\left( x,v\right)
=\left( 0,0\right) $ remain trapped there, or they continue their motion or
they remain trapped during a characteristic time before resuming their
motion in the region $\left\{ \left( x,v\right) \neq \left( 0,0\right) ,\
x>0\right\}.$ 
A more precise formulation of the results is contained in the
following Theorem. This Theorem will be proved rigorously in the companion paper \cite{HVJ2}.

\begin{theorem}
\label{ThIntro} Let  $\mathcal{X}$ be the space of functions $%
C^{\infty }\left( \mathbb{R}^{+}\times \mathbb{R}\times \left( 0,\infty
\right) \right) \cap C\left( \left[ 0,\infty \right) : \right.$ $\left.\mathcal{M}_{+}\left(
\mathbb{R}^{+}\times \mathbb{R}\right) \right) \cap C\left( \overline{%
\mathbb{R}^{+}\times \mathbb{R}}\times \left( 0,\infty \right) \right).$
Suppose that $0<r<r_{c}.$ For any Radon measure $P_{0}\in \mathcal{M}%
_{+}\left( \mathbb{R}^{+}\times \mathbb{R}\right) ,\ $such that $\int_{%
\mathbb{R}^{+}\times \mathbb{R}}P_{0}<\infty ,$ there exist infinitely many
different solutions of the problem (\ref{S0E1})-(\ref{S0E3}). The solutions $%
P=P\left( x,v,t\right) \in \mathcal{X}$ satisfy (\ref{S0E1}), (\ref{S0E2})
in classical sense and (\ref{S0E3}) in the sense of distributions. If $%
r_{c}<r\leq 1$ there exists a unique weak solution of (\ref{S0E1})-(\ref%
{S0E3}) with initial data $P_{0}.$
\end{theorem}

\begin{remark}
We have denoted as $\mathcal{M}_{+}\left( B\right) $ the set of Radon
measures in a given Borel set $B\subset \mathbb{R}^{n},\ n\geq 1.$
\end{remark}




In this paper we compute detailed asymptotic formulas for the solutions of (\ref{S0E1}), (\ref%
{S0E2}) near the point $(x,v)=(0,0)$ as well as its adjoint problem. As indicated above, it turns out that the measure $P$  might have an atom at the singular point $(x,v)=(0,0)$ and the asymptotic behaviour of the solutions is strongly dependent on the form of the interaction of the mass concentrated at the singular point with the mass in the surroundings.

The well-posedness of the Fokker-Planck equation for a general class of boundary conditions has been considered in \cite{Nier}. The approach used in that paper differs in several aspects from the one in this paper and the companion paper \cite{HVJ2}. The functional framework employed in \cite{Nier} uses extensively Sobolev spaces $H^k$ as well as Fourier analysis methods. The specific boundary conditions (\ref{S0E2}) are not considered in \cite{Nier} and neither the nonuniqueness result in Theorem \ref{ThIntro} above. From the functional analysis point of view the approach of this paper differs in a significant manner from the approach in \cite{Nier}. Our approach relies mostly in  the classical theory of Markov semigroups in the space of continuous functions. This allows us to obtain solutions of (\ref{S0E1}), (\ref%
{S0E2}) in the sense of measures, something that is convenient due to the interpretation of our problem in terms of the stochastic dynamics of the particle indicated earlier. 


We will use repeatedly the asymptotic symbol $\sim $ with the following
meaning. Suppose that $f_{1}\left( x,v\right) ,\ f_{2}\left( x,v\right) $
are two functions defined in some region $\mathcal{D}$ whose boundary
contains the origin $\left( x,v\right) =\left( 0,0\right) $. We will say
that $f_{1}\left( x,v\right) \sim f_{2}\left( x,v\right) $ as $\left(
x,v\right) \rightarrow \left( 0,0\right) $ if $\lim_{\left( x,v\right)
\rightarrow \left( 0,0\right) }\frac{f_{1}\left( x,v\right) }{f_{2}\left(
x,v\right) }=1$ with $\left( x,v\right) \in \mathcal{D}$. We will denote as
$|(x,v)| = x + |v|^3 $ the norm of $(x,v)$.

The plan of this paper is the following. In Section 2 we study a simple toy
model, namely the random walks of particles in a one-dimensional half
lattice. This yields a simple diffusion problem in the limit when the size
of the lattice converges to zero for which it is easy to obtain the boundary
conditions at the lower extreme of the lattice. The analysis of this simpler
problem will make more clear the meaning of the boundary conditions imposed
at the singular set of the solutions of (\ref{S0E1})-(\ref{S0E3}). Section 3 describes in
heuristic form the different expected asymptotic behaviours of the solutions
of (\ref{S0E1}), (\ref{S0E2}) near the singular point $\left( x,v\right)
=\left( 0,0\right) .$ These heuristic asymptotics are used to obtain, at the
formal level, a set of adjoint problems to (\ref{S0E1}), (\ref{S0E2}) in
Section 4. The onset of different adjoint problems is due to the possibility
of having different asymptotics for the solutions of (\ref{S0E1}), (\ref%
{S0E2}) if $r<r_{c}.$ 
In Section 5, we give a precise definition of measure valued solutions for the
problem (\ref{S0E1})-(\ref{S0E3}). 



%
%
%

\section{A toy problem: Diffusion in half-line.\label{ToyModel}}

In this section we discuss, using heuristic arguments, an elementary model
which exhibits some similarities with the behaviour of the model (\ref{S1E1}%
), (\ref{S1E2}) in the case $r<r_{c}.$ The results described in this section
are well known, since they just correspond to classical diffusion of a
particle which reaches the boundary of a half-line. However, the results of
this section will help to clarify the nonuniqueness of solutions for (\ref%
{S0E1})-(\ref{S0E3}) 
and the physical meaning of the
different types of solutions obtained for this problem.

Suppose that we consider the stochastic evolution of a particle in a lattice
of points contained in the half-line $\left\{ x\geq 0\right\} .$ If the
distance between the lattice points converges to zero, the probability
density of finding the particle at one specific position can be approximated
by means of the one-dimensional diffusion equation. However, different
boundary conditions must be assumed at $x=0$ for this probability density
depending on the type of dynamics prescribed for the particle which reach
that point.

More precisely, for each $h>0,$ we will denote as $\mathcal{L}_{h}$ the
lattice:%
\begin{equation*}
\mathcal{L}_{h}=\left\{ x_{n}=nh:n=0,1,2,3,...\right\}
\end{equation*}

For each $t_{k}=kh^{2}:k=0,1,2,...,$ we will denote as $X\left( t_{k}\right)
$ the position of one particle moving in the lattice $\mathcal{L}_{h}$ with
the following stochastic dynamics. If $X\left( t_{k}\right) =nh,$ with $n \geq 1$
we will assign the following probability laws for the position of the
particle at the time $t_{k+1}:$%
\begin{equation}
p\left( X\left( t_{k+1}\right) =\left( n+1\right) h~|~X\left( t_{k}\right)
=nh\right) =\frac{1}{2}\ \ ,\ \ p\left( X\left( t_{k+1}\right) =\left(
n-1\right) h~|~X\left( t_{k}\right) =nh\right) =\frac{1}{2}  \label{S2E1}
\end{equation}

On the other hand, we will define three different evolutions for the
particle arriving to the point $x_{0}=0.$ We will say that $X\left( \cdot
\right) $ evolves with trapping boundary conditions at $x=0$ if we assign
the following probability:%
\begin{equation}
p\left( X\left( t_{k+1}\right) =0~|~X\left( t_{k}\right) =0\right) =1
\label{S2E2}
\end{equation}

We will say that $X\left( \cdot \right) $ evolves with nontrapping boundary
conditions at $x=0$ if its evolution after reaching the point $x_{0}=0$ is
given by:%
\begin{equation}
p\left( X\left( t_{k+1}\right) =h~|~X\left( t_{k}\right) =0\right) =\lambda
\ \ ,\ \ p\left( X\left( t_{k+1}\right) =0~|~X\left( t_{k}\right) =0\right)
=1-\lambda  \label{S2E4}
\end{equation}%
where $\lambda \in \left( 0,1\right] $ independent of $h.$ Notice that
trapping boundary conditions are recovered from (\ref{S2E4}) if $\lambda =0$
.

Suppose that we write for short $P_{n}\left( k\right) =p\left( X\left(
t_{k}\right) =nh\right) .$ Then (\ref{S2E1}), (\ref{S2E4}) imply:%
\begin{equation}
P_{n}\left( k+1\right) =\frac{1}{2}\left( P_{n-1}\left( k\right)
+P_{n+1}\left( k\right) \right) \ ,\ \ \ n\geq 2  \label{S2E5}
\end{equation}%
\begin{equation}
P_{1}\left( k+1\right) =\frac{1}{2}P_{2}\left( k\right) +\lambda P_{0}\left(
k\right) \ ,\ \ \ P_{0}\left( k+1\right) =\frac{1}{2}P_{1}\left( k\right)
+\left( 1-\lambda \right) P_{0}\left( k\right)  \label{S2E6}
\end{equation}

The structure of these equations suggests that in time scales where $k$ is
large (i.e. time scales $t>>h^{2}$), the behaviour of the solutions of (\ref{S2E5}), (\ref{S2E6}) should be given locally
by steady states for all the values of $n$ if $\lambda >0$. Indeed, in this
case, particles can leave the point $x_{0}=0,$ and a local equilibrium can
be expected. The local equilibria satisfy:%
\begin{eqnarray}
\varphi _{n} &=&\frac{1}{2}\left( \varphi _{n-1}+\varphi _{n+1}\right) \ \
,\ \ \ n\geq 2  \label{phi_n2} \\
\varphi _{1} &=&\frac{1}{2}\varphi _{2}+\lambda \varphi _{0}\ \ ,\ \ \varphi
_{0}=\frac{1}{2}\varphi _{1}+\left( 1-\lambda \right) \varphi _{0}\ \ ,\ \
\lambda >0  \label{phi_n3}
\end{eqnarray}

The solutions of (\ref{phi_n2}) have the form $\varphi _{n}=A+Bn,$ $n\geq 1.$
Using (\ref{phi_n3}) it then follows that
\begin{equation}
B=0\   \label{S3E1}
\end{equation}

If we define the family of measures $u_{h}\left( x,t\right)
=\sum_{n,k}P_{n}\left( k\right) \delta _{x=x_{n}}\delta _{t=t_{k}}$ we would
expect to have $u_{h}\rightharpoonup U,$ where due to (\ref{S2E5}) it
follows that:%
\begin{equation}
\partial _{t}U=\frac{1}{2}\partial _{xx}U\ \ ,\ \ x>0\ \ ,\ \ t>0\
\label{S3E2}
\end{equation}%
and where (\ref{S3E1}) yields the following boundary condition for $U$ in
the case of nontrapping boundary conditions:%
\begin{equation}
\partial _{x}U\left( 0,t\right) =0  \label{S3E3}
\end{equation}

Suppose now that $\lambda =0.$ In this case, it does not make sense to
assume that $P_{0}\left( k\right) $ is close to an equilibrium value,
because any particle arriving to $x_{0}=0$ remains there for arbitrary times
and therefore $P_{0}\left( \cdot \right) $ is increasing. Solving (\ref%
{phi_n2}) we obtain $\varphi _{n}=A+Bn,$ $n\geq 1$ and using the first
equation in (\ref{phi_n3}) it then follows that $A=0.$ Therefore, in the
case of trapping boundary conditions we obtain that $U$ solves (\ref{S3E2})
with the boundary condition:%
\begin{equation}
U\left( 0,t\right) =0  \label{S3E4}
\end{equation}

Notice that the solutions of (\ref{S3E2}), (\ref{S3E4}) do not preserve $%
\int_{\left( 0,\infty \right) }U\left( x,t\right) dx.$ This seems
contradictory with the fact that $U\left( \cdot ,t\right) $ is a probability
density for all $t\geq 0.$ However, the paradox is solved just taking into
account that the mass lost from $\left( 0,\infty \right) $ is transferred to
the point $x=0.$ Therefore, the result to be expected in this case is:%
\begin{equation}
u_{h}\left( x,t\right) \rightharpoonup U\left( x,t\right) +m\left( t\right)
\delta _{0}\ \ ,\ \ m\left( t\right) =1-\int_{\left( 0,\infty \right)
}U\left( x,t\right) dx  \label{S3E4a}
\end{equation}

It is possible to obtain an intermediate limit choosing $\lambda $ in (\ref%
{S2E6}) depending on $h.$ We will assume that:%
\begin{equation}
\lambda =\mu h\   \label{S3E5a}
\end{equation}%
for some $\mu >0,$ independent of $h.$ The rationale behind this rescaling
is the following. Since we are interested in obtaining $U$ of order one in
regions where $x$ is of order one we can expect to have $P_{n}$ of order $h$
for $n\geq 1.$ On the other hand, we are interested in the case in which
there is a macroscopic fraction of mass at $x=0,$ whence $P_{0}$ is of order
one. The second equation of (\ref{S2E6}) yields then $\left( P_{0}\left(
k+1\right) -P_{0}\left( k\right) \right) =\frac{1}{2}P_{1}\left( k\right)
-\lambda P_{0}\left( k\right) .$ Given that the time scale between jumps is $%
h^{2}$ this suggests the approximation
\begin{equation}
\left[ \partial _{t}m\left( t\right) \right] h^{2}=\frac{1}{2}P_{1}\left(
k\right) -\lambda P_{0}\left( k\right)  \label{S3E5}
\end{equation}

Plugging this into the first equation of (\ref{S2E6}) we would obtain
\begin{equation*}
P_{1}\left( k+1\right) =\frac{1}{2}P_{2}\left( k\right) +\frac{1}{2}%
P_{1}\left( k\right) -\left[ \partial _{t}m\left( t\right) \right] h^{2}
\end{equation*}%
whence%
\begin{equation*}
P_{1}\left( k+1\right) -P_{1}\left( k\right) =\frac{1}{2}\left( P_{2}\left(
k\right) -P_{1}\left( k\right) \right) -\left[ \partial _{t}m\left( t\right) %
\right] h^{2}
\end{equation*}

Therefore, using that $P_{1},\ P_{2}$ are of order $h$, the definition of $%
u_{h}$ and the characteristic scales for the space and time jumps yield:%
\begin{equation}
\frac{1}{2}\partial _{x}U\left( 0,t\right) =\partial _{t}m\left( t\right)
\label{S3E7}
\end{equation}

This equation just provides the global mass conservation for the whole
system. On the other hand, (\ref{S3E5}) yields, if $h^{2}<<\lambda ,$ that $%
\frac{1}{2}P_{1}\left( k\right) \approx \lambda P_{0}\left( k\right) ,$
whence using the scaling (\ref{S3E5a}), as well as the fact that $P_{1}$ is
of order $h$ and $P_{0}$ is like $m\left( t\right) $:%
\begin{equation}
m\left( t\right) =\frac{U\left( 0,t\right) }{2\mu }  \label{S3E6}
\end{equation}

This boundary condition describes the equilibrium between the mass at $x=0$
and the density probability in the surroundings.\ It is possible to
reformulate it as a condition involving only local properties of the
function $U.$ Differentiating (\ref{S3E6}) with respect to $t$ and using (%
\ref{S3E7}) we obtain:%
\begin{equation}
U_{t}\left( 0,t\right) =\mu U_{x}\left( 0,t\right)  \label{S3E8}
\end{equation}

The problem (\ref{S3E2}), (\ref{S3E8}) describes the evolution of the
probability density if $\lambda $ rescales as in (\ref{S3E5a}) in the limit $%
h\rightarrow 0.$ The probability density, including the amount of the mass
at the origin is given by (\ref{S3E4a}). Notice that taking formally the
limit $\mu \rightarrow \infty $ in (\ref{S3E8}) we recover the boundary
condition (\ref{S3E3}). If we take $\mu =0$ in (\ref{S3E8}) and we assume
that $U\left( 0,0\right) =0,$ we recover the boundary condition (\ref{S3E4}).

The conclusion of these arguments is the following.\ If the underlying
particle system described by the parabolic equation (\ref{S3E2}) yields a\
dynamics for which the the particles can arrive to a boundary point of the
domain we need to complete the equation with a suitable boundary condition.
Different particle dynamics yield different boundary conditions. In
particular, trapping boundary conditions yield (\ref{S3E4}) and nontrapping
boundary conditions yield (\ref{S3E3}). Partially trapping boundary
conditions, with an escape rate from $x=0$ given by a coefficient $\lambda ,$
rescaling as (\ref{S3E5a}) yield the boundary condition (\ref{S3E8}).

A way of describing the different types of boundary conditions in an unified
way is the following. Suppose that we look for the following type of
asymptotics for the solutions of (\ref{S3E2}):%
\begin{equation}
U\left( x,t\right) \sim a\left( t\right) F\left( x\right) \ \ \text{as\ \ }%
x\rightarrow 0\   \label{Uorigin}
\end{equation}%
for a smooth function $a$ and a suitable $F$ which behaves like $x^a$ for some power $a$ as $%
x\rightarrow 0.$ Notice that the asymptotics (\ref{Uorigin}) does not make
any assumption about boundary conditions for (\ref{S3E2}). Plugging (\ref%
{Uorigin}) into (\ref{S3E2}) we would obtain that $\partial _{xx}F=0,$
whence we would have that, either $F$ is proportional to $1$ or to $x.$
Then, the asymptotics of $U$ near $x=0$ can be expected to be given by:%
\begin{equation}
U\left( x,t\right) \sim a_{0}\left( t\right) \cdot 1+a_{1}\left( t\right)
\cdot x+...\ \ \text{as\ \ }x\rightarrow 0  \label{Uzero}
\end{equation}

The classical theory of boundary value problems for parabolic equations can
be interpreted as stating that the meaningful boundary conditions for the
solutions of (\ref{S3E2}) are those imposing a relation between $a_{0}\left(
t\right) $ and $a_{1}\left( t\right) .$ For instance, the homogeneous
Dirichlet condition for (\ref{S3E2}) (cf. (\ref{S3E4})) is equivalent to
imposing $a_{0}\left( t\right) =0$ in (\ref{Uzero}), while the homogeneous
Neumann condition would mean $a_{1}\left( t\right) =0$. The boundary
condition (\ref{S3E8}) would yield $\frac{da_{0}}{dt}=\mu a_{1}.$

In the case of diffusion processes on the line, these results (\cite{C}, \cite{Fr2}) are well
known. It is relevant to notice that the boundary conditions (\ref{S3E3}), (%
\ref{S3E4}), (\ref{S3E8}) can be thought as different asymptotics for $U,$
near the boundary point under consideration. In the case of (\ref{S0E1})-(%
\ref{S0E3}), different dynamical laws for the particles reaching the
singular point $\left( x,v\right) =\left( 0,0\right) $ will result in
different asymptotic formulas for $P\left( x,v,t\right) $ as $\left(
x,v\right) \rightarrow \left( 0,0\right) .$ However, some
technical difficulties arise due to the fact that the boundary conditions
must be determined just at the point $\left( x,v\right) =\left( 0,0\right) ,$
although the problem is solved in the two-dimensional domain $\mathbb{R}%
^{+}\times \mathbb{R}$ for $t>0.$


%
%
%

\section{Heuristic description of the solutions of (\protect\ref{S0E1})-(%
\protect\ref{S0E3}) near the singular point.\label{Fsingular}}

\subsection{Asymptotics of the solutions of (\protect\ref{S0E1}), (\protect\ref%
{S0E2}) near the singular point.}

In this Subsection we compute formally the asymptotic behaviour of the
solutions of (\ref{S0E1}), (\ref{S0E2}) as $\left( x,v\right) \rightarrow
\left( 0,0\right) .$ This computation will illustrate the difference between
the cases $r<r_{c}$ and $r>r_{c}$ and it will provide a simple intuitive
explanation of the nonuniqueness arising for $r<r_{c}.$ The insight gained
in this Section will be useful later in order to define suitable concepts of
weak solutions of (\ref{S0E1})-(\ref{S0E3}).

Suppose that $P\left( x,v,t\right) $ is a solution of (\ref{S0E1}), (\ref%
{S0E2}) with the following asymptotic behaviour as $\left( x,v\right)
\rightarrow \left( 0,0\right) :$%
\begin{equation}
P\left( x,v,t\right) \sim a\left( t\right) G\left( x,v\right) \ \ ,\ \ \text{%
as\ \ }\left( x,v\right) \rightarrow \left( 0,0\right) \   \label{S4E1}
\end{equation}%
where $a$ is smooth and $G\left( x,v\right) $ has suitable homogeneity
properties to be determined. Plugging (\ref{S4E1}) into (\ref{S0E1}) we
would obtain the leading order terms
\begin{equation}
\left[ \partial _{t}a\left( t\right) \right] G\left( x,v\right) \sim a\left(
t\right) \left[ -v\partial _{x}G\left( x,v\right) +\partial _{vv}G\left(
x,v\right) \right]  \label{S4E2}
\end{equation}

Notice that, due to the presence of the derivatives with respect to $x$ and $%
v$ we can expect the terms on the left-hand side of (\ref{S4E2}) to be small
compared with the ones on the right as $\left( x,v\right) \rightarrow \left(
0,0\right) .$ It is then natural to assume, using also (\ref{S0E2}), that $%
G\left( x,v\right) $ satisfies:\
\begin{align}
v\partial _{x}G\left( x,v\right) & =\partial _{vv}G\left( x,v\right) \ \ ,\
\ x>0,\ \ v\in \mathbb{R},  \label{steady} \\
G\left( 0,v\right) & =r^{2}G\left( 0,-rv\right) \ \ ,\ \ v<0.
\label{steadyBC}
\end{align}%
This equation has been already found in the physics literature (cf. \cite%
{BK, CSB}).

The invariance of the equations (\ref{steady}), (\ref{steadyBC}) under the
rescaling $x\rightarrow \lambda ^{3}x,\ v\rightarrow \lambda v$ suggests to
look for solutions of these equations with the self-similar form:%
\begin{equation}
G\left( x,v\right) =x^{\gamma }\Phi \left( z\right) \ \ ,\ \ z=-\frac{v^{3}}{%
9x}  \label{Phi}
\end{equation}%
for some suitable $\gamma $ to be determined. Before we give the result on
the existence of such $G$, we start with the following lemma.

\begin{lemma}
\label{alpha(r)} There exists a unique function $\alpha \left( r\right) :$ $%
r\in \mathbb{R}_{+}\rightarrow \left( -\frac{5}{6},\frac{1}{6}\right) $
such that
\begin{equation}
(2+3\alpha \left( r\right) )\log r+\log \left( {2\cos \left( \pi \left(
\alpha \left( r\right) +\frac{1}{3}\right) \right) }\right) =0\;.
\label{exponent}
\end{equation}%
Let%
\begin{equation}
r_{c}=e^{-\frac{\pi }{\sqrt{3}}}  \label{rcrit}
\end{equation}%
The function $r\rightarrow \alpha \left( r\right) $ is non-decreasing in $r\in
\mathbb{R}_{+}.$ We have%
\begin{equation*}
\alpha \left( r_{c}\right) =-\frac{2}{3}\ \ ,\ \ \alpha (1)=0
\end{equation*}%
and $\alpha \left( r\right) \in \left( -\frac{5}{6},-\frac{2}{3}\right) $ if
$r<r_{c}$ and $\alpha \left( r\right) \in \left( -\frac{2}{3},\frac{1}{6}%
\right) $ if $r>r_{c}.$ Moreover,%
\begin{equation*}
\lim_{r\rightarrow 0^{+}}\alpha (r)=-\frac{5}{6}\ \ ,\ \ \lim_{r\rightarrow
\infty }\alpha (r)=\frac{1}{6}
\end{equation*}
\end{lemma}

\begin{proof}
For any given $r\in \mathbb{R}_{+}$, consider
\begin{equation*}
y_{r}\left( x\right) =(2+3x)\log r+\log \left( {2\cos \left( \pi \left( x+%
\frac{1}{3}\right) \right) }\right) .
\end{equation*}%
Then, the functions $y_{r}\left( \cdot \right) $ are smooth for $x\in (-%
\frac{5}{6},\frac{1}{6})$ and they satisfy $y_{r}\left( x\right) \rightarrow
-\infty $ as $x\rightarrow -\frac{5}{6},\frac{1}{6}$. Notice that $%
y_{r}\left( -\frac{2}{3}\right) =0$ for all $r$. Then, the function $\bar{%
\alpha}(r)\equiv -\frac{2}{3}$ satisfies (\ref{exponent}), but not the rest
of conditions imposed to $\alpha \left( r\right) $ in Lemma \ref{alpha(r)}.
In order to show that there exists another solution of (\ref{exponent}), we
first show that the equation $y_{r}\left( x\right) =0$ has a unique solution
$x\in (-\frac{5}{6},\frac{1}{6})$, $x\neq -\frac{2}{3}$ for $r\neq r_{c}$.
Note that
\begin{equation*}
\frac{dy_{r}}{dx}=3\log r-\pi \tan \left( \pi \left( x+\frac{1}{3}\right)
\right) .
\end{equation*}%
Hence $\frac{dy_{r}}{dx}=0$ has a unique solution $x_{rc}:=\frac{1}{\pi }%
\arctan \left( \frac{3}{\pi }\log r\right) -\frac{1}{3}$ and $\frac{dy_{r}}{%
dx}>0$ for $x<x_{rc}$ and $\frac{dy_{r}}{dx}<0$ for $x>x_{rc}$ and thus $%
y_{r}\left( \cdot \right) $ has the maximum value at $x=x_{rc}$. Note that $%
x_{rc}\rightarrow -\frac{5}{6}$ as $r\rightarrow 0^{+}$; $x_{rc}\rightarrow
\frac{1}{6}$ as $r\rightarrow \infty $; $x_{rc}=-\frac{2}{3}$ when $r=r_{c}$%
. Let $r<r_{c}$. Then $x_{rc}<-\frac{2}{3}$. Since $y_{r}(-\frac{2}{3})=0$,
there exists a unique $x_{r}$ so that $-\frac{5}{6}<x_{r}<x_{rc}<-\frac{2}{3}
$ and $y(x_{r})=0$. Similarly, if $r>r_{c}$, there exists a unique $x_{r}$
so that $-\frac{2}{3}<x_{rc}<x_{r}<\frac{1}{6}$ and $y_r(x_{r})=0$. Now if $%
r=r_{c}$, $x_{rc}=-\frac{2}{3}$ and hence the only solution of $y_{r}\left(
\cdot \right) =0$ is $x_{r}=-\frac{2}{3}$. Let us define $\alpha (r):=x_{r}$%
. Then it is now easy to deduce that $\alpha (r)$ satisfies all the
properties in the Lemma.
\end{proof}

We are now ready to state the result on the existence of $G$ of the form (%
\ref{Phi}) satisfying (\ref{steady}), (\ref{steadyBC}).

\begin{proposition}
\label{Fasympt} For any $r>0,$ $r\neq r_{c},$ there are two linearly
independent positive solutions of (\ref{steady}), (\ref{steadyBC}) which
take the form (\ref{Phi}). They are analytic in the domain $\left\{ \left(
x,v\right) :x>0,\ v\in \mathbb{R}\right\} $ and they have the form%
\begin{equation}
G_{\gamma }\left( x,v\right) =x^{\gamma }\Lambda _{\gamma }\left( \zeta
\right) \ ,\ \ \zeta =\frac{v}{\left( 9x\right) ^{\frac{1}{3}}}\ \ \text{%
with }\gamma \in \left\{ -\frac{2}{3},\alpha \left( r\right) \right\}
\label{Glambda}
\end{equation}%
with%
\begin{equation}
\Lambda _{\gamma }\left( \zeta \right) =U\left( -\gamma ,\frac{2}{3};-\zeta
^{3}\right) >0\ \ \text{for }\zeta \in \mathbb{R\ }  \label{LambdaU}
\end{equation}%
where we denote as $U(a,b;z)$ the classical Tricomi confluent hypergeometric
functions (cf. \cite{Ab}). Moreover, the asymptotic behaviour of $\Lambda
_{\gamma }\left( \zeta \right) ,\ \Lambda _{\gamma }^{\prime }\left( \zeta
\right) $ as $\left\vert \zeta \right\vert \rightarrow \infty $ is given by:%
\begin{eqnarray}
\Lambda _{\gamma }(\zeta ) &\sim &%
\begin{cases}
K_{\gamma }|\zeta |^{3\gamma },\quad \zeta \rightarrow \infty , \\
|\zeta |^{3\gamma },\quad \zeta \rightarrow -\infty .%
\end{cases}%
\ \   \label{asymptotics} \\
\Lambda _{\gamma }^{\prime }(\zeta ) &\sim &%
\begin{cases}
3\gamma K_{\gamma }|\zeta |^{3\gamma -1},\quad \zeta \rightarrow \infty , \\
-3\gamma |\zeta |^{3\gamma -1},\quad \zeta \rightarrow -\infty .%
\end{cases}%
\ \   \label{asyDer}
\end{eqnarray}%
where $\gamma \in \left( -\frac{5}{6},\frac{1}{6}\right) $ and $K_{\gamma
}=2\cos \left( \pi \left( \gamma +\frac{1}{3}\right) \right) .$
\end{proposition}

In order to prove Proposition \ref{Fasympt} we will use the following Lemma:

\begin{lemma}
\label{claim} For any $-5/6<\gamma <1/6$, we define:
\begin{equation}
\Lambda _{\gamma }(\zeta )=U(-\gamma ,\frac{2}{3};-\zeta ^{3}),\quad \zeta =%
\frac{v}{\left( 9x\right) ^{\frac{1}{3}}}\in \mathbb{R\ }\text{.}  \label{LU}
\end{equation}%
with $U(a,b;z)$ as in Proposition \ref{Fasympt}. Then:

\begin{itemize}
\item[(i)] $G_{\gamma }(x,v):=x^{\gamma }\Lambda _{\gamma }(\zeta )$
satisfies \eqref{steady}.

\item[(ii)] $\Lambda _{\gamma }(\zeta )$ is analytic in $\zeta \in \mathbb{C}
$ and $\Lambda _{\gamma }(\zeta )>0$ for any $\zeta \in \mathbb{R}$.

\item[(iii)] The asymptotic behaviour of $\Lambda _{\gamma }(\zeta )$ for
large $\left\vert \zeta \right\vert ,$ $\zeta \in \mathbb{R}$ is given by
the formulas in (\ref{asymptotics}).
\end{itemize}
\end{lemma}

\begin{proof}
The proof of (i) is just an elementary computation. We will show (ii) and
(iii) are valid. In order to study the properties of $\Phi (z)$ for negative
values of $z$ we use that (cf. \cite{Ab}, 13.1.3):%
\begin{equation}
U(a,b,z)=\frac{\pi }{\sin (\pi b)}\left( \frac{M(a,b,z)}{\Gamma
(1+a-b)\Gamma (b)}-z^{1-b}\frac{M(1+a-b,2-b,z)}{\Gamma (a)\Gamma (2-b)}%
\right) ,\quad b\notin \mathbb{Z}.  \label{UM}
\end{equation}

The function $M(a,b,z)$ is analytic for all $z\in \mathbb{C}$. Combining (%
\ref{LU}) and (\ref{UM}) we obtain the following representation formula for $%
\Lambda _{\gamma }(\zeta )$:%
\begin{equation}
\Lambda _{\gamma }(\zeta )=\frac{\pi }{\sin (\frac{2}{3}\pi )}\left( \frac{%
M(-\gamma ,\frac{2}{3},-\zeta ^{3})}{\Gamma (\frac{1}{3}-\gamma )\Gamma (%
\frac{2}{3})}+\zeta \frac{M(\frac{1}{3}-\gamma ,\frac{4}{3},-\zeta ^{3})}{%
\Gamma (-\gamma )\Gamma (\frac{4}{3})}\right) \ ,\ \ \zeta \in \mathbb{R}.
\label{Lamb}
\end{equation}%
Formula (\ref{Lamb}) provides a representation formula for $\Lambda (\zeta )$
in terms of the analytic functions $M(-\gamma ,\frac{2}{3},-\zeta ^{3}),$
$\ M(\frac{1}{3}-\gamma ,\frac{4}{3},-\zeta ^{3}).$ Therefore $\Lambda _{\gamma
}(\zeta )$ is analytic in $\zeta \in \mathbb{C}$.

We can compute the asymptotics of $\Lambda _{\gamma }(\zeta )$ as $\zeta
\rightarrow -\infty $ by using (\ref{LU}) and 13.5.2 in \cite{Ab}. Then we
deduce that
\begin{equation}
\Lambda _{\gamma }(\zeta )\sim \left\vert \zeta \right\vert ^{3\gamma }\ \
\text{as\ \ }\zeta \rightarrow -\infty .  \label{Aplus}
\end{equation}%
On the other hand, the formula 13.5.1 in \cite{Ab} yields the asymptotics:%
\begin{equation}
M\left( a,b,z\right) \sim \frac{\Gamma \left( b\right) e^{i\pi a}}{\Gamma
\left( b-a\right) }\left( z\right) ^{-a},\ \ \left\vert z\right\vert
\rightarrow \infty \ ,\ -\frac{\pi }{2}<\arg \left( z\right) <\frac{3\pi }{2}
\label{Masymp}
\end{equation}%
In particular, choosing $z=re^{i\pi }$ we obtain:%
\begin{equation}
M\left( a,b,-r\right) \sim \frac{\Gamma \left( b\right) }{\Gamma \left(
b-a\right) }\left( r\right) ^{-a},\ r\rightarrow \infty .  \label{Mneg}
\end{equation}

We remark that 13.5.1 in \cite{Ab} gives also the asymptotic formula:%
\begin{equation}
M\left( a,b,z\right) \sim \frac{\Gamma \left( b\right) e^{-i\pi a}}{\Gamma
\left( b-a\right) }\left( z\right) ^{-a},\ \ \left\vert z\right\vert
\rightarrow \infty \ ,\ -\frac{3\pi }{2}<\arg \left( z\right) <-\frac{\pi }{2%
}  \label{Masymp2}
\end{equation}

Notice that due to the analyticity of $M\left( a,b,z\right) $ in $\mathbb{C}$
the asymptotic behaviour of $M\left( a,b,z\right) $ obtained along
rectilinear paths approaching infinity and contained in $\left\{ \operatorname{Re}%
\left( z\right) <0\right\} $ must be the same independently on which formula
(\ref{Masymp}) or (\ref{Masymp2}) is used. In particular, it is easy to see
that (\ref{Mneg}) follows from (\ref{Masymp2}) choosing $z=re^{-i\pi }.$
Using (\ref{Mneg}) we obtain:%
\begin{eqnarray}
M(-\gamma ,\frac{2}{3},-\zeta ^{3}) &\sim &\frac{\Gamma \left( \frac{2}{3}%
\right) }{\Gamma \left( \frac{2}{3}+\gamma \right) }\left( \zeta \right)
^{3\gamma }\text{ as}\ \zeta \rightarrow \infty  \label{B1} \\
\zeta M(\frac{1}{3}-\gamma ,\frac{4}{3},-\zeta ^{3}) &\sim &\frac{\Gamma
\left( \frac{4}{3}\right) }{\Gamma \left( 1+\gamma \right) }\left( \zeta
\right) ^{3\gamma }\text{ as}\ \zeta \rightarrow \infty  \label{B2}
\end{eqnarray}

Combining (\ref{Lamb}), (\ref{B1}), (\ref{B2}) we obtain:%
\begin{equation*}
\Lambda _{\gamma }(\zeta )\sim \frac{\pi }{\sin (\frac{2}{3}\pi )}\left[
\frac{1}{\Gamma (\frac{1}{3}-\gamma )}\frac{1}{\Gamma \left( \frac{2}{3}%
+\gamma \right) }+\frac{1}{\Gamma (-\gamma )\Gamma \left( 1+\gamma \right) }%
\right] \left( \zeta \right) ^{3\gamma },\ \zeta \rightarrow \infty
\end{equation*}

Using then that $\Gamma \left( -x\right) \Gamma \left( 1+x\right) =-\frac{%
\pi }{\sin \left( \pi x\right) }$ (cf. 6.1.17\ in \cite{Ab}) as well as the
trigonometric formula $\frac{\left[ \sin \left( \pi \left( \gamma +\frac{2}{3%
}\right) \right) -\sin \left( \pi \gamma \right) \right] }{\sin (\frac{2}{3}%
\pi )}=2\cos \left( \pi \left( \gamma +\frac{1}{3}\right) \right) $ we
obtain
\begin{equation}
\Lambda _{\gamma }(\zeta )\sim K_{\gamma }\left\vert \zeta \right\vert
^{3\gamma }\ \ \text{as\ \ }\zeta \rightarrow \infty  \label{Aminus}
\end{equation}%
with $K_{\gamma }$ as in the statement of the\ Lemma. Notice that $K_{\gamma
}>0$ if $-5/6<\gamma <{1}/{6}.$

To finish the proof of Lemma \ref{claim}, it only remains to prove that $%
\Lambda _{\gamma }(\zeta )>0$ for any $\zeta \in \mathbb{R}$ and the
considered range of values of $\gamma .$ To this end, notice that if $\gamma
\rightarrow 0$ we have $\Lambda _{\gamma }(\zeta )\rightarrow 1>0$ uniformly
in compact sets of $\zeta .$ The functions\ $\Lambda _{\gamma }(\zeta )$
considered as functions of $\gamma ,$ change in a continuous manner. On the
other hand, the asymptotic behaviors (\ref{Aplus}), (\ref{Aminus}) imply
that the functions $\Lambda _{\gamma }(\zeta )$ are positive for large
values of $\left\vert \zeta \right\vert .$ If $\Lambda _{\gamma }(\cdot )$
has a zero at some $\zeta =\zeta _{0}\in \mathbb{R}$ and $-5/6<\gamma <{1}/{6%
}$, then there should exist, by continuity, $-5/6<\gamma _{\ast }<{1}/{6}$
and $\zeta _{\ast }\in \mathbb{R}$ such that $\Lambda _{\gamma _{\ast
}}(\zeta _{\ast })=\frac{d}{d\zeta }{\Lambda _{\gamma _{\ast
}}}(\zeta _{\ast })=0.$ The uniqueness theorem for ODEs then
implies that $\Lambda _{\gamma _{\ast }}(\cdot)=0,$ but
this would contradict the asymptotics (\ref{Aplus}), (\ref{Aminus}), whence
the result follows.
\end{proof}

\begin{proof}[Proof of Proposition \protect\ref{Fasympt}]
Assume $G$ takes the form (\ref{Phi}). Then $\Phi $ satisfies the following
ODE:
\begin{equation}
z\Phi _{zz}+(\frac{2}{3}-z)\Phi _{z}+\gamma \Phi =0.  \label{kummer}
\end{equation}%
It is also convenient to define the following auxiliary variable%
\begin{equation}
\zeta =\frac{v}{\left( 9x\right) ^{\frac{1}{3}}}  \label{zeta}
\end{equation}%
so that $z=-\zeta ^{3}$. Then $\Lambda _{\gamma }(\zeta )\equiv \Phi (-\zeta
^{3})$ satisfies the following ODE
\begin{equation}
\Lambda _{\gamma }^{\prime \prime }\left( \zeta \right) +3\zeta ^{2}\Lambda
_{\gamma }^{\prime }\left( \zeta \right) -9\gamma \zeta \Lambda \left( \zeta
\right) =0.  \label{LambdEq}
\end{equation}%
For each $\gamma ,$ two independent solutions to (\ref{kummer}) are given by
the Kummer function $M(-\gamma ,\frac{2}{3};z)$ and Tricomi confluent
hypergeometric function $U(-\gamma ,\frac{2}{3};z),$ (cf. \cite{Ab}). In
order to obtain a solution of (\ref{kummer}) which behaves algebraically as $%
z\rightarrow \pm \infty ,$ we recall the asymptotic behavior of $M(a,b;z)$
(see 13.1.4 and 13.1.5 in \cite{Ab}):
\begin{equation}
\begin{split}
M(a,b;z)& \sim \frac{\Gamma (b)}{\Gamma (a)}e^{z}z^{a-b}\quad \text{ for }%
z\rightarrow \infty , \\
M(a,b;z)& \sim \frac{\Gamma (b)}{\Gamma (b-a)}(-z)^{-a}\;\text{ for }%
z\rightarrow -\infty .
\end{split}
\label{asympt}
\end{equation}%
On the other hand, $U(-\gamma ,\frac{2}{3};z)$ behaves algebraically as $%
z\rightarrow \infty $ (see 13.5.2 in \cite{Ab}). Therefore, in order to get
the solutions satisfying the boundary condition \eqref{steadyBC}, $\Phi (z)$
should be proportional to the Tricomi confluent hypergeometric function $%
U(-\gamma ,\frac{2}{3};z)$. Due to the linearity of the problem we can
assume that the proportionality constant is one.

We now compute the behavior of $G_{\gamma }(x,v)=x^{\gamma }\Lambda _{\gamma
}(\zeta )$ near the boundary. For $v<0$ and $x\rightarrow 0^{+}$, since $%
\zeta \rightarrow -\infty $,
\begin{equation*}
G_{\gamma }(x,v)=x^{\gamma }\Lambda _{\gamma }(\zeta )\sim x^{\gamma
}\left\vert \frac{v}{(9x)^{\frac{1}{3}}}\right\vert ^{3\gamma }=\frac{%
|v|^{3\gamma }}{9^{\gamma }},\quad v<0\text{ and }x\rightarrow 0^{+}
\end{equation*}%
and for $v>0$ and $x\rightarrow 0^{+}$, since $\zeta \rightarrow \infty $,
\begin{equation*}
G_{\gamma }(x,v)=x^{\gamma }\Lambda _{\gamma }(\zeta )\sim x^{\gamma
}K_{\gamma }\left\vert \frac{v}{(9x)^{\frac{1}{3}}}\right\vert ^{3\gamma
}=K_{\gamma }\frac{|v|^{3\gamma }}{9^{\gamma }},\quad v>0\text{ and }%
x\rightarrow 0^{+}
\end{equation*}%
The boundary condition \eqref{steadyBC} implies that $\frac{|v|^{3\gamma }}{%
9^{\gamma }}=r^{2}K_{\gamma }\frac{|rv|^{3\gamma }}{9^{\gamma }}$ whence $%
r^{2+3\gamma }K_{\gamma }=1.$ Therefore $\gamma $ must satisfy the following
\begin{equation}
r^{2+3\gamma }\cdot 2\cos \left( \pi \left( \gamma +\frac{1}{3}\right)
\right) =1.  \label{Kr}
\end{equation}%
Notice that $\gamma =-\frac{2}{3}$ always fulfills the condition \eqref{Kr}
for all values of $r>0$ and therefore $G_{-\frac{2}{3}}(x,v)$ is a solution
of (\ref{steady}), (\ref{steadyBC}) for all $r>0$. The other value of $%
\gamma $ satisfying \eqref{Kr} is given by $\alpha (r)$ from Lemma \ref%
{alpha(r)} since taking the logarithm of (\ref{Kr}) yields (\ref{exponent}).
Therefore, $G_{\gamma }\left( x,v\right) =x^{\gamma }\Lambda _{\gamma
}\left( \zeta \right) ,\,\zeta =\frac{v}{\left( 9x\right) ^{\frac{1}{3}}}\ \
\text{with }\gamma \in \left\{ -\frac{2}{3},\alpha \left( r\right) \right\} $
are two linearly independent positive solutions of (\ref{steady}), (\ref%
{steadyBC}) for $r>0$ and $r\neq r_{c}$. 
\end{proof}

\begin{remark}
We note that both functions $G_{\gamma }\left( x,v\right) ,$ $\gamma \in
\left\{ -\frac{2}{3},\alpha \left( r\right) \right\} $ are integrable near
the origin, i.e. $\int_{\left\{ 0<x<1,\left\vert v\right\vert <1\right\}
}G_{\gamma }\left( x,v\right) dxdv<\infty .$ To this end we use the estimate
$G_{\gamma }\left( x,v\right) \leq Cx^{\gamma }$ if $\left\vert v\right\vert
\leq x^{1/3}$ and $G_{\gamma }\left( x,v\right) \leq C\left\vert
v\right\vert ^{3\gamma }$ if $\left\vert v\right\vert >x^{1/3}.$ Then, using
the fact that $\gamma >-6/5$ we obtain:%
\begin{eqnarray*}
\int_{\left\{ 0<x<1,\left\vert v\right\vert <1\right\} }G_{\gamma }\left(
x,v\right) dxdv &\leq &C\int_{\left\{ 0<x<1,\left\vert v\right\vert
<x^{1/3}\right\} }x^{\gamma }dxdv+C\int_{\left\{ 0<x<\left\vert v\right\vert
^{3},\left\vert v\right\vert <1\right\} }\left\vert v\right\vert ^{3\gamma
}dxdv \\
&\leq &C\int_{0}^{1}x^{\gamma +\frac{1}{3}}dx+C\int_{-1}^{1}\left\vert
v\right\vert ^{3\gamma +3}dv\leq C
\end{eqnarray*}
\end{remark}

\begin{remark}
We will denote from now on $\alpha \left( r\right) $ as $\alpha $ unless the
dependence on $r$ plays a role in the argument.
\end{remark}

\begin{remark}
\label{AsBeyond}The asymptotics (\ref{asymptotics}) is valid for arbitrary
values of $\gamma ,$ although $K_{\gamma }$ is not necessarily positive if $%
\gamma $ is not contained in the interval $\left( -\frac{5}{6},\frac{1}{6}%
\right) .$ We will occasionally use the asymptotics (\ref{asymptotics}) with
$\gamma =-\frac{2}{3}.$
\end{remark}

We next evaluate the particle fluxes towards the origin associated to $%
G_{\gamma }$ obtained in Proposition \ref{Fasympt}, which will be
importantly used to characterize the boundary conditions at the singular
point. We start with $G_{-\frac{2}{3}}(x,v)$.

We define a family of domains which will be repeatedly used in the following
arguments.

\begin{definition}
\label{domainsR}For any given $r>0$ and any $b>0,$ we define:%
\begin{equation}
\mathcal{R}_{\delta ,b}=\left\{ (x,v):0\leq x^{\frac{1}{3}}\leq b^{\frac{1}{3%
}}\delta ,\ -\delta \leq v\leq r\delta \right\} \   \label{domain}
\end{equation}%
We will denote as $\mathcal{R}_{\delta}$ the domain $\mathcal{R}_{\delta, 1
}.$
\end{definition}

\begin{proposition}
\label{LimitFluxes}Let $r>0$ be given and $\mathcal{R}_{\delta ,b}$ as in (%
\ref{domain}) for some $b>0$. Then
\begin{equation}
\int_{\partial \mathcal{R}_{\delta ,b}\cap \{x>0\}}\left[ -vG_{-\frac{2}{3}%
}n_{x}+\partial _{v}G_{-\frac{2}{3}}n_{v}\right] ds=9^{\frac{2}{3}}\left[
\log \left( r\right) +\frac{\pi }{\sqrt{3}}\right]  \label{Gflux}
\end{equation}%
where $n=(n_{x},n_{v})$ is the unit normal vector to $\partial \mathcal{R}%
_{\delta ,b}$ pointing towards $\mathcal{R}_{\delta ,b}.$
\end{proposition}

We will use the following result in the proof of Proposition \ref%
{LimitFluxes}.

\begin{lemma}
\bigskip \label{claim2} $\Lambda (\zeta )=\Lambda _{-\frac{2}{3}}(\zeta )$
satisfies the following ODE
\begin{equation}
\Lambda ^{\prime }\left( \zeta \right) +3\zeta ^{2}\Lambda \left( \zeta
\right) =3\   \label{lambdaP}
\end{equation}%
and moreover, $\Lambda (\zeta )$ is given by $\Lambda (\zeta
)=3\int_{-\infty }^{\zeta }\exp \left( -\zeta ^{3}+s^{3}\right) ds.$
\end{lemma}

\begin{proof}
Recall the equation for $\Lambda $ \eqref{LambdEq} when $\gamma =-2/3$:
\begin{equation*}
\Lambda ^{\prime \prime }+3\zeta ^{2}\Lambda ^{\prime }+6\zeta \Lambda =0.
\end{equation*}%
This equation is equivalent to $\left( \Lambda ^{\prime }+3\zeta ^{2}\Lambda
\right) ^{\prime }=0$ whence $\Lambda ^{\prime }\left( \zeta \right) +3\zeta
^{2}\Lambda \left( \zeta \right) =\Lambda ^{\prime }\left( 0\right) .$
Solving this equation we obtain:
\begin{equation}
\Lambda \left( \zeta \right) =C\exp \left( -\zeta ^{3}\right) +\Lambda
^{\prime }\left( 0\right) \int_{-\infty }^{\zeta }\exp \left( -\zeta
^{3}+s^{3}\right) ds\   \label{B3}
\end{equation}%
for some constant $C\in \mathbb{R}$. The function $\Lambda \left( \zeta
\right) $ in (\ref{B3}) increases exponentially as $\zeta \rightarrow
-\infty $ if $C\neq 0.$ It then follows from (\ref{asymptotics}) that $C=0,$
whence:%
\begin{equation*}
\Lambda \left( \zeta \right) =\Lambda ^{\prime }\left( 0\right)
\int_{-\infty }^{\zeta }\exp \left( -\zeta ^{3}+s^{3}\right) ds
\end{equation*}%
Then, the asymptotics of $\Lambda \left( \zeta \right) $ as $\zeta
\rightarrow -\infty $ is given by
\begin{equation*}
\Lambda \left( \zeta \right) \sim \frac{\Lambda ^{\prime }\left( 0\right) }{%
3\zeta ^{2}}\ \ \text{as\ \ }\zeta \rightarrow -\infty
\end{equation*}%
whence \eqref{asymptotics} yields $\Lambda ^{\prime }\left( 0\right) =3$ and
the result follows.
\end{proof}

We will use also the following result.

\begin{lemma}
\label{limM} Suppose that $\Lambda \left( \zeta \right) $ is as in Lemma \ref%
{claim2}. Then:%
\begin{equation}
\lim_{M\rightarrow \infty }\int_{-M}^{M}\zeta \Lambda \left( \zeta \right)
d\zeta =\frac{\pi }{\sqrt{3}}.  \label{limit}
\end{equation}
\end{lemma}

\begin{proof}
Using the representation formula for $\Lambda \left( \zeta \right) $
obtained in Lemma \ref{claim2} we obtain:
\begin{equation*}
\ell =\lim_{M\rightarrow \infty }\int_{-M}^{M}\zeta \Lambda \left( \zeta
\right) d\zeta =3\lim_{M\rightarrow \infty }\int_{-M}^{M}\left(
\int_{-\infty }^{\zeta }\exp \left( -\zeta ^{3}+s^{3}\right) ds\right) \zeta
d\zeta
\end{equation*}

Using the changes of variables $X=\zeta ^{3},\;Z=-s^{3}+\zeta ^{3}$ we
obtain:%
\begin{equation*}
\ell =\frac{1}{3}\lim_{M\rightarrow \infty }\int_{0}^{\infty }\exp \left(
-Z\right) dZ\int_{-M^{3}}^{M^{3}}\frac{dX}{X^{\frac{1}{3}}}\frac{1}{\left(
X-Z\right) ^{\frac{2}{3}}}
\end{equation*}%
where it is understood in the following that the function $X\rightarrow X^{%
\frac{1}{3}}$ is defined for $X\in \mathbb{R}$ and it takes negative values
for $X<0.$ In particular $\left( -X\right) ^{\frac{1}{3}}=-X^{\frac{1}{3}}$
for $X\in \mathbb{R}$.

Replacing now the integration in $X$ by integration in $t=\frac{X}{Z}$ we
arrive at:%
\begin{equation}
\ell =\frac{1}{3}\lim_{M\rightarrow \infty }\int_{0}^{\infty }\exp \left(
-Z\right) \varphi \left( \frac{Z}{M}\right) dZ=\frac{1}{3}\lim_{M\rightarrow
\infty }\int_{0}^{\infty }M\exp \left( -MZ\right) \varphi \left( Z\right) dZ
\label{B5}
\end{equation}%
with:%
\begin{equation}
\varphi \left( Z\right) =\int_{-\frac{1}{Z}}^{\frac{1}{Z}}\frac{dt}{t^{\frac{%
1}{3}}}\frac{1}{\left( t-1\right) ^{\frac{2}{3}}}  \label{varphi}
\end{equation}

Notice that if the limit $\varphi \left( 0^{+}\right) =\lim_{Z\rightarrow
0^{+}}\varphi \left( Z\right) $ exists, with $\varphi $ as in (\ref{varphi}%
), it would follow from (\ref{B5}) that $\ell =\frac{\varphi \left(
0^{+}\right) }{3}.$ We then compute $\varphi \left( 0^{+}\right) $ as
follows. We first split the integral in (\ref{varphi}) as $\varphi \left(
Z\right) =\int_{-\frac{1}{Z}}^{0}\left[ \cdot \cdot \cdot \right]
+\int_{0}^{1}\left[ \cdot \cdot \cdot \right] +\int_{1}^{\frac{1}{Z}}\left[
\cdot \cdot \cdot \right] .$ Using then the change of variables $t=-s+1$ in
the first integral and relabelling $s$ as $t$ we obtain:%
\begin{eqnarray}
\varphi \left( Z\right) &=&-\int_{1}^{\frac{1}{Z}+1}\frac{dt}{\left(
t-1\right) ^{\frac{1}{3}}}\frac{1}{\left( t\right) ^{\frac{2}{3}}}%
+\int_{0}^{1}\frac{dt}{t^{\frac{1}{3}}}\frac{1}{\left( t-1\right) ^{\frac{2}{%
3}}}+\int_{1}^{\frac{1}{Z}}\frac{dt}{t^{\frac{1}{3}}}\frac{1}{\left(
t-1\right) ^{\frac{2}{3}}}  \notag \\
&=&\int_{0}^{1}\frac{dt}{t^{\frac{1}{3}}}\frac{1}{\left( t-1\right) ^{\frac{2%
}{3}}}-\int_{\frac{1}{Z}}^{\frac{1}{Z}+1}\frac{dt}{\left( t-1\right) ^{\frac{%
1}{3}}}\frac{1}{\left( t\right) ^{\frac{2}{3}}}+  \label{varphiZ} \\
&&+\int_{1}^{\frac{1}{Z}}\left[ \frac{1}{t^{\frac{1}{3}}}\frac{1}{\left(
t-1\right) ^{\frac{2}{3}}}-\frac{1}{\left( t-1\right) ^{\frac{1}{3}}}\frac{1%
}{\left( t\right) ^{\frac{2}{3}}}\right] dt  \notag
\end{eqnarray}

The first integral on the right-hand side of (\ref{varphiZ}) can be computed
using Beta functions. The second one can be estimated by $-\log (1-Z)$ and
therefore it converges to zero as $Z\rightarrow 0^{+}.$ The third integral
on the right-hand side of (\ref{varphiZ}) converges to an integral in $%
\left( 1,\infty \right) $ as $Z\rightarrow 0^{+}.$ Then, using that $B\left(
\frac{2}{3},\frac{1}{3}\right) =\frac{2}{3}\sqrt{3}\pi $ we obtain:
\begin{equation*}
\varphi \left( 0^{+}\right) =\frac{2}{3}\sqrt{3}\pi +\int_{1}^{\infty }\left[
\frac{1}{t^{\frac{1}{3}}}\frac{1}{\left( t-1\right) ^{\frac{2}{3}}}-\frac{1}{%
\left( t-1\right) ^{\frac{1}{3}}}\frac{1}{\left( t\right) ^{\frac{2}{3}}}%
\right] dt
\end{equation*}

The integral on the right can be transformed to a more convenient form using
the change of variables $t=\frac{1}{x}.$ Then:%
\begin{equation}
\varphi (0^{+})=\frac{2\sqrt{3}\pi }{3}+I,  \label{varphi(0)}
\end{equation}%
where $I=\int_{0}^{1}\frac{dx}{x}\left[ \frac{1}{\left( 1-x\right) ^{\frac{2%
}{3}}}-\frac{1}{\left( 1-x\right) ^{\frac{1}{3}}}\right] .$ We can compute $%
I $ writing:
\begin{eqnarray*}
I &=&\lim_{\varepsilon \rightarrow 0^{+}}\int_{0}^{1}\frac{dx}{%
x^{1-\varepsilon }}\left[ \frac{1}{\left( 1-x\right) ^{\frac{2}{3}}}-\frac{1%
}{\left( 1-x\right) ^{\frac{1}{3}}}\right] \\
&=&\lim_{\varepsilon \rightarrow 0^{+}}\left[ B\left( \varepsilon ,\frac{1}{3%
}\right) -B\left( \varepsilon ,\frac{2}{3}\right) \right] =\lim_{\varepsilon
\rightarrow 0^{+}}\left( \Gamma \left( \varepsilon \right) \left( \frac{%
\Gamma \left( \frac{1}{3}\right) }{\Gamma \left( \frac{1}{3}+\varepsilon
\right) }-\frac{\Gamma \left( \frac{2}{3}\right) }{\Gamma \left( \frac{2}{3}%
+\varepsilon \right) }\right) \right) \\
&=&\left( \frac{\Gamma ^{\prime }\left( \frac{2}{3}\right) }{\Gamma \left(
\frac{2}{3}\right) }-\frac{\Gamma ^{\prime }\left( \frac{1}{3}\right) }{%
\Gamma \left( \frac{1}{3}\right) }\right) \lim_{\varepsilon \rightarrow
0^{+}}\left( \varepsilon \Gamma \left( \varepsilon \right) \right)
\end{eqnarray*}

Using then 6.1.3 and 6.3.7\ in \cite{Ab}, we obtain $I=\frac{\pi }{\sqrt{3}}%
. $ Plugging this into (\ref{varphi(0)}) we obtain $\varphi (0^{+})=\sqrt{3}%
\pi .$ It then follows from (\ref{B5}) that $\ell =\frac{\pi }{\sqrt{3}}$
and therefore, the Lemma follows.
\end{proof}

The previous Lemmas allow now to prove Proposition \ref{LimitFluxes}.

\begin{proof}[Proof of Proposition \protect\ref{LimitFluxes}]
We write
\begin{equation*}
Q\left( \delta ,b\right) =\int_{\partial \mathcal{R}_{\delta ,b}\cap \{x>0\}}%
\left[ -vG_{-\frac{2}{3}}n_{x}+\partial _{v}G_{-\frac{2}{3}}n_{v}\right] ds.
\end{equation*}%
The homogeneity of $G_{-\frac{2}{3}}$ and the definition of the domains $%
\mathcal{R}_{\delta ,b}$ implies that $Q\left( \delta ,b\right) =Q\left(
1,b\right) .$ On the other hand, Gauss Theorem and the differential equation
(\ref{steady}) yield that $Q\left( 1,b\right) $ is independent of $b.$
Therefore, we just need to show that
\begin{equation}
\lim_{b\rightarrow 0^{+}}Q\left( 1,b\right) =9^{\frac{2}{3}}\left[ \log
\left( r\right) +\frac{\pi }{\sqrt{3}}\right]  \label{limQ}
\end{equation}

Notice that the normal vector to $\partial \mathcal{R}_{1,b}$ is given by $%
n=(-1,0)$ if $x=b,\ -1\leq v\leq r$; $n=(0,-1)$ if $v=r,\ 0<x\leq b,$ $%
n=(0,1)$ if $v=-1,\ 0<x\leq b.$ Therefore,
\begin{eqnarray}
&&Q\left( 1,b\right)  \label{flux} \\
&=&-\int_{0}^{b}\partial _{v}G_{-\frac{2}{3}}\left( x,r\right)
dx+\int_{0}^{b}\partial _{v}G_{-\frac{2}{3}}\left( x,-1\right)
dx+\int_{-1}^{r}G_{-\frac{2}{3}}\left( b,v\right) vdv  \notag \\
&=&\left( I_{11}+I_{12}\right) +I_{2}=I_{1}+I_{2}  \notag
\end{eqnarray}

We will compute $I_{1}$: the first two integrals $I_{11}+I_{12}$. From %
\eqref{Glambda}, $\partial _{v}G_{-\frac{2}{3}}\left( x,v\right) =\frac{9^{-%
\frac{1}{3}}}{x}\Lambda ^{\prime }\left( \zeta \right) $ where $\Lambda
\equiv \Lambda _{-\frac{2}{3}}$, and thus
\begin{equation}
I_{1}=\frac{1}{9^{\frac{1}{3}}}\int_{0}^{b}\Lambda ^{\prime }\left( -\frac{1%
}{\left( 9x\right) ^{\frac{1}{3}}}\right) \frac{dx}{x}-\frac{1}{9^{\frac{1}{3%
}}}\int_{0}^{b}\Lambda ^{\prime }\left( \frac{r}{\left( 9x\right) ^{\frac{1}{%
3}}}\right) \frac{dx}{x}\   \label{B4}
\end{equation}

Using the change of variables $x\rightarrow r^{3}x$ in the last integral of (%
\ref{B4}) and splitting the resulting integral in the interval $\left( 0,%
\frac{b}{r^{3}}\right) $ in the integrals $\int_{0}^{b}\left[ \cdot \cdot
\cdot \right] +\int_{b}^{\frac{b}{r^{3}}}\left[ \cdot \cdot \cdot \right] $
we obtain:
\begin{eqnarray}
I_{1} &=&\frac{1}{9^{\frac{1}{3}}}\int_{0}^{b}\left[ \Lambda ^{\prime
}\left( -\frac{1}{\left( 9x\right) ^{\frac{1}{3}}}\right) -\Lambda ^{\prime
}\left( \frac{1}{\left( 9x\right) ^{\frac{1}{3}}}\right) \right] \frac{dx}{x}%
-\frac{1}{9^{\frac{1}{3}}}\int_{b}^{\frac{b}{r^{3}}}\Lambda ^{\prime }\left(
\frac{1}{\left( 9x\right) ^{\frac{1}{3}}}\right) \frac{dx}{x}  \notag \\
&=&I_{1a}+I_{1b}  \notag
\end{eqnarray}%
Notice that $I_{1b}$ depends on $r$ whereas $I_{1a}$ does not depend on $r$.
We first estimate $I_{1a}$. Since $\Lambda \left( \zeta \right) =O\left(
\left\vert \zeta \right\vert ^{-2}\right) $ as $\left\vert \zeta \right\vert
\rightarrow \infty ,$ and due to the analyticity properties of this function
from Proposition \ref{Fasympt}, we have also $\Lambda ^{\prime }\left( \zeta
\right) =O\left( \left\vert \zeta \right\vert ^{-3}\right) \ \text{as\ }%
\left\vert \zeta \right\vert \rightarrow \infty $. It then follows that $%
\Lambda ^{\prime }\left( \frac{1}{\left( 9x\right) ^{\frac{1}{3}}}\right) $,
$\Lambda ^{\prime }\left( -\frac{1}{\left( 9x\right) ^{\frac{1}{3}}}\right) $
are bounded for $0<x\leq b\leq 1,$ and
\begin{equation*}
\left\vert \Lambda ^{\prime }\left( \frac{1}{\left( 9x\right) ^{\frac{1}{3}}}%
\right) \right\vert +\left\vert \Lambda ^{\prime }\left( -\frac{1}{\left(
9x\right) ^{\frac{1}{3}}}\right) \right\vert \leq Cx\ \ ,\ \ 0<x\leq b\leq
1,\ \ 0<r\leq 1
\end{equation*}%
which yields
\begin{equation}
\left\vert I_{1a}\right\vert \leq Cb,  \label{Iab}
\end{equation}%
where $C$ is a uniform constant. Therefore
\begin{equation*}
\lim_{b\rightarrow 0}I_{1a}=0,~~\lim_{b\rightarrow 0}Q\left( 1,b\right)
=\lim_{b\rightarrow 0}\left( I_{1b}+I_{2}\right) .
\end{equation*}
In the rest of the Proof of this Proposition we compute this limit. Using (%
\ref{lambdaP}) we can rewrite $I_{1b}$ as:
\begin{equation}
\begin{split}
I_{1b}& =-\frac{3}{9^{\frac{1}{3}}}\int_{b}^{\frac{b}{r^{3}}}\frac{dx}{x}+%
\frac{3}{9^{\frac{1}{3}}}\int_{b}^{\frac{b}{r^{3}}}\frac{1}{\left( 9x\right)
^{\frac{2}{3}}}\Lambda \left( \frac{1}{\left( 9x\right) ^{\frac{1}{3}}}%
\right) \frac{dx}{x} \\
& =9^{\frac{2}{3}}\log (r)+9^{\frac{2}{3}}\int_{r\left( 9b\right) ^{-\frac{1%
}{3}}}^{\left( 9b\right) ^{-\frac{1}{3}}}\Lambda \left( \zeta \right) \zeta
d\zeta ,
\end{split}
\label{Ic}
\end{equation}%
where we have changed the variables to $\zeta =(9x)^{-\frac{1}{3}}$ for the
second integral. On the other hand, using (\ref{Glambda}) we can write $%
I_{2} $ in (\ref{flux}) as $I_{2}=\int_{-1}^{r}\left( b\right) ^{-\frac{2}{3}%
}\Lambda \left( \frac{v}{\left( 9b\right) ^{\frac{1}{3}}}\right) vdv=9^{%
\frac{2}{3}}\int_{-\left( 9b\right) \ \ ^{-\frac{1}{3}}}^{r\left( 9b\right)
\ \ ^{-\frac{1}{3}}}\Lambda \left( \zeta \right) \zeta d\zeta .$ Then:%
\begin{equation*}
\lim_{b\rightarrow 0}\left( I_{1b}+I_{2}\right) =9^{\frac{2}{3}}\log (r)+9^{%
\frac{2}{3}}\lim_{b\rightarrow 0}\int_{-\left( 9b\right) \ ^{-\frac{1}{3}%
}}^{r\left( 9b\right) \ ^{-\frac{1}{3}}}\Lambda \left( \zeta \right) \zeta
d\zeta
\end{equation*}%
and using Lemma \ref{limM} we obtain (\ref{limQ}) and the Proposition
follows.
\end{proof}

\begin{proposition}
\label{FluxesGalpha}Let $r>0$ be given, $r\neq r_{c}$ with $r_{c}$ as in (%
\ref{rcrit})$.$ Let $\mathcal{R}_{\delta ,b}$ be as in (\ref{domain}) with $%
b>0.$ Suppose that $G_{\alpha }$ is as in (\ref{Glambda}) with $\gamma
=\alpha =\alpha \left( r\right) ,$ where $\alpha \left( r\right) $ is as in
Lemma \ref{alpha(r)}. Then:%
\begin{equation}
\lim_{\delta \rightarrow 0}\int_{\partial \mathcal{R}_{\delta ,b}\cap
\{x>0\}}\left[ -vG_{\alpha }n_{x}+\partial _{v}G_{\alpha }n_{v}\right] ds=0
\label{FGalp}
\end{equation}%
where $n=(n_{x},n_{v})$ is the unit normal vector to $\partial \mathcal{R}%
_{\delta ,b}$ pointing towards $\mathcal{R}_{\delta ,b}.$
\end{proposition}

\begin{proof}
We write $Q\left( \delta ,b\right) =\int_{\partial \mathcal{R}_{\delta
,b}\cap \{x>0\}}\left[ -vG_{\alpha }n_{x}+\partial _{v}G_{\alpha }n_{v}%
\right] ds.$ Using Gauss Theorem we obtain that the function \thinspace $%
Q\left( \delta ,b\right) $ is independent of $b,$ whence $Q\left( \delta
,b\right) =Q\left( \delta ,1\right) .$ On the other hand, using Gauss
Theorem, as well as the boundary condition (\ref{steadyBC}) we obtain that $%
Q\left( \delta ,1\right) $ is independent of $\delta .$ Moreover, the
homogeneity of $G_{\alpha }$ and $\mathcal{R}_{\delta ,b}$ imply that $%
Q\left( \delta ,1\right) =C\delta ^{2+3\alpha },$ for some suitable constant
$C\in \mathbb{R}$. The independence of $Q\left( \delta ,1\right) $ of $%
\delta $ then implies $C=0$ and the result follows.
\end{proof}

\subsection{The case $r<r_{c}:$ Trapping, nontrapping and partially trapping
boundary conditions.\label{subcritical_bc}}

The main heuristic idea behind the nonuniqueness results in this paper for $%
r<r_{c}$ as well as the role of the critical parameter $r_{c}$ can be seen
as follows. Proposition \ref{Fasympt} suggests that an integrable
nonnegative solution of (\ref{S0E1}), (\ref{S0E2}) in $\left( x,v\right) \in
\mathbb{R}_{+}\times \mathbb{R}$ has the following asymptotic behaviour (cf.
(\ref{S4E1})):%
\begin{equation}
P\left( x,v,t\right) \sim a_{\alpha }\left( t\right) G_{\alpha }\left(
x,v\right) +a_{-\frac{2}{3}}\left( t\right) G_{-\frac{2}{3}}\left(
x,v\right) \ \ \text{as\ \ }\left( x,v\right) \rightarrow \left( 0,0\right)
\label{Forigin}
\end{equation}%
for suitable functions $a_{-\frac{2}{3}}\left( t\right) ,\ a_{\alpha }\left(
t\right) .$ Notice that for $r<r_{c}$ the most singular term in the
right-hand side of (\ref{Forigin}) is $a_{\alpha }\left( t\right) G_{\alpha
}\left( x,v\right) $ (assuming that $a_{\alpha }\left( t\right) \neq 0$). If
we assume also that the functions $a_{\alpha }\left( t\right) ,\ a_{-\frac{2%
}{3}}\left( t\right) $ are differentiable, we can expect to have corrective
terms in (\ref{Forigin}) of order $\left( x+\left\vert v\right\vert
^{3}\right) ^{\beta +\frac{2}{3}}$ with $\beta =\min \left\{ \alpha ,-\frac{%
2}{3}\right\} .$ Given that $\beta +\frac{2}{3}>\max \left\{ \alpha ,-\frac{2%
}{3}\right\} $ it then follows that such corrective terms would be
negligible compared with the two terms on the right-hand side of (\ref%
{Forigin}).

By analogy with the one-dimensional diffusion process considered in Section %
\ref{ToyModel} we can expect to be able to impose boundary conditions for $P$
at the point $\left( x,v\right) =\left( 0,0\right) $ by means of a
relationship between $a_{-\frac{2}{3}}\left( t\right) $ and $a_{\alpha
}\left( t\right) .$ In order to understand the meaning of those conditions
we remark that the solution $a_{-\frac{2}{3}}\left( t\right) G_{-\frac{2}{3}%
}\left( x,v\right) $ is associated to particle fluxes towards $\left(
x,v\right) =\left( 0,0\right) $ in the same manner as the term $a_{1}\left(
t\right) \cdot x$ in (\ref{Uzero}) is associated to fluxes towards $x=0$ for
the solutions of (\ref{S3E2}). This can be seen by means of the following
computation. Suppose that $P$ has the asymptotics (\ref{Forigin}) and it
decreases fast enough as $\left\vert \left( x,v\right) \right\vert
\rightarrow \infty ,$ in order to avoid particle fluxes towards infinity.
Let us denote as $\mathcal{U}$ the set:\
\begin{equation}
\mathcal{U}=\left\{ \left( x,v\right) :x\geq 0,\ v\in \mathbb{R},\ \left(
x,v\right) \neq \left( 0,0\right) \right\} .  \label{S8E6}
\end{equation}

We will use also the following notation for the boundary of $\mathcal{U}$:%
\begin{equation}
\partial \mathcal{U}=\left\{ \left( 0,v\right) :\ v\in \mathbb{R}\right\}  \label{S8E6a}
\end{equation}

We compute $\partial _{t}\left( \int_{\mathcal{U}}Pdxdv\right) $ using (\ref%
{S0E1}), (\ref{S0E2}):%
\begin{equation}
\partial _{t}\left( \int_{\mathcal{U}}Pdxdv\right) =\int_{\mathcal{U}}\left[
-v\partial _{x}P+\partial _{vv}P\right] dxdv=\lim_{\delta \rightarrow
0}\int_{\mathcal{U}\diagdown \mathcal{R}_{\delta }}\left[ -v\partial
_{x}P+\partial _{vv}P\right] dxdv  \label{S5E1}
\end{equation}%
where $\mathcal{R}_{\delta }$ is as in Definition \ref{domainsR}.

We can transform the integral on the right-hand side of (\ref{S5E1}) using
Gauss Theorem. Then, the right-hand side of (\ref{S5E1}) can be transformed
in:

\begin{equation}
\lim_{\delta \rightarrow 0}\int_{\left\{ v<-\delta ^{3}\text{ or}\ v>r\delta
^{3}\right\} }vP\left( 0,v,t\right) dv+\lim_{\delta \rightarrow
0}\int_{\partial \mathcal{R}_{\delta }\cap \left\{ x>0\right\} }\left[
-vPn_{x}+\partial _{v}Pn_{v}\right] ds\   \label{S5E2}
\end{equation}%
where $n=\left( n_{x},n_{v}\right) $ is the normal vector to $\partial
\mathcal{R}_{\delta }$ away from $\mathcal{R}_{\delta }$. The first integral
in (\ref{S5E2}) vanishes due to (\ref{S0E2}).\ Using the asymptotics (\ref%
{Forigin}) we can then write the left-hand side of (\ref{S5E1}) as:%
\begin{eqnarray}
&&a_{-\frac{2}{3}}\left( t\right) \lim_{\delta \rightarrow 0}\int_{\partial
\mathcal{R}_{\delta }\cap \left\{ x>0\right\} }\left[ -vG_{-\frac{2}{3}%
}n_{x}+\partial _{v}G_{-\frac{2}{3}}n_{v}\right] ds  \label{KeyIntegrals} \\
&&+a_{\alpha }\left( t\right) \lim_{\delta \rightarrow 0}\int_{\partial
\mathcal{R}_{\delta }\cap \left\{ x>0\right\} }\left[ -vG_{\alpha
}n_{x}+\partial _{v}G_{\alpha }n_{v}\right] ds  \notag
\end{eqnarray}

Using Propositions \ref{LimitFluxes}, \ref{FluxesGalpha} we can compute the
limits in (\ref{KeyIntegrals}) whence:%
\begin{equation}
\partial _{t}\left( \int_{\mathcal{U}}Pdxdv\right) =-\kappa a_{-\frac{2}{3}%
}\left( t\right) \ \text{with }\kappa =-9^{\frac{2}{3}}\left[ \log \left(
r\right) +\frac{\pi }{\sqrt{3}}\right] .  \label{S8E1}
\end{equation}

We will now indicate how to define different types of boundary conditions
for $P,$ assuming that we have the asymptotics (\ref{Forigin}). Taking into
account (\ref{S8E1}) it is natural to assume in the case of nontrapping
boundary conditions:%
\begin{equation}
a_{-\frac{2}{3}}\left( t\right) =0  \label{S5E3}
\end{equation}

Notice that in the asymptotics (\ref{Forigin}) the most singular term is $%
a_{\alpha }\left( t\right) G_{\alpha }\left( x,v\right) .$ Since $G_{\alpha
}\left( x,v\right) >0$ (cf. Proposition \ref{Fasympt}), and $P\geq 0$ we
must have $a_{\alpha }\left( t\right) \geq 0.$

On the other hand, notice that in the Toy model considered in Section 2 
the probability density $P$ takes smallest values than in the other two cases
if we impose trapping boundary conditions, since in that case the boundary condition $P(0)=0$ holds. 
Therefore $P$ is smaller than in the case of reflecting or mixed boundary conditions 
due to the maximum principle. It is natural, arguing by analogy, 
to define then trapping boundary condition by means of:

\begin{equation}
a_{\alpha }\left( t\right) =0  \label{S5E4}
\end{equation}

In principle, in order to show that the condition (\ref{S5E4}) is the one
associated to trapping boundary conditions, one should study in detail the
properties of a stochastic process in which particles can arrive to $\left(
x,v\right) =\left( 0,0\right) $ in finite time and to impose that those
particles remain there for later times. Alternatively, some heuristic
justification of (\ref{S5E4}) by means of a discrete particle model, in the
spirit of the Toy model considered in Section \ref{ToyModel} could be given.
We will not do neither of them in this paper. However, it is possible to
provide some justification by means of PDE arguments of the fact that (\ref%
{S5E4}) is the condition that must be imposed in order to obtain a
particle density $P\left( x,v,t\right) $ with trapping boundary
conditions at $\left( x,v\right) =\left( 0,0\right) .$ Indeed, notice that (%
\ref{Forigin}) and (\ref{S5E4}) as well as the fact that $P\left(
x,v,t\right) \geq 0$ imply that $a_{-\frac{2}{3}}\left( t\right) \geq 0.$
Using then (\ref{S8E1}) as well as the fact that $r<r_{c}$ it follows that
for any nonnegative solution of (\ref{S0E1}), (\ref{S0E2}) satisfying (\ref%
{Forigin}), (\ref{S5E4}) we have $\partial _{t}\left( \int_{\mathcal{U}%
}Pdxdv\right) \leq 0,$ i.e. for these solutions the mass could be
transferred from $\mathcal{U}$ to $\left( x,v\right) =\left( 0,0\right) ,$
but not in the reverse way, as it would be expected for trapping boundary
conditions. Moreover, for the class of solutions satisfying (\ref{Forigin}),
the boundary condition (\ref{S5E4}) is necessary in order to have $\partial
_{t}\left( \int_{\mathcal{U}}Pdxdv\right) \leq 0.$ Indeed, suppose that a
positive solution of (\ref{S0E1}), (\ref{S0E2}) satisfies (\ref{Forigin})
with $a_{\alpha }\left( t\right) >0,\ a_{-\frac{2}{3}}\left( t\right) <0$
for some time interval$.$ For such distributions we would have, due to (\ref%
{S8E1}) that $\partial _{t}\left( \int_{\mathcal{U}}Pdxdv\right) >0.$
Therefore, (\ref{S5E4}) is the boundary condition that must be imposed for
trapping boundary conditions if we assume that the solutions of (\ref{S0E1}%
), (\ref{S0E2}) satisfy (\ref{Forigin}).

We now derive the boundary condition playing a role 
similar to the mixed boundary condition obtained for the Toy model in Section 2.
For these boundary conditions, the flux of mass from $(x,v)=(0,0)$ 
to the region $\mathcal(U)$ should be proportional to the amount of mass at the point $(x,v)=(0,0).$
Given that the fluxes take place in a very small region close to the origin, 
we can expect this region to be in local equilibrium, including in this equilibrium the
transfer of mass from $(x,v)=(0,0)$ towards $\mathcal(U)$. We notice that $P$ is approximately given to the leading order in the region $\mathcal(U)$ in a neighbourhood of $(0,0)$ by $a_{\alpha}(t) G_{\alpha}(x,v). $
Therefore, the local equilibrium condition would require a condition of the form 

\begin{equation}
a_{\alpha }\left( t\right) =\mu _{\ast }m\left( t\right)  \label{S5E6}
\end{equation}%
for some $\mu _{\ast }\geq 0.$ The nonnegativity of $\mu _{\ast }$ is due to
the fact that $a_{\alpha }\left( t\right) \geq 0,$ $m\left( t\right) \geq 0.$
Notice that in the case $\mu _{\ast }=0$ (\ref{S5E6}) reduces to the
trapping boundary condition (\ref{S5E4}).

Differentiating (\ref{S5E6}), and using that, due to (\ref{S8E1}) we have
\begin{equation*}
\frac{dm}{dt}=-9^{\frac{2}{3}}\left[ \log \left( r\right) +\frac{\pi }{\sqrt{%
3}}\right] a_{-\frac{2}{3}}\left( t\right) ,
\end{equation*}
and we obtain the following boundary condition%
\begin{equation}
\frac{da_{\alpha }}{dt}=-9^{\frac{2}{3}}\left[ \log \left( r\right) +\frac{%
\pi }{\sqrt{3}}\right] \mu _{\ast }a_{-\frac{2}{3}}  \label{S5E7}
\end{equation}

We can interpret the solutions of (\ref{S0E1}%
), (\ref{S0E2}), (\ref{S5E7})
as the particle density associated to a particle system with the property
that a particle reaching the point $\left( x,v\right) =\left(
0,0\right)$ returns to $\mathcal{U}$ to the rate $\mu
_{\ast }.$.
Notice that in the limit case $\mu _{\ast }=\infty $ (\ref{S5E7})
formally reduces to the nontrapping boundary condition (\ref{S5E3}).

\subsection{The case $r>r_{c}:$ nontrapping boundary conditions and particle
fluxes from $\left( 0,0\right) $ to $\mathcal{U}$.\label{supercritical_bc}}

In the case $r>r_{c}$ (including $r\geq 1$), the situation changes
completely with respect to the case $r<r_{c}$ due to the fact that we now
have $\alpha =\alpha \left( r\right) >-\frac{2}{3},$ and also because the
right-hand side of (\ref{Gflux}) becomes now nonnegative.

Suppose that the asymptotics (\ref{Forigin}) holds. Then, to the leading
order $P$ can be approximated near $\left( x,v\right) =\left( 0,0\right) $
by means of $a_{-\frac{2}{3}}\left( t\right) G_{-\frac{2}{3}}\left(
x,v\right) .$ The nonnegativity of $P$ implies that $a_{-\frac{2}{3}}\left(
t\right) \geq 0,$ and therefore (\ref{S8E1}) yields $\partial _{t}\left(
\int_{\mathcal{U}}Pdxdv\right) \geq 0.$ Moreover, if $a_{-\frac{2}{3}}\left(
t\right) >0$ we would obtain $\partial _{t}\left( \int_{\mathcal{U}%
}Pdxdv\right) >0.$ Therefore, if $r>r_{c}$ the only boundary condition which
is compatible with solutions defining a probability measure in $\mathcal{U}%
\cup \left\{ \left( 0,0\right) \right\} $ for general initial data is the
nontrapping boundary condition, namely:

\begin{equation}
a_{-\frac{2}{3}}\left( t\right) =0  \label{T1E1}
\end{equation}

Notice that the previous argument does not imply that it is impossible to
construct solutions of the PDE problem (\ref{S0E1}), (\ref{S0E2}) satisfying
(\ref{Forigin}) and having $a_{-\frac{2}{3}}\left( t\right) \neq 0.$ The
problem is that those solutions cannot be understood in general as
probability distributions. Indeed, if $a_{-\frac{2}{3}}\left( t\right) <0$
we would have $P\left( x,v,t\right) <0$ for $\left( x,v\right) \in \mathcal{U%
}$ small, and then, the resulting function $P$ would not be a probability
density. On the other hand it is possible to have positive solutions of (\ref%
{S0E1}), (\ref{S0E2}) satisfying (\ref{Forigin}) with $a_{-\frac{2}{3}%
}\left( t\right) >0.$ However, for such solutions $\int_{\mathcal{U}}Pdxdv$
is an increasing function of $t$. It would be possible to obtain conserved
measures $m\delta _{\left( 0,0\right) }+P$ assuming that initially $m>0$ and
$m+\int_{\mathcal{U}}Pdxdv=1,$ and having $a_{-\frac{2}{3}}\left( t\right)
>0 $ only during the range of times in which $\int_{\mathcal{U}}Pdxdv<1.$
Since those boundary conditions do not allow to interpret $P$ as a
particle density they will not be considered in this paper, although they
could be useful in some problems. We could also have $P>0,\ a_{-\frac{2}{3}%
}\left( t\right) >0$ if we do not impose that $P$ is the restriction of a
probability measure to $\mathcal{U}$ but, say, a particle density. The
corresponding solutions would represent then, particle densities for which
there is a flux of particles from $\left( x,v\right) =\left( 0,0\right) $ to
$\mathcal{U}$. It is interesting to remark that the same fluxes from $\left(
x,v\right) =\left( 0,0\right) $ to $\mathcal{U}$ can be obtained if $%
r<r_{c}, $ although in that case, we must assume that $a_{-\frac{2}{3}%
}\left( t\right) <0$ in order to obtain an increasing number of particles in
$\mathcal{U}$. Nevertheless, we will restrict our attention in the following
just to the solutions of (\ref{S0E1}), (\ref{S0E2}) satisfying (\ref{Forigin}%
) as well as one of the boundary conditions (\ref{S5E3}), (\ref{S5E4}), (\ref%
{S5E7}) if $r<r_{c}$ and (\ref{T1E1}) if $r>r_{c}.$ Indeed, for those
boundary conditions we have natural interpretations for $P$ as the
particle density describing the evolution of a particle with nontrapping,
trapping or partially trapping boundary conditions.

\subsection{Definition of a Probability Measure in $\mathcal{V}=\overline{%
\mathcal{U}}$\label{defMeasure}}

Suppose that $P$ is a solution of (\ref{S0E1})-(\ref{S0E3}) in $\mathcal{U}$
satisfying (\ref{Forigin}) as well as one of the boundary conditions (\ref%
{S5E4}), (\ref{S5E7}), (\ref{S5E3}) (if $r<r_{c}$) or (\ref{T1E1}) (if $%
r>r_{c}$). In order to have mass conservation, it is convenient to define a
probability measure $f$ in the domain $\mathcal{V}=\overline{\mathcal{U}}$
as follows:%
\begin{equation}
f\left( x,v,t\right) =P\left( x,v,t\right) +m\left( t\right) \delta _{\left(
0,0\right) }\left( x,v\right) \   \label{meas1}
\end{equation}%
where:%
\begin{equation}
\int_{\mathcal{U}}P\left( \cdot ,t\right) dxdv+m\left( t\right) =1\
\label{meas2}
\end{equation}%
for $t\geq 0.$

Notice that due to (\ref{meas2}) the measure $f$ is a probability measure.
In the case $r>r_{c},$ where the boundary condition which we need to impose
is (\ref{T1E1}), or in the case $r<r_{c},$ if we assume the boundary
condition (\ref{S5E3}) we have that $\int_{\mathcal{U}}P\left( \cdot
,t\right) dxdv$ is constant. If we assume that $\int_{\mathcal{U}%
}P_{0}\left( \cdot \right) dxdv=1,$ it then follows that in those cases $%
m\left( t\right) =0,\ t\geq 0.$

Notice that in order to define the measure $f$ we do not need to have (\ref%
{Forigin}). Since we will use later functions $P$ for which the detailed
asymptotics (\ref{Forigin}) will not be rigorously proved, we remark that
the definition (\ref{meas1}), (\ref{meas2}) is meaningful if we have, say $%
P\in L^{1}\left( \left( 0,T\right) :L^{1}\left( \mathcal{U}\right) \right) $
for $0<T\leq \infty .$ In particular the detailed boundary conditions (\ref%
{S5E3}), (\ref{S5E4}), (\ref{S5E7}) or (\ref{T1E1}) are not needed. We can
formulate this more precisely as follows.

\begin{definition}
\label{DefMeas}Given $P\in L^{1}\left( \left( 0,T\right) :L^{1}\left(
\mathcal{U}\right) \right) $ for $0<T\leq \infty $ we define $f\in
L^{1}\left( \left( 0,T\right) :\mathcal{M}\left( \mathcal{V}\right) \right) $
by means of (\ref{meas1}) with $m\left( t\right) $ as in (\ref{meas2}) for $%
a.e.$ $t\in \left( 0,T\right) $.
\end{definition}


%
%
%

\section{Formulation of the adjoint problems.\label{AdjProblems}}

Instead of studying the problem (\ref{S0E1})-(\ref{S0E3}) directly, it is more convenient to study an adjoint version of it. 
The adjoint problem will have several advantages over the original one. 
First the adjoint problem has a simple stochastic interpretation because it is closely related to the generator of the stochastic process which describes the dynamics of the particle. On the other hand, it has good maximum principle properties which will be crucially used to prove Theorem \ref{ThIntro}. Due to this the adjoint equation can be studied using the theory of Markov semigroups in the Banach spaces $C(X)$. 
The analysis of the adjoint problems formulated in this Section will allow us to obtain 
measure valued solutions of
(\ref{S0E1})-(\ref{S0E3}) satisfying the conditions (\ref{S5E4}), (\ref{S5E7}%
), (\ref{S5E3}), (\ref{T1E1}) by duality. 

We begin with the formal derivation of the adjoint problem of (\ref{S0E1})-(\ref{S0E3}) for each of the boundary conditions  (\ref{S5E4}), (\ref{S5E7}), (\ref{S5E3}), (\ref%
{T1E1}).



\begin{definition}
\label{FormAdj}Suppose that $P$ is a smooth function in $\mathcal{U}$ which
satisfies (\ref{Forigin}) and it solves (\ref{S0E1}), (\ref{S0E2}), with one
of the boundary conditions (\ref{S5E3}), (\ref{S5E4}), (\ref{S5E7}), (\ref%
{T1E1}). We will say that the operator $\mathbb{A}$ defined in a set of
functions $\mathcal{D}\left( \mathbb{A}\right) \subset C^{2}\left( \mathbb{R}%
_{+}^{2}\right) $ is the adjoint of the evolution given by (\ref{S0E1}), (%
\ref{S0E2}) and the corresponding boundary condition if for any function
smooth outside the origin $P$ solving (\ref{S0E1}), (\ref{S0E2}) with the
corresponding boundary condition, and for functions $\varphi =\varphi \left(
x,v,t\right) ,$ with $\left( \left\vert x\right\vert +\left\vert
v\right\vert ^{3}\right) ^{\alpha }\left\vert \partial _{t}\varphi
\right\vert ,\ \left( \left\vert x\right\vert +\left\vert v\right\vert
^{3}\right) ^{\alpha }\left\vert \mathbb{A}\left( \varphi \right)
\right\vert \in L^{1}\left( \left[ 0,T\right] \times \mathbb{R}%
_{+}^{2}\right) ,$ $\varphi \left( \cdot ,t\right) \in \mathcal{D}\left(
\mathbb{A}\right) $ for any $t\in \left[ 0,T\right] ,$ the identity:%
\begin{equation}
\int_{\mathcal{V}}\varphi \left( x,v,T\right) f\left( dxdv,T\right) -\int_{%
\mathcal{V}}\varphi \left( x,v,0\right) f\left( dxdv,0\right)
=\int_{0}^{T}\int_{\mathcal{U}}P\left( \varphi _{t}+\mathbb{A}\left( \varphi
\right) \right) dxdvdt  \label{Ad1}
\end{equation}%
holds, where $f$ is as in Definition \ref{DefMeas}.
\end{definition}

\begin{remark}
We will say that $P$ is smooth if $P_{vv},\ P_{x},\ P_{t}$ exist and are
continuous in $\mathcal{U\times }\left( 0,T\right) .$
\end{remark}

\begin{remark}
Notice that (\ref{Ad1}) implies the following. If $\varphi $ satisfies the
conditions in Definition \ref{FormAdj} and, in addition, it has the property
that $\varphi _{t}+\mathbb{A}\left( \varphi \right) =0$ in $\mathcal{U}
\times \left( 0,T \right)$,
then:%
\begin{equation*}
\int_{\mathcal{V}}\varphi \left( x,v,T\right) f\left( dxdv,T\right) =\int_{%
\mathcal{V}}\varphi \left( x,v,0\right) f\left( dxdv,0\right)
\end{equation*}
\end{remark}


In order to derive the form of the operator $\mathbb{A}$ we need to assume
that the class of solutions of the problem (\ref{S0E1}), (\ref{S0E2}), with
any of the boundary conditions (\ref{S5E3}), (\ref{S5E4}), (\ref{S5E7}), (%
\ref{T1E1}), is large enough. This will be formalized in the following
property which we state here for further reference.

\begin{definition}
\label{PropertyP}We denote as $\mathcal{S}$ the set of smooth solutions of
any of the problems (\ref{S0E1}), (\ref{S0E2}) with one of the boundary
conditions (\ref{S5E3}), (\ref{S5E4}), (\ref{S5E7}), (\ref{T1E1}) in an
interval $t\in \left[ 0,T\right] $. We will say that the set $\mathcal{S}$
has the Property (P) if the identity:%
\begin{equation*}
\int_{\mathcal{U\times }\left( 0,T\right) }Phdxdvdt=0\ \ \text{for any }P\in
\mathcal{S}
\end{equation*}%
implies $h=0$ in $\mathcal{U\times }\left( 0,T\right) \ $and
\begin{equation*}
\int_{\left\{ x=0,\ v>0\right\} \mathcal{\times }\left( 0,T\right)
}hPdvdt=0\ \ \text{for any }P\in \mathcal{S}
\end{equation*}%
implies $h=0$ in $\left\{ x=0,\ v>0\right\} \mathcal{\times }\left(
0,T\right) .$
\end{definition}

\begin{remark}
We will not prove in this paper that the solutions of the problem (\ref{S0E1}%
), (\ref{S0E2}) with any of the boundary conditions (\ref{S5E3}), (\ref{S5E4}%
), (\ref{S5E7}), (\ref{T1E1}) has the Property (P) above. We will only use
this definition to derive the form of the adjoint operators $\mathbb{A}$ in
Definition \ref{FormAdj}. We just remark that the Property (P) is a natural
assumption for any equation having a set of solutions sufficiently large.
\end{remark}

\begin{remark}
Smooth solutions in Definition \ref{PropertyP} just means that all the
derivatives appearing in (\ref{S0E1}) exist and are continuous functions in $%
\mathcal{U}$.
\end{remark}

\subsection{Derivation of the adjoint equation and the boundary conditions away
from the singular point.}

It turns out that the adjoint operators $\mathbb{A}$ defined in Definition %
\ref{FormAdj} are given by a second order differential operator in $\mathcal{%
U}$, for functions $\varphi $ such that $\operatorname{supp}\left( \varphi \left(
\cdot ,t\right) \right) \cap \left\{ \left( 0,0\right) \right\} =\varnothing
$ for $t\in \left[ 0,T\right] .$ Moreover, we can obtain also a set of
boundary conditions for the functions in $\mathcal{D}\left( \mathbb{A}%
\right) $ at the boundaries $\partial \mathcal{U}\times \left( 0,\infty \right)
=\left\{ \left( x,v,t\right) =\left( 0,v,t\right) :v\in \mathbb{R}\text{ ,\ }%
t>0\right\},$ where $\mathcal{U}$ and $\partial \mathcal{U}$ are defined as in
(\ref{S8E6}) and (\ref{S8E6a}).
It is worth to remark that the action of the operators $%
\mathbb{A}$, in functions $\varphi $ with $\operatorname{supp}\left( \varphi
\right) \cap \left\{ \left( 0,0\right) \right\} =\varnothing $ as well as
the corresponding boundary conditions, are the same for all the set of
boundary conditions (\ref{S5E3}), (\ref{S5E4}), (\ref{S5E7}), (\ref{T1E1}).

\begin{proposition}
Suppose that $\mathbb{A}$ is the adjoint of the evolution (\ref{S0E1}), (\ref%
{S0E2}) with any of the boundary conditions (\ref{S5E3}), (\ref{S5E4}), (\ref%
{S5E7}), (\ref{T1E1}) in the sense of Definition \ref{FormAdj}. Suppose that
for any $t\in \left[ 0,T\right] $ we have $\varphi \left( \cdot ,t\right)
\in C^{2}\left( \mathcal{V}\right) ,$ with $\operatorname{supp}\left( \varphi
\left( \cdot ,t\right) \right) \cap \left\{ \left( 0,0\right) \right\}
=\varnothing $ for any $t\in \left[ 0,T\right] $ and that $\left\vert
\partial _{t}\varphi \right\vert ,\ \left\vert \mathbb{A}\left( \varphi
\right) \right\vert $ satisfy the integrability conditions in Definition \ref%
{FormAdj}. Suppose that the set of solutions $\mathcal{S}$ defined in
Definition \ref{PropertyP} satisfies Property (P). Then we have:%
\begin{equation}
\mathbb{A}\varphi \left( x,v,t\right) =v\partial _{x}\varphi \left(
x,v,t\right) +\partial _{vv}\varphi \left( x,v,t\right)  \label{AdC1}
\end{equation}%
\begin{equation}
\varphi \left( 0,rv,t\right) =\varphi \left( 0,-v,t\right) \ \ ,\ \ v>0\ \
,\ \ t>0  \label{S6E3}
\end{equation}
\end{proposition}

\begin{proof}
Suppose that $\varphi =\varphi \left( x,v,t\right) $ is a test function
whose support is contained in $\left( x,v,t\right) \in \mathcal{U\times }%
\left( 0,\infty \right) $. Let $P \in \mathcal{S}$.
Multiplying (\ref{S0E1}) by $\varphi $ and
integrating by parts in $\mathcal{U\times }\left( 0,\infty \right) $\ we
obtain:%
\begin{equation}
\int_{\mathcal{U\times }\left( 0,\infty \right) }P\left( -\partial
_{t}\varphi -v\partial _{x}\varphi -\partial _{vv}\varphi \right)
dxdvdt+\int_{\partial \mathcal{U}\times \left( 0,\infty \right) }v\varphi Pdvdt=0
\label{S6E0}
\end{equation}

Suppose first that the support of $\varphi $ does not intersect $L^{\ast
}\times \left( 0,\infty \right) .$ Since $P$ is an arbitrary solution of (%
\ref{S0E1}), (\ref{S0E2}) with one of the conditions (\ref{S5E3}), (\ref%
{S5E4}), (\ref{S5E7}), (\ref{T1E1}) it then follows that:%
\begin{equation}
\partial _{t}\varphi \left( x,v,t\right) +v\partial _{x}\varphi \left(
x,v,t\right) +\partial _{vv}\varphi \left( x,v,t\right) =0\ \ ,\ \ x>0,\ \
v\in \mathbb{R}\ \ ,\ \ t>0  \label{S6E1}
\end{equation}

Now for general $\varphi $, it follows from (\ref{S0E2}) and (\ref{S6E0}) that:%
\begin{equation*}
\int_{0}^{\infty }dt\int_{0}^{\infty }dv\left[ \varphi \left( 0,v,t\right)
-\varphi \left( 0,-\frac{v}{r},t\right) \right] vP\left( 0,v,t\right) =0
\end{equation*}%
Since this identity holds for arbitrary solutions of (\ref{S0E1}), (\ref%
{S0E2}) we can use again Property (P) to obtain (\ref{S6E3}).
\end{proof}

The problem (\ref{S6E1}), (\ref{S6E3}) defines the adjoint problem of (\ref%
{S0E1}), (\ref{S0E2}) for test functions $\varphi $ whose support does not
contain the singular point.\ However, the problem (\ref{S6E1}), (\ref{S6E3})
does not define uniquely an evolution semigroup if $r<r_{c}$ and additional
conditions concerning the asymptotics of $\varphi $ as $\left( x,v\right)
\rightarrow \left( 0,0\right) $ are required in order to prescribe uniquely
an evolution problem for $\varphi .$ In order to determine this set of
boundary conditions we first study the possible asymptotics of the solutions
of (\ref{S6E1}), (\ref{S6E3}) near the singular point. The arguments used to
derive the asymptotics of (\ref{S6E1}), (\ref{S6E3}) will be formal, close
in spirit to those yielding (\ref{Forigin}). Rigorous asymptotic expansions
for the functions $\varphi $ will be made precise and obtained later.

\subsection{Asymptotics of the solutions of (\protect\ref{S6E1}), (\protect\ref%
{S6E3}) near the singular point.}

We compute formally the possible asymptotic behaviour of the solutions of (%
\ref{S6E1}), (\ref{S6E3}) with a method similar to the one used in Section %
\ref{Fsingular} to compute the asymptotics of the solutions of (\ref{S0E1}),
(\ref{S0E2}). More precisely, we look for solutions of (\ref{S6E1}), (\ref%
{S6E3}) with the form:%
\begin{equation}
\varphi \left( x,v,t\right) \sim a\left( t\right) F\left( x,v\right) \ \
\text{as\ \ }\left( x,v\right) \rightarrow \left( 0,0\right)  \label{S6E4}
\end{equation}%
where $F$ behaves algebraically near the singular point and $a\left(
t\right) $ is a smooth function. Arguing as in Section \ref{Fsingular} it
then follows that $F$ must be a solution of the stationary problem:

\begin{align}
v\partial _{x}F+\partial _{vv}F& =0,\;\;x>0,\;\;v\in \mathbb{R},
\label{steadyFPa} \\
F(0,rv)& =F(0,-v),\;\;v>0.  \label{BCa}
\end{align}

The invariance of (\ref{steadyFPa}), (\ref{BCa}) under the rescaling $%
x\rightarrow \lambda ^{3}x,\ v\rightarrow \lambda v$ suggests to look for
solutions of this problem with the form:%
\begin{equation}
F_{\beta }\left( x,v\right) =x^{\beta }\Phi \left( y\right) \ \ ,\ \ y=\frac{%
v^{3}}{9x}  \label{Fa}
\end{equation}

\begin{lemma}
\label{LemmaBeta} There exists a function $\beta \left( r\right) :r\in
\mathbb{R}_{+}\rightarrow \left( -\frac{5}{6},\frac{1}{6}\right) $ such that%
\begin{equation}
-3\beta \left( r\right) \log \left( r\right) +\log \left( 2\sin \left( \pi
\left( \frac{1}{6}-\beta \left( r\right) \right) \right) \right) =0\ \ ,\ \
\beta \left( r\right) \in \left( -\frac{5}{6},\frac{1}{6}\right)
\label{S8E2}
\end{equation}

We have:\
\begin{equation*}
\beta \left( r_{c}\right) =0\ \ ,\ \ \beta (1)=-\frac{2}{3}
\end{equation*}%
with $r_{c}$ as in (\ref{rcrit}). Moreover:%
\begin{equation*}
\lim_{r\rightarrow 0^{+}}\beta (r)=\frac{1}{6}\ \ ,\ \ \lim_{r\rightarrow
\infty }\beta (r)=-\frac{5}{6}
\end{equation*}

The function $\beta (r)$ is related with the function $\alpha \left(
r\right) $ obtained in Proposition \ref{alpha(r)} by means of:%
\begin{equation}
\beta (r)=-\alpha \left( r\right) -\frac{2}{3}  \label{S8E3}
\end{equation}
\end{lemma}

\begin{proof}
The equation (\ref{S8E2}) can be transformed into (\ref{exponent}) by means
of the change of variables (\ref{S8E3}). The result then follows from Lemma %
\ref{alpha(r)}.
\end{proof}

\begin{proposition}
\label{FbetaSt}For any $r>0$ the function $F_{0}\left( x,v\right) =1$ is a
solution of (\ref{steadyFPa}), (\ref{BCa}). Moreover, for any $r>0,$ $r\neq
r_{c},$ there exists another linearly independent positive solution $%
F_{\beta }$ with the form (\ref{Fa}) with $\beta =\beta \left( r\right) $ as
in Lemma \ref{LemmaBeta}. The function $F_{\beta }$ is analytic in $\left\{
\left( x,v\right) :x>0,\ v\in \mathbb{R}\right\} $ and it has the form:%
\begin{equation}
F_{\beta }\left( x,v\right) =x^{\beta }\Phi _{\beta }\left( y\right) \text{
with\ }\Phi _{\beta }\left( y\right) =U(-\beta ,\frac{2}{3};y)\ \ ,\ \ y=%
\frac{v^{3}}{9x}  \label{Fbeta}
\end{equation}

The asymptotics of $\Phi _{\beta }\left( y\right) $ is given by:\
\begin{eqnarray}
\Phi _{\beta }\left( y\right) &\sim &\left\vert y\right\vert ^{\beta }\text{
\ \ \ \ as }y\rightarrow \infty ,\   \label{S8E4} \\
\Phi _{\beta }\left( y\right) &\sim &K\left\vert y\right\vert ^{\beta }\text{
\ \ as }y\rightarrow -\infty  \label{S8E5}
\end{eqnarray}%
where $K=r^{3\beta }.$
\end{proposition}

\begin{proof}
We look for solutions of (\ref{steadyFPa}), (\ref{BCa}) with the form (\ref%
{Fa}). Then $\Phi $ satisfies:
\begin{equation*}
y\Phi _{yy}+\left( \frac{2}{3}-y\right) \Phi _{y}+\beta \Phi =0.
\end{equation*}

Notice that this equation is the same as (\ref{kummer}). The solution of
this equation yielding algebraic behaviour as $\left\vert y\right\vert
\rightarrow \infty $ and satisfying the normalization (\ref{S8E4}) is $\Phi
_{\beta }\left( y\right) =U(-\beta ,\frac{2}{3};y).$ The asymptotics (\ref%
{S8E5}) with $K=r^{3\beta }$ follows from Proposition \ref{Fasympt}. It then
follows from these asymptotic formulas combined with the fact that $F_{\beta
}\left( x,v\right) =x^{\beta }\Phi _{\beta }\left( y\right) $ that $F_{\beta
}\left( 0^{+},v\right) =\frac{1}{9^{\beta }}v^{3\beta }$\ for $v>0$ and $%
F_{\beta }\left( 0^{+},v\right) =\frac{K}{9^{\beta }}\left\vert v\right\vert
^{3\beta }$ for $v<0,$ whence the boundary condition (\ref{BCa})\ follows
using the value of $K.$
\end{proof}

Notice that (\ref{S6E4}) and Lemma \ref{LemmaBeta} suggest the following
asymptotics for the function $\varphi $ near the singular point:%
\begin{equation}
\varphi \left( x,v,t\right) \sim b_{0}\left( t\right) F_{0}\left( x,v\right)
+b_{\beta }\left( t\right) F_{\beta }\left( x,v\right) +...\ \ \text{as\ \ }%
\left( x,v\right) \rightarrow \left( 0,0\right)  \label{S6E5}
\end{equation}

Using the formal asymptotics (\ref{Forigin}), (\ref{S6E5}) we can obtain
precise formulations for the adjoint problems of the problems defined by
means of (\ref{S0E1}), (\ref{S0E2}) with one of the boundary conditions (\ref%
{S5E3}), (\ref{S5E4}), (\ref{S5E7}), (\ref{T1E1}). More precisely, we will
encode the asymptotics (\ref{S6E5}) in the domains of the operators $\mathbb{%
A}$ in Definition \ref{FormAdj}. Our next goal is to define several
operators $\Omega _{\sigma }$ depending on the boundary conditions under
consideration, which will be proved to be the adjoints $\mathbb{A}$ defined
in Definition \ref{FormAdj} for the different sets of boundary conditions
under consideration. We first compute an integral which will be used in the
derivation of some of the adjoint operators. This computation is a bit
tedious and technical, although it just uses classical tools of Complex
Analysis, like contour deformations and the computation of suitable limits.

\subsubsection{Computation of an integral related to particle fluxes.}

We will need to compute the following limit:%
\begin{equation}
C_{\ast }=\lim_{\delta \rightarrow 0}\int_{\partial \mathcal{R}_{\delta }}%
\left[ G_{\alpha }\left( n_{v}D_{v}F_{\beta }+n_{x}vF_{\beta }\right)
-F_{\beta }D_{v}G_{\alpha }n_{v}\right] ds  \label{W1E5}
\end{equation}%
where $\mathcal{R}_{\delta }$ is as in (\ref{domain}) and $n=\left(
n_{x},n_{v}\right) $ is the normal vector to $\partial \mathcal{R}_{\delta }$
pointing towards $\mathcal{R}_{\delta }.$ We have the following result.

\begin{proposition}
\label{Cstar}Suppose that $0<r<r_{c}$ and $\mathcal{R}_{\delta }$ is as in (%
\ref{domain}). Let us assume also that $\alpha $ is as in Lemma \ref%
{alpha(r)}, $\beta $ as in Lemma \ref{LemmaBeta}, $G_{\alpha }$ is defined
as in (\ref{Glambda}) with $\gamma =\alpha $, $F_{\beta }$ as in \
Proposition \ref{FbetaSt} and $n$ is the normal vector to $\partial \mathcal{%
R}_{\delta }$ pointing towards $\mathcal{R}_{\delta }$. Then the constant $%
C_{\ast }$ defined in (\ref{W1E5}) takes the following value:%
\begin{equation}
\frac{C_{\ast }}{9^{\frac{2}{3}}}=\frac{\pi }{3}\left( \sin \left( \pi
\alpha \right) +\sqrt{3}\cos \left( \pi \alpha \right) \right) -2\cos \left(
\pi \left( \beta +\frac{1}{3}\right) \right) \log \left( r\right)
\label{W1E5a}
\end{equation}
\end{proposition}

We will use that $C_{\ast }<0$ for $0<r<r_{c}.$

\begin{lemma}
\label{CstarSign}The constant $C_{\ast }$ defined in (\ref{W1E5}) is
strictly negative for $0<r<r_{c}.$\bigskip
\end{lemma}

\begin{proof}
We can rewrite (\ref{W1E5a}) in an equivalent form. Using (\ref{S1E3}), (\ref%
{W1E5a}), $\alpha +\beta +\frac{2}{3}=0$ and simple trigonometric formulas
we obtain:\
\begin{equation*}
\frac{C_{\ast }}{9^{\frac{2}{3}}}=\frac{4\pi }{3}\sin \left( \pi \left(
\alpha +\frac{2}{3}\right) \right) -2\cos \left( \pi \left( \alpha +\frac{1}{%
3}\right) \right) \log \left( \frac{r}{r_{c}}\right)
\end{equation*}%
Using then (\ref{exponent}) we can rewrite this expression as:%
\begin{equation}
\frac{C_{\ast }}{9^{\frac{2}{3}}}=-\frac{4\pi }{3}\sin \left( \pi \zeta
\right) -2\cos \left( \pi \left( \zeta +\frac{1}{3}\right) \right) \left[
\frac{\log \left( 2\cos \left( \pi \left( \zeta +\frac{1}{3}\right) \right)
\right) }{3\zeta }+\frac{\sqrt{3}\pi }{3}\right]  \label{W1E5b}
\end{equation}%
where $\zeta =-\left( \alpha +\frac{2}{3}\right) .$ Using Lemma \ref%
{alpha(r)} it follows that $C_{\ast }$ would be negative for $0<r<r_{c}$ if
the right-hand side of (\ref{W1E5b}) is negative for $\zeta \in \left( 0,%
\frac{1}{6}\right) .$ This negativity can be proved as follows. The
convexity of the function $\phi \left( x\right) =x\log \left( x\right) -x+1$
for $x\geq 0$ implies the inequality $x\log \left( x\right) \geq x-1$ for $%
x\in \left( 0,1\right) ,$ whence $\log \left( 2A\right) \geq \frac{2A-1}{2A}$
for $A\in \left( 0,\frac{1}{2}\right) .$ Since $\cos \left( \pi \left( \zeta
+\frac{1}{3}\right) \right) \in \left( 0,\frac{1}{2}\right) $ for $\zeta \in
\left( 0,\frac{1}{6}\right) $ we then obtain $\log \left( 2\cos \left( \pi
\left( \zeta +\frac{1}{3}\right) \right) \right) \geq \frac{2\cos \left( \pi
\left( \zeta +\frac{1}{3}\right) \right) -1}{2\cos \left( \pi \left( \zeta +%
\frac{1}{3}\right) \right) }.$ Using this inequality in (\ref{W1E5b}) we
obtain $\frac{C_{\ast }}{9^{\frac{2}{3}}}\leq -\frac{4\pi }{3}\Phi \left(
\zeta \right) $ with
\begin{equation*}
\Phi \left( \zeta \right) =\sin \left( \pi \zeta \right) +\frac{2\cos \left(
\pi \left( \zeta +\frac{1}{3}\right) \right) -1}{4\pi \zeta }+\frac{\sqrt{3}%
}{2}\cos \left( \pi \left( \zeta +\frac{1}{3}\right) \right)
\end{equation*}%
The concavity of the function $\cos \left( \pi \left( \zeta +\frac{1}{3}%
\right) \right) $ for $\zeta \in \left( 0,\frac{1}{6}\right) $ implies 
$$
2\cos \left( \pi \left( \zeta +\frac{1}{3}\right) \right) -1\geq -2\pi \zeta
\sin \left( \pi \left( \zeta +\frac{1}{3}\right) \right) .$$ Then:%
\begin{equation*}
\Phi \left( \zeta \right) >\sin \left( \pi \zeta \right) -\frac{\sin \left(
\pi \left( \zeta +\frac{1}{3}\right) \right) }{2}+\frac{\sqrt{3}}{2}\cos
\left( \pi \left( \zeta +\frac{1}{3}\right) \right) =0\ \ \text{for }\zeta
\in \left( 0,\frac{1}{6}\right)
\end{equation*}%
whence $C_{\ast }<0$ for $0<r<r_{c}.$
\end{proof}

The proof of Proposition \ref{Cstar} is based on elementary arguments such
as the representation of hypergeometric functions in terms of integral
formulas and suitable contour deformations. However, the arguments are
relatively cumbersome and the proof will be split in a sequence of Lemmas.
We first derive a representation formula for $C_{\ast }$ in terms of
confluent hypergeometric functions.

\begin{lemma}
\label{LemmaCstar}Under the assumptions of Proposition \ref{Cstar} we have:%
\begin{equation}
\frac{C_{\ast }}{\left( 9\right) ^{\frac{2}{3}}}=-\lim_{R\rightarrow \infty
}\int_{0}^{R}w\Delta _{\alpha }\left( w\right) dw-2\cos \left( \pi \left(
\beta +\frac{1}{3}\right) \right) \log \left( r\right)  \label{X2}
\end{equation}%
where:%
\begin{equation}
\Delta _{\alpha }=U\left( -\alpha ,\frac{2}{3};-w^{3}\right) U\left( -\beta ,%
\frac{2}{3};w^{3}\right) -U\left( -\alpha ,\frac{2}{3};w^{3}\right) U\left(
-\beta ,\frac{2}{3};-w^{3}\right)  \label{X6}
\end{equation}
\end{lemma}

\begin{proof}
The homogeneity of the integral in (\ref{W1E5}) implies that the integrals $%
\int_{\partial \mathcal{R}_{\delta }}\left[ \cdot \cdot \cdot \right] $ are
independent of $\delta .$ Moreover, integrating by parts and using (\ref%
{steadyFPa}), (\ref{BCa}) we obtain that the integrals $\int_{\partial
\mathcal{R}_{1,b}}\left[ \cdot \cdot \cdot \right] $ take the same value for
$0<b<\infty .$ Therefore, using the form of $n:$%
\begin{eqnarray}
C_{\ast } &=&-\int_{0}^{b}\left[ G_{\alpha }\left( x,r\right) D_{v}F_{\beta
}\left( x,r\right) -F_{\beta }\left( x,r\right) D_{v}G_{\alpha }\left(
x,r\right) \right] dx+  \notag \\
&&+\int_{0}^{b}\left[ G_{\alpha }\left( x,-1\right) D_{v}F_{\beta }\left(
x,-1\right) -F_{\beta }\left( x,-1\right) D_{v}G_{\alpha }\left( x,-1\right) %
\right] dx-  \notag \\
&&-\int_{-1}^{r}vG_{\alpha }\left( b,v\right) F_{\beta }\left( b,v\right) dv
\label{W1E6}
\end{eqnarray}%
for any $b>0.$

Using (\ref{Glambda}), (\ref{LambdaU}) and (\ref{Fbeta}) we obtain:\
\begin{equation*}
\left[ G_{\alpha }\left( x,v\right) D_{v}F_{\beta }\left( x,v\right)
-F_{\beta }\left( x,v\right) D_{v}G_{\alpha }\left( x,v\right) \right] =%
\frac{x^{\alpha +\beta -\frac{1}{3}}}{3}\left( \frac{v}{x^{\frac{1}{3}}}%
\right) ^{2}\Psi \left( \frac{v^{3}}{9x}\right)
\end{equation*}%
where%
\begin{equation}
\Psi \left( s\right) =\left[ U\left( -\alpha ,\frac{2}{3};-s\right) DU\left(
-\beta ,\frac{2}{3};s\right) +U\left( -\beta ,\frac{2}{3};s\right) DU\left(
-\alpha ,\frac{2}{3};-s\right) \right]  \label{W1E8}
\end{equation}

Notice that the two first integrals on the right-hand side of (\ref{W1E6})
have the form:%
\begin{equation}
\int_{0}^{b}\frac{x^{\alpha +\beta -\frac{1}{3}}}{3}\left( \frac{v}{x^{\frac{%
1}{3}}}\right) ^{2}\Psi \left( \frac{v^{3}}{9x}\right) dx\   \label{W1E7}
\end{equation}%
where $v$ takes the values $r$ and $\left( -1\right) $ respectively. Using
the change of variables $x=by$ as well as the fact that $\alpha +\beta +%
\frac{2}{3}=0$ (cf. (\ref{S8E3})) we can transform the integral (\ref{W1E7})
in:\
\begin{equation}
\frac{9^{\frac{2}{3}}}{3}\int_{0}^{1}\left( \frac{L^{3}}{9y}\right) ^{\frac{2%
}{3}}\Psi \left( \frac{L^{3}}{9y}\right) \frac{dy}{y}\   \label{W1E9}
\end{equation}%
where $L=\frac{v}{b^{\frac{1}{3}}}.$ Using (\ref{LambdaU}), (\ref%
{asymptotics}) we obtain $\left\vert \Psi \left( s\right) \right\vert \leq
C\left\vert s\right\vert ^{\alpha +\beta -1}$ whence $\left\vert s^{\frac{2}{%
3}}\Psi \left( s\right) \right\vert \leq C\left\vert s\right\vert ^{\alpha
+\beta +\frac{2}{3}-1}=\frac{C}{\left\vert s\right\vert }.$ Therefore the
integral in (\ref{W1E9}) can be estimated as $\frac{C}{L^{3}}$ whence the
limit of the integrals in (\ref{W1E7}) converges to zero as $b\rightarrow 0.$
It then follows from (\ref{W1E6}) that:%
\begin{equation}
C_{\ast }=-\lim_{b\rightarrow 0}\int_{-1}^{r}vG_{\alpha }\left( b,v\right)
F_{\beta }\left( b,v\right) dv  \label{W2E1}
\end{equation}

We now notice that (\ref{Glambda}), (\ref{LambdaU}), (\ref{Fbeta}) as well
as $\alpha +\beta +\frac{2}{3}=0$ yield:%
\begin{equation*}
vG_{\alpha }\left( b,v\right) F_{\beta }\left( b,v\right) =\frac{9^{\frac{2}{%
3}}}{\left( 9b\right) ^{\frac{1}{3}}}U\left( -\alpha ,\frac{2}{3};-\frac{%
v^{3}}{9b}\right) U\left( -\beta ,\frac{2}{3};\frac{v^{3}}{9b}\right) \frac{v%
}{\left( 9b\right) ^{\frac{1}{3}}}
\end{equation*}

Then, using the change of variables $w=\frac{v}{\left( 9b\right) ^{\frac{1}{3%
}}}$ and writing $R=\frac{1}{\left( 9b\right) ^{\frac{1}{3}}}$ we obtain:%
\begin{equation*}
C_{\ast }=-\left( 9\right) ^{\frac{2}{3}}\lim_{R\rightarrow \infty
}\int_{-R}^{rR}U\left( -\alpha ,\frac{2}{3};-w^{3}\right) U\left( -\beta ,%
\frac{2}{3};w^{3}\right) wdw
\end{equation*}

Splitting the integral as $\int_{-R}^{rR}\left[ \cdot \cdot \cdot \right]
=\int_{0}^{rR}\left[ \cdot \cdot \cdot \right] +\int_{-rR}^{0}\left[ \cdot
\cdot \cdot \right] +\int_{-R}^{-rR}\left[ \cdot \cdot \cdot \right] $ and
using the change of variables $w\rightarrow \left( -w\right) $ in the last
two integrals we obtain:\
\begin{eqnarray}
C_{\ast } &=&-\left( 9\right) ^{\frac{2}{3}}\lim_{R\rightarrow \infty
}\int_{0}^{rR}[U\left( -\alpha ,\frac{2}{3};-w^{3}\right) U\left( -\beta ,%
\frac{2}{3};w^{3}\right)  \notag \\
&&~\ ~\ ~~~~\ ~\ ~\ ~~~~\ ~-U\left( -\alpha ,\frac{2}{3};w^{3}\right)
U\left( -\beta ,\frac{2}{3};-w^{3}\right) ]wdw  \notag \\
~~~ &&-\lim_{R\rightarrow \infty }\int_{rR}^{R}U\left( -\alpha ,\frac{2}{3}%
;w^{3}\right) U\left( -\beta ,\frac{2}{3};-w^{3}\right) wdw  \label{W2E2}
\end{eqnarray}

Using (\ref{LambdaU}), (\ref{asymptotics}) as well as $\alpha +\beta +\frac{2%
}{3}=0,$ we obtain:%
\begin{equation*}
U\left( -\alpha ,\frac{2}{3};w^{3}\right) U\left( -\beta ,\frac{2}{3}%
;-w^{3}\right) w\sim 2\cos \left( \pi \left( \beta +\frac{1}{3}\right)
\right) \left( w\right) ^{-1}\text{ as }w\rightarrow \infty
\end{equation*}%
whence:%
\begin{equation*}
\lim_{R\rightarrow \infty }\int_{rR}^{R}U\left( -\alpha ,\frac{2}{3}%
;w^{3}\right) U\left( -\beta ,\frac{2}{3};-w^{3}\right) wdw=-2\cos \left(
\pi \left( \beta +\frac{1}{3}\right) \right) \log \left( r\right)
\end{equation*}

Plugging this identity in (\ref{W2E2}) and replacing $rR$ by $R$ in the
first limit on the right we obtain (\ref{X2}) and the result follows.
\end{proof}

In order to compute the value of $\Delta _{\alpha }$ in (\ref{X6}) we will
use some representation formulas for the functions $U\left( -\alpha ,\frac{2%
}{3};-w^{3}\right) ,\ U\left( -\beta ,\frac{2}{3};w^{3}\right) :$

\begin{lemma}
\label{Urep}Suppose that $0<r<r_{c}$ and $\alpha $ , $\beta $ are as in the
Lemmas \ref{alpha(r)}, \ref{LemmaBeta} respectively. Then, the following
representation formulas hold:%
\begin{eqnarray}
U\left( -\alpha ,\frac{2}{3};w^{3}\right) &=&\frac{e^{w^{3}}}{\Gamma \left(
-\alpha \right) }Q_{2}\left( w;\alpha +1\right) \ \ ,\ \ w>0\   \label{W2E3}
\\
U\left( -\beta ,\frac{2}{3};w^{3}\right) &=&\frac{we^{w^{3}}}{\Gamma \left(
\frac{1}{3}-\beta \right) }Q_{1}\left( w;\beta +\frac{2}{3}\right) \ \ ,\ \
w>0  \label{W2E4}
\end{eqnarray}%
where the functions $Q_{n}$ are defined by means of:%
\begin{equation}
Q_{n}\left( w;a\right) =\int_{1}^{\infty }e^{-w^{3}t}\left( \frac{t}{t-1}%
\right) ^{a}\frac{dt}{t^{\frac{2n}{3}}}\ \ ,\ \ 0<a<1,\ \ n>0\ \ ,\ \ w>0
\label{W2E5}
\end{equation}
\end{lemma}

\begin{remark}
Notice that the functions $Q_{n}$ are defined for noninteger values of $n.$
\end{remark}

\begin{proof}
We have:%
\begin{eqnarray}
U\left( a,b,z\right) &=&\frac{e^{z}}{\Gamma \left( a\right) }%
\int_{1}^{\infty }e^{-zt}\left( t-1\right) ^{a-1}\left( t\right) ^{b-a-1}dt\
,\ a>0\ ,\ \operatorname{Re}\left( z\right) >0\   \label{Ab1} \\
U\left( a,b,z\right) &=&z^{1-b}U\left( 1+a-b,2-b,z\right) \   \label{Ab2}
\end{eqnarray}%
(cf. 13.2.6 and 13.1.29 from \cite{Ab}). Using that $\alpha <0$ we then
obtain from (\ref{Ab1}):%
\begin{equation}
U\left( -\alpha ,\frac{2}{3};w^{3}\right) =\frac{e^{w^{3}}}{\Gamma \left(
-\alpha \right) }\int_{1}^{\infty }e^{-w^{3}t}\left( t-1\right) ^{-\alpha
-1}\left( t\right) ^{\alpha -\frac{1}{3}}dt\ \ ,\ \ w>0  \label{X4}
\end{equation}%
whence (\ref{W2E3}) follows. On the other hand (\ref{Ab2}) yields
\begin{equation*}
U\left( -\beta ,\frac{2}{3};w^{3}\right) =wU\left( \frac{1}{3}-\beta ,\frac{4%
}{3};w^{3}\right) .
\end{equation*}
Using then that $\beta <\frac{1}{3}$ we obtain from (\ref{Ab1}):%
\begin{equation}
U\left( -\beta ,\frac{2}{3};w^{3}\right) =\frac{we^{w^{3}}}{\Gamma \left(
\frac{1}{3}-\beta \right) }\int_{1}^{\infty }e^{-w^{3}t}\left( t-1\right)
^{-\left( \beta +\frac{2}{3}\right) }t^{\beta }dt\ \ ,\ \ w>0  \label{X5}
\end{equation}%
and (\ref{W2E4}) follows.
\end{proof}

We now prove that the functions $Q_{n}\left( w;a\right) $ can be extended
analytically for $w\neq 0$ and derive suitable representation formulas. To
this end we need to give a precise definition of some analytic functions
with branch points.

\begin{definition}
\label{branches}We define a branch of the function $\left( \frac{t}{t-1}%
\right) ^{a}$ analytic in $\mathbb{C}\diagdown \left[ 0,1\right] ,$ which
will be denoted as $\Phi \left( t;a\right) ,$ prescribing that $\left(
s\right) ^{a}=\left\vert s\right\vert ^{a}e^{ia\arg \left( s\right) },$ with
$\arg \left( s\right) \in \left( -\pi ,\pi \right) $ for $s\in \mathbb{C}%
\diagdown \left[ 0,\infty \right) .$ On the other hand, we define the
functions $t^{\frac{2n}{3}}$ with $n=1,2$ in a subset of a Riemann surface
given by $\mathcal{S}=\left\{ t=\left\vert t\right\vert e^{i\theta
}:\left\vert t\right\vert \neq 0,\ \theta \in \left[ -3\pi ,0\right]
\right\} $. We set:%
\begin{equation*}
\left( \left\vert t\right\vert e^{i\theta }\right) ^{\frac{2n}{3}%
}=\left\vert t\right\vert ^{\frac{2n}{3}}e^{\frac{2n\theta i}{3}}
\end{equation*}

There exists a natural projection from the Riemann surface $\mathcal{S}$ to $%
\mathbb{C}$. We will say that a contour $\Lambda $ defined in $\mathcal{S}$
surrounds the interval $\left[ 0,1\right] $ if the projection of $\Lambda $
into $\mathbb{C}$ surrounds the interval $\left[ 0,1\right] .$
\end{definition}

\begin{remark}
We use the notation $t=\left\vert t\right\vert e^{i\theta }$ to denote
points in the Riemann surface $\mathcal{S}$ with a value of the phase $%
\theta .$ Notice that then the points $t_{0}=1\in \mathcal{S}$ and $%
t_{0}=e^{-2\pi i}\in \mathcal{S}$ are different points.
\end{remark}

We can then obtain the following:

\begin{lemma}
\label{Qanalytic} The functions $Q_{n}\left( w;a\right) $ defined in (\ref%
{W2E5}) for $0<a<1,\ n=1,2$ and $w>0$ can be extended analytically to the
set $\left\{ w\in \mathbb{C}:\operatorname{Im}\left( w\right) >0\right\} $ and
continuously to the set $\left\{ \operatorname{Im}\left( w\right) \geq 0,\ w\neq
0\right\} .$\ Moreover, the following representation formulas hold for $w<0:$%
\begin{equation}
Q_{n}\left( w;a\right) =\int_{\gamma }e^{-w^{3}t}\Phi \left( t;a\right)
\frac{dt}{t^{\frac{2n}{3}}}\ \ ,\ \ n>0\ ,\ 0<a<1  \label{W3E3}
\end{equation}%
where the functions $\Phi \left( t\right) $ and $t^{\frac{2n}{3}}$ are as in
Definition \ref{branches} and $\gamma $ is a contour starting at $t=1$
contained in the Riemann surface $\mathcal{S}$, surrounding the interval $%
\left[ 0,1\right] $ and approaching asymptotically to $t=\infty \cdot
e^{-3\pi i}$.
\end{lemma}

\begin{figure}[tbph]
\centering \includegraphics[width=0.5\linewidth]{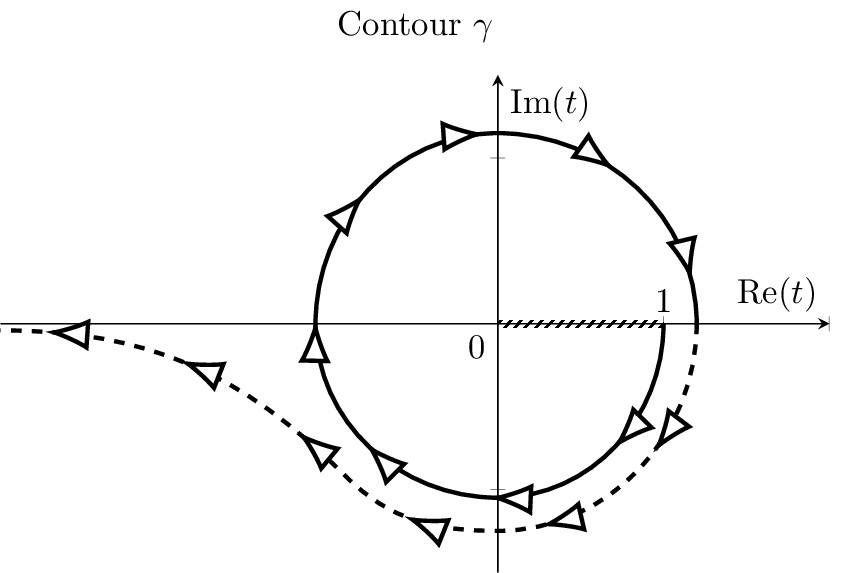}
\caption{Contour $\protect\gamma $ }
\label{fig:figure-1}
\end{figure}

\begin{remark}
Notice that Lemma \ref{Qanalytic} yields a representation formula for
\begin{equation*}
U\left( -\alpha ,\frac{2}{3};w^{3}\right) ,\ U\left( -\beta ,\frac{2}{3}%
;w^{3}\right) \text{ with }w<0
\end{equation*}
by means of (\ref{W2E3}), (\ref{W2E4}).
\end{remark}

\begin{proof}
We can then rewrite (\ref{W2E5}) as:
\begin{equation}
Q_{n}\left( w;a\right) =\int_{C}e^{-w^{3}t}\left( \frac{t}{t-1}\right) ^{a}%
\frac{dt}{t^{\frac{2n}{3}}}\ \ ,\ \ w>0\   \label{X7}
\end{equation}%
where $C\subset \mathcal{S}$ is the line connecting $1$ and $\infty $ given
by $\left\{ t:t=1+\rho e^{i\theta },\ \rho \in \left[ 0,\infty \right) ,\
\theta =0\right\} .$ In order to extend analytically the function $%
Q_{n}\left( w;a\right) $ to the region $\left\{ \operatorname{Im}\left( w\right)
>0\right\} $ we use (\ref{X7}) with $w=\left\vert w\right\vert e^{\varphi i}$
with $\varphi \in \left[ 0,\pi \right] .$ In order to ensure the convergence
of the integrals we can modify the contour of integration to a new contour $%
C_{\varphi }$ which connects $t=1$ with $t=\infty $ and for sufficiently
large $\left\vert t\right\vert ,$ say $\left\vert t\right\vert \geq 3$ is
just the line $\left\{ t:t=\left\vert t\right\vert e^{-3\varphi i}\right\} .$
Therefore, if $\varphi $ varies from $0$ to $\pi $ we would obtain that this
line which describes the asymptotics of the new contours $C_{\varphi }$ is
just the line $\left\{ t:t=\left\vert t\right\vert e^{-3\pi i}\right\} $ of
the set $\mathcal{S}$. Notice that this contour deformation must be made
avoiding intersections of the new contours $C_{\varphi }$ with the interval $%
\left[ 0,1\right] .$ The function $\Phi \left( t;a\right) $ can be defined
in $\mathcal{S}$ in a natural manner using just its definition in $\mathbb{C}%
\diagdown \left[ 0,1\right] .$ After concluding the deformation of the
contours we obtain a representation formula for $Q_{n}\left( w;a\right) $
with $w<0$ having the form (\ref{W3E3}).
\end{proof}

\bigskip

\begin{figure}[htbp]
\centering \includegraphics[width=1.0\linewidth]{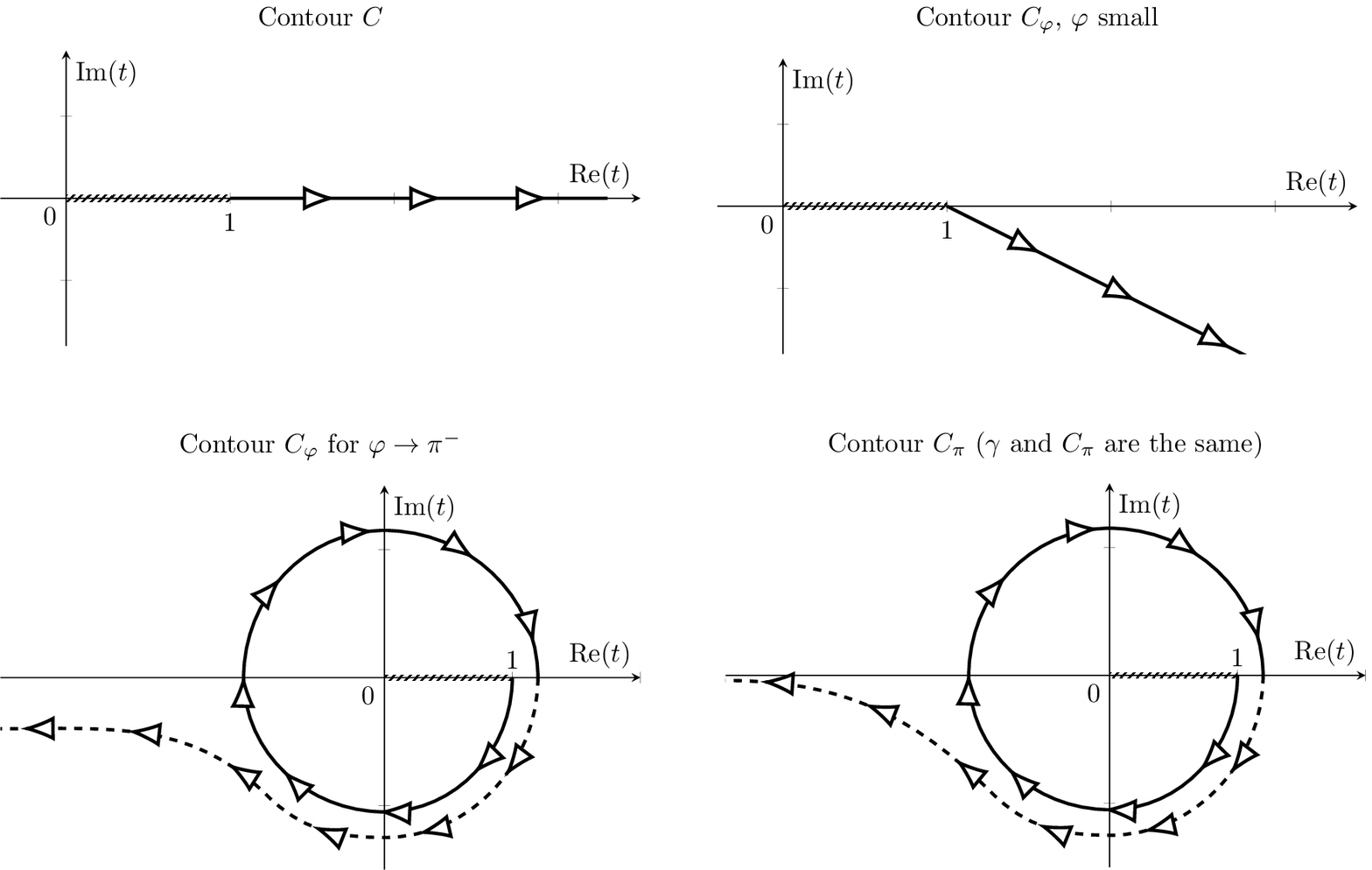}
\caption{Contours $C_{\protect\varphi }$\texttt{\ for some values of }$%
\protect\varphi$ }
\label{fig:figure-2}
\end{figure}

We can now find a representation formula for the limit $\lim_{R\rightarrow
\infty }\int_{0}^{R}w\Delta _{\alpha }\left( w\right) dw$ in (\ref{X2}). We
first define the following family of auxiliary functions:

\begin{definition}
\label{PsiAB}Suppose that $0<A<1,\ 0<B<1,\ n>0,\ m>0.$ We define:%
\begin{equation*}
\Psi \left( t,x;A,B,n,m\right) =\frac{\Phi \left( t;A\right) \Phi \left(
x;B\right) }{t^{\frac{2m}{3}}x^{\frac{2n}{3}}}+\frac{\Phi \left( t;B\right)
\Phi \left( x;A\right) }{x^{\frac{2m}{3}}t^{\frac{2n}{3}}},\ t\in \mathcal{S}%
\text{ \ ,\ }x\in \mathcal{S}
\end{equation*}%
where the functions $\Phi $ are as in Definition \ref{branches} and the
power laws $s\rightarrow \left( s\right) ^{a}$ are computed in the portion
of Riemann surface $\mathcal{S}$ as in Definition \ref{branches}.
\end{definition}

We then have:

\begin{lemma}
\label{LemmaExp}Suppose that we define the following real functions for $%
x>1,\ 0<t<1:$\
\begin{eqnarray}
\Psi _{1}\left( t,x;A,B,n,m\right) &=&\left( \frac{t}{1-t}\right) ^{A}\left(
\frac{x}{x-1}\right) ^{B}\frac{1}{t^{\frac{2m}{3}}x^{\frac{2n}{3}}}
\label{W4E6a} \\
\Psi _{2}\left( t,x;A,B,n,m\right) &=&\left( \frac{t}{1-t}\right) ^{B}\left(
\frac{x}{x-1}\right) ^{A}\frac{1}{t^{\frac{2n}{3}}x^{\frac{2m}{3}}}
\label{W4E6b}
\end{eqnarray}

We define a contour $\Lambda $ in $\mathcal{S}$ as $\Lambda =\Lambda
_{1}\cup \Lambda _{2}$ where:%
\begin{equation}
\Lambda _{1}=\left\{ t\in \mathcal{S}:t=\frac{e^{i\theta }}{2},\ \theta \in %
\left[ 0,-2\pi \right] \right\} \   \label{W4E1}
\end{equation}%
and $\Lambda _{2}$ is a contour connecting $\frac{e^{-2\pi i}}{2}$ with $%
t=\infty \cdot e^{-3\pi i}.$\ We define:%
\begin{equation}
K_{m}\left( A\right) =\left( -e^{i\pi A}+e^{-i\pi A}e^{\frac{4\pi mi}{3}%
}-e^{i\pi A}e^{\frac{4\pi mi}{3}}\right) \ ,\ 0<A<1,\ m>0\   \label{W4E2}
\end{equation}%
Let $A=\alpha +1,\ B=\beta +\frac{2}{3},\ n=1,\ m=2.$ Let
\begin{equation}
G\left( t,x;A,B,n,m\right) =K_{m}\left( A\right) \Psi _{1}+K_{n}\left(
B\right) \Psi _{2}  \label{W4E2a}
\end{equation}%
Then the function $\frac{\left\vert G\left( t,x\right) \right\vert }{%
\left\vert t-x\right\vert }$ is integrable in $\left( t,x\right) \in \left[
\frac{1}{2},1\right] \times \left( 1,\infty \right) $ and we have:
\begin{eqnarray}
&&\lim_{R\rightarrow \infty }\int_{0}^{R}w\Delta _{\alpha }\left( w\right) dw
\label{W4E5} \\
&=&-\frac{\sin \left( \pi \alpha \right) }{3\pi }\left[ \int_{1}^{\infty
}dx\int_{\Lambda }dt\frac{\Psi \left( t,x;A,B,n,m\right) }{\left( x-t\right)
}+\right.  \notag \\
&&\left. +\int_{1}^{\infty }dx\int_{\frac{1}{2}}^{1}dt\frac{G\left(
t,x;A,B,n,m\right) }{\left( x-t\right) }\right]  \notag
\end{eqnarray}
\end{lemma}

\begin{remark}
The rationale behind the definition of the contour $\Lambda $ is to avoid
the contour of integration approaching the origin $t=0,$ because the
function $\Psi $ is not integrable there for the range of parameters
required.
\end{remark}

\begin{remark}
\label{contour}We introduce here a contour $\tilde{C}$ for further
reference. This contour consists in the limit of contours approaching the
segment $\left[ \frac{1}{2},1\right] ,$ connecting the points $t=1$ with $t=%
\frac{1}{2}$ with $t=\left\vert t\right\vert e^{i0}$ and $\operatorname{Im}\left(
t\right) \rightarrow 0^{-}.$ We continue the contour by means of the contour
$\Lambda _{1}$ defined in (\ref{W4E1}). It is then followed by a contour
obtained as limit of contours converging to the interval $\left[ \frac{1}{2}%
,1\right] $ with $\operatorname{Im}\left( t\right) >0,\ t=\left\vert t\right\vert
e^{-2\pi i}$ and connecting the points $t=\frac{1}{2},1.$ We then continue
the contour by a segment converging to $\left[ \frac{1}{2},1\right] $ with $%
\operatorname{Im}\left( t\right) <0,$ $t=\left\vert t\right\vert e^{-2\pi i},$
connecting $t=1,\frac{1}{2}.$ The last part of the contour is then the
contour $\Lambda _{2}$ in the statement of the Lemma.
\end{remark}

\bigskip
\begin{figure}[htbp]
\centering \includegraphics[width=0.6
\linewidth]{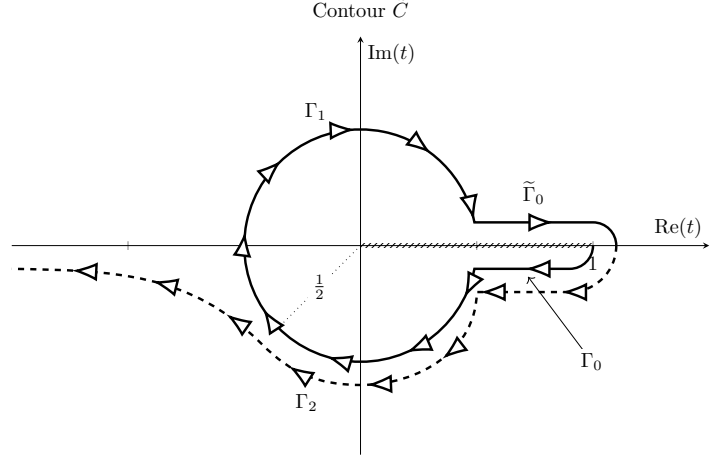}
\caption{Contour $\tilde{C}$ }
\label{fig:figure-3}
\end{figure}

\begin{proof}
Using (\ref{X6}), (\ref{W2E3}), (\ref{W2E4}), (\ref{W3E3}) we obtain the
representation formula:%
\begin{eqnarray}
\Delta _{\alpha } &=&\frac{w}{\Gamma \left( -\alpha \right) \Gamma \left(
\frac{1}{3}-\beta \right) }\left[ \int_{1}^{\infty }e^{-w^{3}x}\left( \frac{x%
}{x-1}\right) ^{\beta +\frac{2}{3}}\frac{dx}{x^{\frac{2}{3}}}\int_{\gamma
}e^{w^{3}t}\Phi \left( t;\alpha +1\right) \frac{dt}{t^{\frac{4}{3}}}\right.
\label{W4E3} \\
&&\left. +\int_{1}^{\infty }e^{-w^{3}x}\left( \frac{x}{x-1}\right) ^{\alpha
+1}\frac{dx}{x^{\frac{4}{3}}}\int_{\gamma }e^{w^{3}t}\Phi \left( t;\beta +%
\frac{2}{3}\right) \frac{dt}{t^{\frac{2}{3}}}\right]  \notag \\
&\equiv &\frac{w}{\Gamma \left( -\alpha \right) \Gamma \left( \frac{1}{3}%
-\beta \right) }\left( J_{1}+J_{2}\right)  \notag
\end{eqnarray}%
where $\gamma \in \mathcal{S}$ is as in Lemma \ref{Qanalytic}. We can deform
the contour of integration $\gamma $ to the contour $\tilde{C}$ introduced
in Remark \ref{contour}. Using the Definition of the functions $\Phi \left(
t;\alpha +1\right) ,\ t^{\frac{4}{3}}$ we obtain the following formula for
the integral $\int_{\gamma }=\int_{\tilde{C}}$ in the definition of $J_{1}:$%
\begin{eqnarray*}
&&\int_{\gamma }e^{w^{3}t}\Phi \left( t;\alpha +1\right) \frac{dt}{t^{\frac{4%
}{3}}} \\
&=&\int_{\Lambda }e^{w^{3}t}\Phi \left( t;\alpha +1\right) \frac{dt}{t^{%
\frac{4}{3}}} \\
&&+\left( -e^{i\pi \left( \alpha +1\right) }+e^{-i\pi \left( \alpha
+1\right) }e^{\frac{8\pi i}{3}}-e^{i\pi \left( \alpha +1\right) }e^{\frac{%
8\pi i}{3}}\right) \int_{\frac{1}{2}}^{1}e^{w^{3}t}\left( \frac{t}{1-t}%
\right) ^{\alpha +1}\frac{dt}{t^{\frac{4}{3}}} \\
&=&\int_{\Lambda }e^{w^{3}t}\Phi \left( t;\alpha +1\right) \frac{dt}{t^{%
\frac{4}{3}}}+K_{2}\left( \alpha +1\right) \int_{\frac{1}{2}%
}^{1}e^{w^{3}t}\left( \frac{t}{1-t}\right) ^{\alpha +1}\frac{dt}{t^{\frac{4}{%
3}}}
\end{eqnarray*}%
where we use that:%
\begin{eqnarray}
\lim_{\varepsilon \rightarrow 0^{+}}\Phi \left( x+\varepsilon i\right)
&=&\left( \frac{x}{1-x}\right) ^{a}e^{-i\pi a}\ \ ,\ \ x\in \left( 0,1\right)
\label{W2E8a} \\
\lim_{\varepsilon \rightarrow 0^{+}}\Phi \left( x-\varepsilon i\right)
&=&\left( \frac{x}{1-x}\right) ^{a}e^{i\pi a}\ \ ,\ \ \ x\in \left(
0,1\right)  \label{W2E8b}
\end{eqnarray}

A similar argument gives:%
\begin{eqnarray*}
&&\int_{\gamma }e^{w^{3}t}\Phi \left( t;\beta +\frac{2}{3}\right) \frac{dt}{%
t^{\frac{2}{3}}} \\
&=&\int_{\Lambda }e^{w^{3}t}\Phi \left( t;\beta +\frac{2}{3}\right) \frac{dt%
}{t^{\frac{2}{3}}}+K_{1}\left( \beta +\frac{2}{3}\right) \int_{\frac{1}{2}%
}^{1}e^{-w^{3}t}\left( \frac{t}{1-t}\right) ^{\beta +\frac{2}{3}}\frac{dt}{%
t^{\frac{2}{3}}}
\end{eqnarray*}

We then obtain:%
\begin{eqnarray}
J_{1}+J_{2} &=&\int_{1}^{\infty }dx\int_{\Lambda }dte^{-w^{3}\left(
x-t\right) }\Psi \left( t,x;A,B,n,m\right) +  \label{W4E4} \\
&&+\int_{1}^{\infty }dx\int_{\frac{1}{2}}^{1}dte^{-w^{3}\left( x-t\right)
}G\left( t,x;A,B,n,m\right)  \notag
\end{eqnarray}

Notice that the integrability of $\frac{\left\vert G\left( t,x\right)
\right\vert }{\left\vert t-x\right\vert }$ near $x=1$, $t=1$ is a
consequence of the fact that we have $K_{2}\left( \alpha +1\right)
+K_{1}\left( \beta +\frac{2}{3}\right) =0$ and this implies the inequality:%
\begin{eqnarray}
&&\left\vert G\left( t,x;\alpha +1,\beta +\frac{2}{3},1,2\right) \right\vert
\label{W4E8} \\
&\leq &C\left( \frac{1}{\left( 1-t\right) ^{\alpha +1}\left( x-1\right)
^{\beta +\frac{2}{3}}}+\frac{1}{\left( x-1\right) ^{\alpha +1}\left(
1-t\right) ^{\beta +\frac{2}{3}}}\right) \left\vert x-t\right\vert  \notag
\end{eqnarray}%
and the desired integrability follows since $\left( \alpha +1\right) <1,$ $%
\left( \beta +\frac{2}{3}\right) <1.$ The fact that $K_{2}\left( \alpha
+1\right) +K_{1}\left( \beta +\frac{2}{3}\right) =0$ follows from the
following computation which is a consequence of the fact that $\beta +\frac{2%
}{3}=-\alpha :$\
\begin{equation*}
K_{2}\left( \alpha +1\right) +K_{1}\left( \beta +\frac{2}{3}\right) =\left(
e^{i\pi \alpha }-e^{-i\pi \alpha }\right) \left( 1+e^{\frac{4\pi i}{3}}+e^{%
\frac{8\pi i}{3}}\right) =0
\end{equation*}

Integrating (\ref{W4E3}) in $\left[ 0,R\right] $ after multiplying by $w,$
using (\ref{W4E4}) and taking the limit $R\rightarrow \infty $ we obtain (%
\ref{W4E5}), using also the formula $\Gamma \left( \beta +\frac{2}{3}\right)
\Gamma \left( \frac{1}{3}-\beta \right) =\frac{\pi }{\sin \left( \pi \left(
\beta +\frac{2}{3}\right) \right) }=-\frac{\pi }{\sin \left( \pi \alpha
\right) }$ (cf. 6.1.17 in \cite{Ab})\ we obtain (\ref{W4E5}).
\end{proof}

We now study the properties of the functions defined by means of the
right-hand side of (\ref{W4E5}). More precisely, we define the following
subset of $\mathbb{R}^{4}$:%
\begin{equation}
\mathcal{K}=\{\left( A,B,n,m\right) :0<A,B<1,\ n>0,\ m>0,\ \frac{2\left(
n+m\right) }{3}>1,K_{m}\left( A\right) +K_{n}\left( B\right) =0\ \}
\label{W4E7}
\end{equation}%
where the functions $K_{m}\left( A\right) ,\ K_{n}\left( B\right) $ are
defined as in (\ref{W4E2}). Therefore:

\begin{lemma}
\label{DefFunct}The formula:%
\begin{eqnarray}
Q\left( A,B,n,m\right) &=&\left[ \int_{1}^{\infty }dx\int_{\Lambda }dt\frac{%
\Psi \left( t,x;A,B,n,m\right) }{\left( x-t\right) }+\right.  \label{W4E9} \\
&&\left. +\int_{1}^{\infty }dx\int_{\frac{1}{2}}^{1}dt\frac{G\left(
t,x;A,B,n,m\right) }{\left( x-t\right) }\right] \   \notag
\end{eqnarray}%
defines a continuous function in the set $\mathcal{K}$ in (\ref{W4E7}) where
the contour $\Lambda $ is as in Lemma \ref{LemmaExp}, $\Psi $ is as in
Definition \ref{PsiAB}, and $G\left( t,x;A,B,n,m\right) $ as in (\ref{W4E2a}%
) with $K_{m}\left( A\right) ,\ K_{n}\left( B\right) $ as in (\ref{W4E2})
and $\Psi _{1},\ \Psi _{2}$ are as in (\ref{W4E6a}), (\ref{W4E6b}).

Suppose that in addition to $\left( A,B,n,m\right) \in \mathcal{K}$ we have $%
A+B<1.$ Then the following representation formula holds:%
\begin{equation}
Q\left( A,B,n,m\right) =\int_{1}^{\infty }dx\int_{\tilde{C}}dt\frac{\Psi
\left( t,x;A,B,n,m\right) }{\left( x-t\right) }\   \label{W5E1}
\end{equation}%
where $\tilde{C}$ is the contour introduced in Remark \ref{contour}.
\end{lemma}

\begin{proof}
Since $\left( A,B,n,m\right) \in \mathcal{K}$ it follows from the definition
of $G,\ K_{m}\left( A\right) ,\ K_{n}\left( B\right) ,$ and$\ \Psi _{1},\
\Psi _{2}$ that the inequality (\ref{W4E8}) holds. Therefore the formula (%
\ref{W4E9}) defines a continuous function $Q$ in $\mathcal{K}$. Suppose that
$A+B<1.$ Then, since the contour $\tilde{C}$ is contained in $\operatorname{Re}%
\left( t\right) \leq 1$ we obtain, using Definition \ref{PsiAB} the
following estimate for $t\in \Lambda \cap \left\{ \left\vert t\right\vert
\leq 2\right\} ,\ x\in \left\{ \left\vert x\right\vert \leq 2\right\} :$%
\begin{equation*}
\left\vert \frac{\Psi \left( t,x;A,B,n,m\right) }{x-t}\right\vert \leq
K\left( \frac{1}{\left( 1-t\right) ^{A+A^{\ast }}\left( x-1\right)
^{B+B^{\ast }}}+\frac{1}{\left( 1-t\right) ^{B+B^{\ast }}\left( x-1\right)
^{A+A^{\ast }}}\right)
\end{equation*}%
where $A^{\ast }+B^{\ast }=1,$ $A+A^{\ast }<1,\ B+B^{\ast }<1$ and $K$ is
independent of $A,B,n,m.$ Notice that the choice of $A^{\ast },\ B^{\ast }$
is possible because $A+B<1.$ Then, the function $\left\vert \frac{\Psi
\left( t,x;A,B,n,m\right) }{x-t}\right\vert $ is integrable in a
neighbourhood of $t=x=1,$ and using the definition of the functions $\Phi ,\
\left( t\right) ^{a}$ in Definition \ref{branches} we obtain exactly the
numerical factors in the functions $K_{m}\left( A\right) ,\ K_{n}\left(
B\right) $ in the integral in the contour $\tilde{C}.$ This gives the
representation formula (\ref{W5E1}).
\end{proof}

We now describe the structure of the set $\mathcal{K}$ in the neighbourhood
of a point $\left( A_{0},B_{0},1,2\right) $ with $A_{0}+B_{0}=1,$ $A_{0}>0,\
B_{0}>0.$ Notice that a computation similar to (\ref{W5E1}) shows that such
points are contained in $\mathcal{K}$. We have:

\begin{lemma}
\label{localK}For each positive $A_{0},\ B_{0}$ satisfying $A_{0}+B_{0}=1$
there exists $\delta >0$ and two differentiable functions $\Theta _{1}\left(
A,B\right) ,\ \Theta _{2}\left( A,B\right) $ defined in $\left\vert
A-A_{0}\right\vert +\left\vert B-B_{0}\right\vert <\delta $ such that $%
\Theta _{1}\left( A_{0},B_{0}\right) =1,\ \Theta _{2}\left(
A_{0},B_{0}\right) =2$ and such that $\left( A,B,\Theta _{1}\left(
A,B\right) ,\Theta _{2}\left( A,B\right) \right) $ belong to $\mathcal{K}$.
Moreover, we have:%
\begin{eqnarray}
\lim_{\varepsilon \rightarrow 0}\frac{\Theta _{1}\left(
A_{0},B_{0}-\varepsilon \right) -1}{\varepsilon } &=&-\frac{3}{8}\left( 1+%
\sqrt{3}\cot \left( \pi A_{0}\right) \right) \ \   \label{W5E3} \\
\ \ \lim_{\varepsilon \rightarrow 0}\frac{\Theta _{2}\left(
A_{0},B_{0}-\varepsilon \right) -2}{\varepsilon } &=&-\frac{3}{8}\left( 1+%
\sqrt{3}\cot \left( \pi A_{0}\right) \right) \   \label{W5E3b}
\end{eqnarray}
\end{lemma}

\begin{proof}
The set $\mathcal{K}$ is defined by means of the equation $K_{m}\left(
A\right) +K_{n}\left( B\right) =0.$ Taking the real and imaginary part of
this equation we obtain that it is equivalent to the two real equations $%
F_{1}=F_{2}=0$ with%
\begin{eqnarray*}
&&F_{1}\left( A,B,n,m\right) \\
&=&-\cos \left( \pi A\right) +\cos \left( \frac{4\pi m}{3}-\pi A\right)
-\cos \left( \frac{4\pi m}{3}+\pi A\right) \\
&&-\cos \left( \pi B\right) +\cos \left( \frac{4\pi n}{3}-\pi B\right) -\cos
\left( \frac{4\pi n}{3}+\pi B\right)
\end{eqnarray*}%
\begin{eqnarray*}
&&F_{2}\left( A,B,n,m\right) \\
&=&-\sin \left( \pi A\right) +\sin \left( \frac{4\pi m}{3}-\pi A\right)
-\sin \left( \frac{4\pi m}{3}+\pi A\right) \\
&&-\sin \left( \pi B\right) +\sin \left( \frac{4\pi n}{3}-\pi B\right) -\sin
\left( \frac{4\pi n}{3}+\pi B\right)
\end{eqnarray*}

The existence of the functions $\Theta _{1},\ \Theta _{2}$ is just a
consequence of the Implicit Funcion Theorem. Indeed, since we have $%
F_{1}\left( A_{0},B_{0},1,2\right) =F_{2}\left( A_{0},B_{0},1,2\right) =0$
we only need to check that $\frac{\partial \left( F_{1},F_{2}\right) }{%
\partial \left( n,m\right) }\left( A_{0},B_{0},1,2\right) \neq 0.$ This is
equivalent to proving that:%
\begin{equation*}
\det \left(
\begin{array}{cc}
-\sin \left( \frac{4\pi }{3}-\pi B_{0}\right) +\sin \left( \frac{4\pi }{3}%
+\pi B_{0}\right) & \cos \left( \frac{4\pi }{3}-\pi B_{0}\right) -\cos
\left( \frac{4\pi }{3}+\pi B_{0}\right) \\
-\sin \left( \frac{8\pi }{3}-\pi A_{0}\right) +\sin \left( \frac{8\pi }{3}%
+\pi A_{0}\right) & \cos \left( \frac{8\pi }{3}-\pi A_{0}\right) -\cos
\left( \frac{8\pi }{3}+\pi A_{0}\right)%
\end{array}%
\right) \neq 0
\end{equation*}

Using that $B_{0}=\left( 1-A_{0}\right) $ and elementary trigonometric
formulas we obtain that this condition equivalent to:%
\begin{equation*}
0\neq \det \left(
\begin{array}{cc}
2\sin \left( \pi \left( 1-A_{0}\right) \right) \cos \left( \frac{4\pi }{3}%
\right) & 2\sin \left( \frac{4\pi }{3}\right) \sin \left( \pi \left(
1-A_{0}\right) \right) \\
2\cos \left( \frac{8\pi }{3}\right) \sin \left( \pi A_{0}\right) & 2\sin
\left( \frac{8\pi }{3}\right) \sin \left( \pi A_{0}\right)%
\end{array}%
\right) =-2\sqrt{3}\sin ^{2}\left( \pi A_{0}\right)
\end{equation*}%
which holds for $A_{0}\in \left( 0,1\right) .$

In order to obtain the asymptotics (\ref{W5E3}), (\ref{W5E3b}) we argue as
follows. We write $n=1+\delta _{1},\ m=2+\delta _{2}.$ Using Taylor's to
approximate the equation $K_{m}\left( A\right) +K_{n}\left( B\right) =0$ for
$A=A_{0},\ B=B_{0}-\varepsilon ,$ we obtain the following approximation to
the linear order:\
\begin{eqnarray*}
0 &=&-\frac{8}{3}\left( -\frac{1}{2}+\frac{\sqrt{3}}{2}i\right) \sin \left(
\pi A_{0}\right) \delta _{2}+\left( \sin \left( \pi A_{0}\right) +\sqrt{3}%
\cos \left( \pi A_{0}\right) \right) \varepsilon \\
&&-\frac{8}{3}\left( -\frac{1}{2}-\frac{\sqrt{3}}{2}i\right) \sin \left( \pi
A_{0}\right) \delta _{1}
\end{eqnarray*}

The imaginary part of this equation implies, to the linear order, that $%
\delta _{1}=\delta _{2}.$ We then obtain the approximation $\delta
_{1}=\delta _{2}=-\frac{3}{8}\left( 1+\sqrt{3}\cot \left( \pi A_{0}\right)
\right) \varepsilon $ as $\varepsilon \rightarrow 0.$ This gives (\ref{W5E3}%
), (\ref{W5E3b}). This computation can be made fully rigorous by means of a
standard application of the Implicit Function Theorem. We will skip the
details.
\end{proof}

We can obtain now a representation formula for $Q\left( A,B,n,m\right) $
with $\left( A,B,n,m\right) \in \mathcal{K}$, $A+B<1.$

\begin{lemma}
\label{RepQ}Let $\left( A,B,n,m\right) \in \mathcal{K}$ with $A+B<1$.\
Suppose that $Q\left( A,B,n,m\right) $ is as in Lemma \ref{DefFunct}. Then
the following representation formula holds:%
\begin{equation}
Q\left( A,B,n,m\right) =-2\pi i\int_{1}^{\infty }\left( 1+e^{\frac{4\pi ni}{3%
}}+e^{\frac{4\pi mi}{3}}\right) \left( \frac{t}{t-1}\right) ^{A+B}\frac{dt}{%
t^{\frac{2\left( n+m\right) }{3}}}  \label{W5E4}
\end{equation}

\begin{remark}
Notice that $\frac{2\left( n+m\right) }{3}>1$ if $\left( A,B,n,m\right) \in
\mathcal{K}$. Therefore the integral on the right of (\ref{W5E4}) is well
defined.
\end{remark}
\end{lemma}

\begin{proof}
Under the assumptions of the Lemma, (\ref{W5E1}) holds. We can now deform
the contour of integration $\tilde{C}\in \mathcal{S}$ to a new contour given
by the line $L=\left\{ t\in \mathcal{S}:1+re^{-i\delta }:r\geq 0\right\} $
where $\delta >0$ is a small number. The deformation is made by means of a
family of contours which behave asymptotically for large $\left\vert
t\right\vert $ as the line $\left\{ t=re^{i\theta }\right\} $ with $\theta $
varying from $-3\pi $ to $-\delta .$ In the process of deformation the
deforming contours must cross any point $t=x$ with $x\in \left( 1,\infty
\right) .$ This gives a contribution due to Residue Theorem equal to $-2\pi
i\Psi \left( e^{-2\pi i}x,x\right) ,$ where the notation $\Psi \left(
e^{-2\pi i}x,x\right) $ indicates that the function $\Psi $ must be
evaluated at that particular point of the Riemann surface $\mathcal{S}$.
Then:%
\begin{equation}
Q\left( A,B,n,m\right) =-2\pi i\int_{1}^{\infty }\Psi \left( e^{-2\pi
i}x,x;;A,B,n,m\right) dx-\int_{1}^{\infty }dx\int_{L}dt\frac{\Psi \left(
t,x;A,B,n,m\right) }{\left( t-x\right) }  \label{W5E5}
\end{equation}

In order to compute the last integral in (\ref{W5E5}) we argue as follows.
Exchanging the role of $x$ and $t$ we obtain:%
\begin{equation*}
\int_{1}^{\infty }dx\int_{L}dt\frac{\Psi \left( t,x;A,B,n,m\right) }{\left(
t-x\right) }=\int_{1}^{\infty }dt\int_{L}dx\frac{\Psi \left(
x,t;A,B,n,m\right) }{\left( x-t\right) }
\end{equation*}

Using now that $\Psi \left( x,t;A,B,n,m\right) =\Psi \left(
t,x;A,B,n,m\right) $ we obtain:%
\begin{equation}
\int_{1}^{\infty }dx\int_{L}dt\frac{\Psi \left( t,x;A,B,n,m\right) }{\left(
t-x\right) }=\int_{1}^{\infty }dt\int_{L}dx\frac{\Psi \left(
t,x;A,B,n,m\right) }{\left( x-t\right) }  \label{W5E6}
\end{equation}

In order to transform the last integral in the original one we deform the
contour $L$ to $\left( 1,\infty \right) $ and the contour $\left( 1,\infty
\right) $ to $L.$ In the process of deformation one the contour $L$ must
cross the point $x=t$ yielding a contribution due to Residue Theorem. Then:%
\begin{eqnarray}
&&\int_{1}^{\infty }dt\int_{L}dx\frac{\Psi \left( t,x;A,B,n,m\right) }{%
\left( x-t\right) }  \label{W5E7} \\
&=&\int_{L}dt\int_{1}^{\infty }dx\frac{\Psi \left( t,x;A,B,n,m\right) }{%
\left( x-t\right) }+2\pi i\int_{1}^{\infty }\Psi \left( t,t;A,B,n,m\right) dt
\notag
\end{eqnarray}

Combining (\ref{W5E6}), (\ref{W5E7}) we obtain:%
\begin{equation*}
\int_{1}^{\infty }dx\int_{L}dt\frac{\Psi \left( t,x;A,B,n,m\right) }{\left(
t-x\right) }=\pi i\int_{1}^{\infty }\Psi \left( t,t;A,B,n,m\right) dt
\end{equation*}

Plugging this formula into (\ref{W5E5}) and relabelling the name of the
variable in the first integral of the right of (\ref{W5E5}) we obtain:%
\begin{equation*}
Q\left( A,B,n,m\right) =-2\pi i\int_{1}^{\infty }\Psi \left( e^{-2\pi
i}t,t;;A,B,n,m\right) dt-\pi i\int_{1}^{\infty }\Psi \left(
t,t;A,B,n,m\right) dt
\end{equation*}

We now use that, taking into account the analyticity of the function $\Phi
\left( t;a\right) =\left( \frac{t}{t-1}\right) ^{a}$ in $C\diagdown \left[
0,1\right] $ we have:
\begin{equation*}
\Psi \left( t,t;A,B,n,m\right) =2\left( \frac{t}{t-1}\right) ^{A+B}\frac{1}{%
t^{\frac{2\left( n+m\right) }{3}}},\ \ \arg \left( t\right) \in \left[
-\delta ,0\right]
\end{equation*}%
\begin{equation*}
\Psi \left( e^{-2\pi i}t,t;;A,B,n,m\right) =\left[ e^{\frac{4\pi ni}{3}}+e^{%
\frac{4\pi mi}{3}}\right] \left( \frac{t}{t-1}\right) ^{A+B}\frac{1}{t^{%
\frac{2\left( n+m\right) }{3}}}\ \ ,\ \ \arg \left( t\right) \in \left[
-\delta ,0\right]
\end{equation*}%
if $\delta >0$ is sufficiently small, whence (\ref{W5E4}) follows.
\end{proof}

We now take the limit of $Q\left( A,B,n,m\right) $ as $\left( A+B\right)
\rightarrow 1^{-}$.

\begin{lemma}
\label{LimitQ}Suppose that $0<A_{0}<1,\ 0<B_{0}<1$ satisfy $A_{0}+B_{0}=1.$
Then:%
\begin{equation*}
Q\left( A_{0},B_{0},1,2\right) =\pi ^{2}\left( 1+\sqrt{3}\cot \left( \pi
A_{0}\right) \right)
\end{equation*}
\end{lemma}

\begin{proof}
Lemmas \ref{DefFunct}, \ref{localK} imply:%
\begin{equation}
Q\left( A_{0},B_{0},1,2\right) =\lim_{\varepsilon \rightarrow 0^{+}}Q\left(
A_{0},B_{0}-\varepsilon ,\Theta _{1}\left( A_{0},B_{0}-\varepsilon \right)
,\Theta _{2}\left( A_{0},B_{0}-\varepsilon \right) \right)  \label{W5E8}
\end{equation}

Using (\ref{W5E3}), (\ref{W5E3b}) as well as the fact that $1+e^{\frac{4\pi i%
}{3}}+e^{\frac{8\pi i}{3}}=0$ we obtain:%
\begin{eqnarray}
\frac{1}{\varepsilon }\left( 1+e^{\frac{4\pi \Theta _{1}\left(
A_{0},B_{0}-\varepsilon \right) i}{3}}+e^{\frac{4\pi \Theta _{2}\left(
A_{0},B_{0}-\varepsilon \right) i}{3}}\right) &\rightarrow &-\frac{3\cdot
4\pi i\left( e^{\frac{4\pi i}{3}}+e^{\frac{8\pi i}{3}}\right) }{3\cdot 8}%
\left( 1+\sqrt{3}\cot \left( \pi A_{0}\right) \right)  \notag \\
&=&\frac{\pi i}{2}\left( 1+\sqrt{3}\cot \left( \pi A_{0}\right) \right) \
\label{W5E9}
\end{eqnarray}%
as $\varepsilon \rightarrow 0.$ Notice that we use $e^{\frac{4\pi i}{3}}+e^{%
\frac{8\pi i}{3}}=-1.$ We can now take the limit in (\ref{W5E8}) using (\ref%
{W5E4}). Notice that the main contribution to the integral is due to the
region where $t$ is close to one. We can then split the integral as $%
\int_{1}^{1+\delta }+\int_{1+\delta }^{\infty }$ with $\delta >0$ small. The
second term is bounded as $C_{\delta }\varepsilon $ for each $\delta >0$ due
to (\ref{W5E9}) and in the first we can approximate $t$ by $1$ in $t^{\frac{%
2\left( n+m\right) }{3}}$ and $t^{A+B}.$ Then, using again (\ref{W5E9}) we
obtain:\
\begin{eqnarray*}
&&\lim_{\varepsilon \rightarrow 0^{+}}Q\left( A_{0},B_{0}-\varepsilon
,\Theta _{1}\left( A_{0},B_{0}-\varepsilon \right) ,\Theta _{2}\left(
A_{0},B_{0}-\varepsilon \right) \right) \\
&=&-2\pi i\lim_{\varepsilon \rightarrow 0^{+}}\int_{1}^{1+\delta }\left(
1+e^{\frac{4\pi \Theta _{1}i}{3}}+e^{\frac{4\pi \Theta _{2}i}{3}}\right)
\left( \frac{t}{t-1}\right) ^{A_{0}+B_{0}-\varepsilon }\frac{dt}{t^{\frac{%
2\left( n+m\right) }{3}}} \\
&=&\pi ^{2}\left( 1+\sqrt{3}\cot \left( \pi A_{0}\right) \right)
\lim_{\varepsilon \rightarrow 0^{+}}\varepsilon \int_{1}^{1+\delta }\left(
\frac{1}{t-1}\right) ^{1-\varepsilon } \\
&=&\pi ^{2}\left( 1+\sqrt{3}\cot \left( \pi A_{0}\right) \right)
\end{eqnarray*}%
and the result follows.
\end{proof}

We can now compute $C_{\ast }$ and finish the Proof of Proposition \ref%
{Cstar}

\begin{proof}[End of the Proof of Proposition \protect\ref{Cstar}]
Using (\ref{W4E5}), (\ref{W4E9}) we obtain
\begin{equation*}
\lim_{R\rightarrow \infty }\int_{0}^{R}w\Delta _{\alpha }\left( w\right) dw=-%
\frac{\sin \left( \pi \alpha \right) }{3\pi }Q\left( 1+\alpha ,\beta +\frac{2%
}{3},1,2\right) .
\end{equation*}
Lemma \ref{LimitQ} yields
\begin{eqnarray*}
\lim_{R\rightarrow \infty }\int_{0}^{R}w\Delta _{\alpha }\left( w\right) dw
&=&-\frac{\pi \sin \left( \pi \alpha \right) }{3}\left( 1+\sqrt{3}\cot
\left( \pi \left( 1+\alpha \right) \right) \right) \\
&=&-\frac{\pi }{3}\left( \sin \left( \pi \alpha \right) +\sqrt{3}\cos \left(
\pi \alpha \right) \right) ,
\end{eqnarray*}
whence:
\begin{equation*}
\frac{C_{\ast }}{\left( 9\right) ^{\frac{2}{3}}}=\frac{\pi }{3}\left( \sin
\left( \pi \alpha \right) +\sqrt{3}\cos \left( \pi \alpha \right) \right)
-2\cos \left( \pi \left( \beta +\frac{1}{3}\right) \right) \log \left(
r\right)
\end{equation*}
\end{proof}

\bigskip

\subsection{Definition of some differential operators.}

Our goal now is to obtain suitable adjoint operators in the sense of
Definition \ref{FormAdj} for (\ref{S0E1}), (\ref{S0E2}), with the boundary
conditions (\ref{S5E3}), (\ref{S5E4}), (\ref{S5E7}), (\ref{T1E1}). These
operators will act over the class of continuous functions on a topological
space $X.$ The action of the operators will be given by a differential
operator $\mathcal{L}$ with suitable boundary conditions. In this Section we
give the precise definitions of $X$ and $\mathcal{L}$.

\begin{definition}
\label{DefX}We define as $X_{0}$ the set obtained identifying the subset of
points $\left[ 0,\infty \right) \times \left( -\infty ,\infty \right) $ such
that $\left( x,v\right) =\left( 0,-v\right) $ and $\left( x,v\right) =\left(
0,rv\right) $, $v>0.$ We then define $X=X_{0}\cup \left\{ \infty \right\} ,$
and we endowed it with the natural topology inherited from $\mathbb{R}^{2}$
complemented with the following set of neighbourhoods of the point $\infty :$%
\begin{equation*}
\mathcal{O}_{M}=\left\{ \left( x,v\right) \in \left[ 0,\infty \right) \times
\left( -\infty ,\infty \right) :v<-M\text{ or }v>rM\ \text{or\ }x>M\right\}
\ \ ,\ \ M>0
\end{equation*}
\end{definition}

The set $X$ is a topological compact set. The continuous functions of this
space can be identified with the bounded continuous functions $\varphi $ in $%
\left[ 0,\infty \right) \times \left( -\infty ,\infty \right) $ such that
\begin{equation}
\varphi \left( 0,-v\right) =\varphi \left( 0,rv\right) ,\ v>0  \label{comp1}
\end{equation}%
and such that the limit $\lim_{x+\left\vert v\right\vert \rightarrow \infty
}\varphi \left( x,v\right) $ exists. We will denote this set of functions as
$C\left( X\right) .$ Notice that a function $\varphi \in C\left( X\right) $
defines a function in $C\left( \mathcal{U}\right) $ satisfying (\ref{comp1}%
). We will use the same notation $\varphi $ to refer to both functions for
the sake of simplicity.

We need to introduce some local directionality in a neighbourhood of each
point of $X\diagdown \left\{ \left( 0,0\right) ,\infty \right\} $ in order
to compute directional limits.

\begin{definition}
\label{defLeft}Given two points $\left( x_{1},v_{1}\right) ,\left(
x_{2},v_{2}\right) \in X\diagdown \left[ \left\{ \left( x,0\right) :x\geq
0\right\} \cup \left\{ \infty \right\} \right] .$ We will say that $\left(
x_{1},v_{1}\right) $ is to the left of $\left( x_{2},v_{2}\right) $ and we
will write $\left( x_{1},v_{1}\right) \ll \left( x_{2},v_{2}\right) $ if $x_{1}%
\operatorname{sgn}\left( v_{1}\right) <x_{2}\operatorname{sgn}\left( v_{2}\right) .$\
\end{definition}

\bigskip Notice that the previous definition just means that $\left(
x_{1},v_{1}\right) \ll \left( x_{2},v_{2}\right) $ in one of the following
three cases: (i)\ If $v_{1}>0$\ and $v_{2}>0$\ we have $x_{1}<x_{2}.$ (ii) $%
v_{1}<0<v_{2}.$ (iii) If $v_{1}<0$\ and $v_{2}<0$\ we have $x_{1}>x_{2}.$

\begin{definition}
\label{LeftNei}Given a point $\left( x_{0},v_{0}\right) \in X\setminus
\left\{ \left( x,0\right) :x\geq 0\right\} \cup \left\{ \infty \right\} $
and a neighbourhood $\mathcal{B}$ of $\left( x_{0},v_{0}\right) $ in the
topological space $X$ we define the left neighbourhood $\mathcal{B}%
^{-}\left( x_{0},v_{0}\right) $ as:%
\begin{equation*}
\mathcal{B}^{-}\left( x_{0},v_{0}\right) =\left\{ \left( x,v\right) \in
\mathcal{B}:\left( x,v\right) \ll \left( x_{0},v_{0}\right) \right\}
\end{equation*}
\end{definition}

\begin{remark}
Notice that the neighbourhood $\mathcal{B}$ must be understood as a
neighbourhood of the topological space $X.$ In particular, if $\left(
x_{0},v_{0}\right) =\left( 0,v_{0}\right) $ any neighbourhood of $\left(
x_{0},v_{0}\right) $ contains points $\left( x,v\right) $ with $v>0$ and $%
v<0.$
\end{remark}

\begin{definition}
\label{LineLeft}We will say that $L\subset X$ is a vertical segment if it
has the form $\left\{ \left( x_{0},v\right) \in X: v\in \left( \alpha ,\beta \right) \right\} $ 
for some $x_{0}\geq 0,\ \alpha ,\beta \in \mathbb{R}$
with $\alpha \cdot \beta >0,$ $\alpha <\beta .$ Given two vertical segments $%
L_{1},\ L_{2}$ we will say that $L_{1}$ is to the left of $L_{2}$ if for any
$\left( x_{k},v_{k}\right) \in L_{k}$ with $k=1,2$ we have\ $\left(
x_{1},v_{1}\right) \ll \left( x_{2},v_{2}\right) .$ We will then write $%
L_{1}\ll L_{2}$. We will say that $L\subset X$ is a horizontal segment if it
has the form $\left\{ \left( x,v_{0}\right) \in X:x\in \left( \alpha ,\beta
\right) \right\} $ for some $v_{0}\in \mathbb{R}$ $\alpha ,\beta \in \mathbb{%
R}$ with $0\leq \alpha <\beta .$
\end{definition}

It will be convenient to define a suitable concept of convergence in the set
of segments.

\begin{definition}
\label{Distance}Given two segments $L_{1},\ L_{2}$ in $X$ we define a
distance between them as:%
\begin{equation}
\operatorname{dist}{}_{H}\left( L_{1},L_{2}\right) =\inf \left\{ \operatorname{dist}%
\left( \left( x,v\right) ,L_{2}\right) :\left( x,v\right) \in L_{1}\right\}
\label{distH}
\end{equation}
\end{definition}

The action of the operators $\Omega _{\sigma }$ on smooth functions
supported in $\mathcal{U}$ (cf. (\ref{S8E6})) is given by the differential
operator
\begin{equation}
\mathcal{L}=D_{v}^{2}+vD_{x}  \label{diffOp}
\end{equation}%
where the operator $\mathcal{L}$ will be defined in the sense of
distributions as indicated later. Nevertheless, the operators $\Omega
_{\sigma }$ will differ in the different cases (cf. (\ref{S5E3}), (\ref{S5E4}%
), (\ref{S5E6}), (\ref{T1E1})) in its domain of definition which will encode
the asymptotic behaviour of $\varphi $ near the singular point $\left(
x,v\right) =\left( 0,0\right) .$

We endow the set $C\left( X\right) $ with a Banach space structure using the
norm:
\begin{equation}
\left\Vert \varphi \right\Vert =\sup_{\left( x,v\right) \in X}|\varphi (x,v)|
\label{Norm}
\end{equation}

We need to impose suitable regularity and compatibility conditions in the
class of test functions in order to take into account the compatibility
conditions imposed in $X.$ Let $V$ be an open subset of $\ \mathcal{U}$. We
will consider functions $\zeta \in C\left( \bar{V}\right) $ satisfying:%
\begin{equation}
\zeta \left( 0,-v\right) =r^{2}\zeta \left( 0,rv\right) \ ;\ v>0\ ,\ \text{%
if }\left( 0,-v\right) ,\left( 0,rv\right) \in \bar{V}  \label{comp1t}
\end{equation}%
\begin{equation}
\text{There exist }\zeta _{x},\zeta _{vv}\in C\left( \bar{V}\right)
\label{comp3}
\end{equation}%
\begin{equation}
\zeta _{x}\left( 0,-v\right) =r^{2}\zeta _{x}\left( 0,rv\right) \ ;\ v>0%
\text{ if }\left( 0,-v\right) ,\left( 0,rv\right) \in \bar{V}\
\label{comp4}
\end{equation}%
\begin{equation}
\operatorname{supp}\left( \zeta \right) \cap \left[ \left\{ x+\left\vert
v\right\vert \geq R\right\} \cup \left\{ 0,0\right\} \right] =\varnothing \
\ \text{for some }R>0  \label{comp2}
\end{equation}

We can define a set of functions
\begin{equation}
\mathcal{F}\left( V\right) \mathcal{=}\left\{ \zeta \in C\left( \bar{V}%
\right) :\text{(\ref{comp1t}),\ (\ref{comp3}), (\ref{comp4}), (\ref{comp2})
hold}\right\}  \label{defF}
\end{equation}

We now define the action of the operator $\mathcal{L}$ in a subset of $%
C\left( X\right) .$

\begin{definition}
\label{LadjDef}Suppose that $W$ is any open subset of $X.$ Given $\varphi
\in C\left( W\right) ,$ we will say that $\mathcal{L}\varphi $ is defined if
there exists $w\in \mathcal{M} \left( W\right) $ such that for any $\zeta \in \mathcal{F%
}\left( \mathcal{U}\cap W\right) $ we have:%
\begin{equation}
\int_{\mathcal{U}\cap W}\varphi \mathcal{L}^{\ast }\left( \zeta \right)
dxdv=\int_{\mathcal{U}\cap W}w\zeta dxdv  \label{Ladjoint}
\end{equation}%
where $\mathcal{L}^{\ast }=D_{v}^{2}-vD_{x}.$ We will then write $w=\mathcal{%
L}\varphi .$ Given $V\subset \mathcal{U}$, we will say that $\mathcal{L}%
\varphi $ is defined in $V$ if there exists $w\in C\left( V\right) $
satisfying (\ref{comp1}) such that for any $\zeta \in \mathcal{F}\left(
\mathcal{U}\right) $ such that $\operatorname{supp}\left( \zeta \right) \cap
\left( \partial V\right) =\varnothing ,$ (\ref{Ladjoint}) holds.
\end{definition}


\subsection{Definition of the operators $\Omega _{\protect\sigma }.$ Domains $%
\mathcal{D}(\Omega _{\protect\sigma })$. \label{Operators}}

We now define some operators $\Omega _{\sigma }$ for the different boundary
conditions (nontrapping, trapping, partially trapping), described in Section \ref{Fsingular}, where $%
\sigma $ is the subindex labelling each set of boundary conditions. In all
the cases the operator $\Omega _{\sigma }$ acts over continuous functions
defined on a compact topological space $X$ in Definition \ref{DefX}.

We will assume that the functions $\varphi \in \mathcal{D}(\Omega _{\sigma
}) $ have the following asymptotic behaviour near the singular set:%
\begin{equation}
\varphi \left( x,v\right) =\varphi \left( 0,0\right) +\mathcal{A}\left(
\varphi \right) F_{\beta }\left( x,v\right) +\psi \left( x,v\right) ,\ \
\lim_{R\rightarrow 0}\frac{\sup_{x+\left\vert v\right\vert ^{3}=R}\left\vert
\psi \left( x,v\right) \right\vert }{R^{\beta }}=0  \label{phi_decom}
\end{equation}%
where $\mathcal{A}\left( \varphi \right) \in \mathbb{R}$ and with $F_{\beta
} $ as in (\ref{Fbeta}).

\subsubsection{The case $r<r_{c}.$ Trapping boundary conditions.\label%
{Absorbing}}

In this case we will define $\Omega _{t,sub}\varphi $ as follows. We
consider the domain:%
\begin{equation}
\mathcal{D}(\Omega _{t,sub})=\{\varphi ,\mathcal{L}\varphi \in C\left(
X\right) :\varphi \text{ satisfies (\ref{phi_decom}),}\lim_{\left(
x,v\right) \rightarrow \left( 0,0\right) }\left( \mathcal{L}\varphi \right)
\left( x,v\right) =0\}  \label{S7E1}
\end{equation}

We then define:

\begin{eqnarray}
\left( \Omega _{t,sub}\varphi \right) \left( x,v\right) &=&\left( \mathcal{L}%
\varphi \right) \left( x,v\right) \ ,\ \ \text{if}\ \ \left( x,v\right) \neq
\left( 0,0\right) ,\infty  \label{S7E1a} \\
\left( \Omega _{t,sub}\varphi \right) \left( 0,0\right) &=&\lim_{\left(
x,v\right) \rightarrow \left( 0,0\right) }\left( \mathcal{L}\varphi \right)
\left( x,v\right) =0\ \ ,\ \ \left( \Omega _{t,sub}\varphi \right) \left(
\infty \right) =\left( \mathcal{L}\varphi \right) \left( \infty \right)
\notag
\end{eqnarray}%
with $\mathcal{L}\varphi $ as in Definition \ref{LadjDef}. We remark that in
this case as well as in the following three cases, we have that $\mathcal{L}%
\varphi \in C\left( X\right) $ for the functions $\varphi $ in the domains
and then the limit $\lim_{\left( x,v\right) \rightarrow \left( 0,0\right)
}\left( \mathcal{L}\varphi \right) \left( x,v\right) $ exists.

\subsubsection{The case $r<r_{c}.$ Nontrapping boundary conditions.\label%
{Reflecting}}

We define $\Omega _{nt,sub}\varphi $ by means of (\ref{diffOp}) in the
domain:%
\begin{equation}
\mathcal{D}(\Omega _{nt,sub})=\left\{ \varphi ,\mathcal{L}\varphi \in
C\left( X\right) :\varphi \text{ satisfies (\ref{phi_decom})\ and\ }\mathcal{%
A}\left( \varphi \right) =0\right\}  \label{S7E2}
\end{equation}

We then define:

\begin{eqnarray}
\left( \Omega _{nt,sub}\varphi \right) \left( x,v\right) &=&\left( \mathcal{L%
}\varphi \right) \left( x,v\right) \ ,\ \ \text{if}\ \ \left( x,v\right)
\neq \left( 0,0\right) ,\ \infty  \label{S7E2a} \\
\left( \Omega _{nt,sub}\varphi \right) \left( 0,0\right) &=&\lim_{\left(
x,v\right) \rightarrow \left( 0,0\right) }\left( \mathcal{L}\varphi \right)
\left( x,v\right) \ \ ,\ \ \left( \Omega _{nt,sub}\varphi \right) \left(
\infty \right) =\left( \mathcal{L}\varphi \right) \left( \infty \right) \
\notag
\end{eqnarray}%
with $\mathcal{L}\varphi $ as in Definition \ref{LadjDef}.

\subsubsection{The case $r<r_{c}.$ Partially trapping boundary conditions.
\label{Mixed}}

Given any $\mu _{\ast }>0,$ we define $\Omega _{pt,sub}\varphi $ by means of
(\ref{diffOp}) in the domain:\
\begin{equation}
\mathcal{D}(\Omega _{pt,sub})=\left\{ \varphi ,\mathcal{L}\varphi \in
C\left( X\right) :\varphi \text{ satisfies (\ref{phi_decom}), there exists\ }%
\lim_{\left( x,v\right) \rightarrow \left( 0,0\right) }\left( \mathcal{L}%
\varphi \right) \left( x,v\right) =\mu _{\ast }\left\vert C_{\ast
}\right\vert \mathcal{A}\left( \varphi \right) \ \right\}  \label{S7E3}
\end{equation}%
where $C_{\ast }$ is as in (\ref{W1E5}) (cf. also Proposition \ref{Cstar}).
We will denote from now on $\lim_{\left( x,v\right) \rightarrow \left(
0,0\right) }\left( \mathcal{L}\varphi \right) \left( x,v\right) =\left(
\mathcal{L}\varphi \right) \left( 0,0\right) .$ We then define:%
\begin{eqnarray}
\left( \Omega _{pt,sub}\varphi \right) \left( x,v\right) &=&\left( \mathcal{L%
}\varphi \right) \left( x,v\right) \ ,\ \ \text{if}\ \ \left( x,v\right)
\neq \left( 0,0\right) ,\ \infty  \label{S7E3a} \\
\left( \Omega _{pt,sub}\varphi \right) \left( 0,0\right) &=&\left( \mathcal{L%
}\varphi \right) \left( 0,0\right) =\mu _{\ast }\left\vert C_{\ast
}\right\vert \mathcal{A}\left( \varphi \right) \ \ ,\ \ \left( \Omega
_{pt,sub}\varphi \right) \left( \infty \right) =\left( \mathcal{L}\varphi
\right) \left( \infty \right)  \notag
\end{eqnarray}

We will not make explicit the dependence of the operators $\Omega _{pt,sub}$
in $\mu _{\ast }$ for the sake of simplicity. Notice that trapping boundary
conditions reduce formally to the case $\mu _{\ast }=0$ and nontrapping
boundary conditions to the case $\mu _{\ast }=\infty .$

\subsubsection{The case $r>r_{c}.$\label{Supercrit}}

In this case we only consider the case of nontrapping boundary conditions.
We then define $\Omega _{\sup }\varphi $ by means of (\ref{diffOp}) in the
domain:%
\begin{equation}
\mathcal{D}(\Omega _{sup})=\left\{ \varphi ,\mathcal{L}\varphi \in C\left(
X\right) :\text{ there exists\ }\lim_{\left( x,v\right) \rightarrow \left(
0,0\right) }\left( \mathcal{L}\varphi \right) \left( x,v\right) \right\}
\label{S7E4}
\end{equation}

Notice that in this case the condition (\ref{phi_decom}) does not make sense
if\ $\varphi \in C\left( X\right) $ because $\beta <0.$

We then define

\begin{eqnarray}
\left( \Omega _{sup}\varphi \right) \left( x,v\right) &=&\left( \mathcal{L}%
\varphi \right) \left( x,v\right) \ ,\ \ \text{if}\ \ \left( x,v\right) \neq
\left( 0,0\right) ,\ \infty  \label{S7E4a} \\
\left( \Omega _{sup}\varphi \right) \left( 0,0\right) &=&\lim_{\left(
x,v\right) \rightarrow \left( 0,0\right) }\left( \mathcal{L}\varphi \right)
\left( x,v\right) ,\ \left( \Omega _{sup}\varphi \right) \left( \infty
\right) =\left( \mathcal{L}\varphi \right) \left( \infty \right)  \notag
\end{eqnarray}

\subsection{Formulation of the adjoint problems if $r<r_{c}.$}

In this Subsection we prove the following characterizations of the adjoint
operators $\mathbb{A}$ for the boundary conditions (\ref{S5E3}), (\ref{S5E4}%
), (\ref{S5E7}).

\begin{proposition}
\label{AdjAbs}Let $r<r_{c}.$ The operator $\Omega _{t,sub}$ defined in
Subsection \ref{Absorbing} is an\ adjoint operator $\mathbb{A}$ for the problem
(\ref{S0E1}), (\ref{S0E2}) with boundary conditions (\ref{S5E4}) in the
sense of Definition \ref{FormAdj}.
\end{proposition}

\begin{proposition}
\label{AdjRef}Let $r<r_{c}.$ The operator $\Omega _{nt,sub}$ defined in
Subsection \ref{Reflecting} is an adjoint operator $\mathbb{A}$ for the problem
(\ref{S0E1}), (\ref{S0E2}) with boundary conditions (\ref{S5E3}) in the
sense of Definition \ref{FormAdj}.
\end{proposition}

\begin{proposition}
\label{AdjMix}Let $r<r_{c}.$ The operator $\Omega _{pt,sub}$ defined in
Subsection \ref{Mixed} is an adjoint operator $\mathbb{A}$ for the problem (\ref%
{S0E1}), (\ref{S0E2}) with boundary conditions (\ref{S5E7}) in the sense of
Definition \ref{FormAdj}.
\end{proposition}

\begin{proof}[Proof of Propositions \protect\ref{AdjAbs}, \ref{AdjRef}, \ref%
{AdjMix}]

We use Definition \ref{FormAdj}. We will assume in the following that $P$ is
a smooth function outside the singular point satisfying (\ref{S0E1}), (\ref%
{S0E2}) and (\ref{Forigin}). Suppose first that $P$ satisfies also (\ref%
{S5E4}) and that $\varphi \left( \cdot ,t\right) \in \mathcal{D}(\Omega
_{t,sub})$ for all $t\in \left[ 0,T\right] .$ We then compute the integral:%
\begin{equation}
\mathcal{I}=\int_{\left[ 0,T\right] }\int_{\mathcal{U}}P\left( \partial
_{t}\varphi +\Omega _{t,sub}\varphi \right) dxdvdt  \label{S8E9}
\end{equation}

Notice that the assumptions on $\varphi $ in Definition \ref{FormAdj} as
well as the asymptotics (\ref{Forigin}) and the integrability of $P$ imply that $\left\vert P\right\vert
\left\vert \partial _{t}\varphi \right\vert $ and $\left\vert P\right\vert
\left\vert \Omega _{t,sub}\varphi \right\vert $ belong to $L^{1}\left( \left[
0,T\right] \times \mathbb{R}_{+}^{2}\right) .$ Therefore:%
\begin{eqnarray}
\mathcal{I} &=&\lim_{\delta \rightarrow 0}\int_{\left[ 0,T\right] }\int_{%
\mathcal{U}\diagdown \mathcal{R}_{\delta }}P\left( \partial _{t}\varphi
+\Omega _{t,sub}\varphi \right) dxdvdt  \notag \\
&=&\int_{\left[ 0,T\right] }\left( \lim_{\delta \rightarrow 0}\int_{\mathcal{%
U}\diagdown \mathcal{R}_{\delta }}P\left( \partial _{t}\varphi +\Omega
_{t,sub}\varphi \right) dxdv\right) dt\   \label{S9E1}
\end{eqnarray}%
with $\mathcal{R}_{\delta }$ as in (\ref{domain}) with $b=1.$ Notice that $%
\Omega _{t,sub}\varphi =\mathcal{L}\varphi $ for $\left( x,v\right) \neq
\left( 0,0\right) .$ Hypoellipticity in \cite{RS}, \cite{HVJ2} implies that $D_{v}\varphi $ is
continuous outside the singular point $\left( x,v\right) =\left( 0,0\right)
. $ Then:%
\begin{equation}
\int_{\mathcal{U}\diagdown \mathcal{R}_{\delta }}P\Omega _{t,sub}\varphi
dxdv=\mathcal{J}_{\delta }+\int_{\mathcal{U}\diagdown \mathcal{R}_{\delta
}}\varphi \left( D_{v}^{2}P-vD_{x}P\right) dxdv\   \label{S8E7}
\end{equation}%
with%
\begin{equation}
\mathcal{J}_{\delta }\left( t\right) =\int_{\partial \mathcal{R}_{\delta }}%
\left[ P\left( n_{v}D_{v}\varphi +n_{x}v\varphi \right) -\varphi D_{v}Pn_{v}%
\right] ds\   \label{S9E2}
\end{equation}%
where $ds$ is the arc-length of $\partial \mathcal{R}_{\delta }$ and the
normal vector is pointing towards $\mathcal{R}_{\delta }.$ Notice that in
the derivation of this formula we have used the existence of the derivative $%
D_{v}P.$ \ \ On the other hand, using that $\left\vert P\right\vert
\left\vert \partial _{t}\varphi \right\vert \in L^{1}\left( \left[ 0,T\right]
\times \mathbb{R}_{+}^{2}\right) $ we obtain:%
\begin{eqnarray}
\int_{\left[ 0,T\right] }\int_{\mathcal{U}\diagdown \mathcal{R}_{\delta
}}P\partial _{t}\varphi dxdvdt &=&\int_{\mathcal{U}\diagdown \mathcal{R}%
_{\delta }}\left( P\varphi \right) \left( \cdot ,T\right) dxdv-\int_{%
\mathcal{U}\diagdown \mathcal{R}_{\delta }}\left( P\varphi \right) \left(
\cdot ,0\right) dxdv  \notag \\
&&-\int_{\left[ 0,T\right] }\int_{\mathcal{U}\diagdown \mathcal{R}_{\delta
}}\varphi \partial _{t}Pdxdvdt  \label{S8E8}
\end{eqnarray}

Combining (\ref{S9E1}), (\ref{S8E7}), (\ref{S8E8}) and using also (\ref{S0E1}%
), (\ref{S0E2}) we obtain:\
\begin{equation*}
\mathcal{I}=\lim_{\delta \rightarrow 0}\int_{\left[ 0,T\right] }\mathcal{J}%
_{\delta }\left( t\right) dt+\lim_{\delta \rightarrow 0}\int_{\mathcal{U}%
\diagdown \mathcal{R}_{\delta }}\left( P\varphi \right) \left( \cdot
,T\right) dxdv-\lim_{\delta \rightarrow 0}\int_{\mathcal{U}\diagdown
\mathcal{R}_{\delta }}\left( P\varphi \right) \left( \cdot ,0\right) dxdv
\end{equation*}%
whence:\
\begin{equation}
\mathcal{I}=\lim_{\delta \rightarrow 0}\int_{\left[ 0,T\right] }\mathcal{J}%
_{\delta }\left( t\right) dt+\int_{\mathcal{U}}\left( P\varphi \right)
\left( \cdot ,T\right) dxdv-\int_{\mathcal{U}}\left( P\varphi \right) \left(
\cdot ,0\right) dxdv  \label{S9E4}
\end{equation}

We compute $\lim_{\delta \rightarrow 0}\int_{\left[ 0,T\right] }\mathcal{J}%
_{\delta }\left( t\right) dt$ as follows. Using the asymptotics (\ref%
{Forigin}), (\ref{phi_decom}) as well as (\ref{S5E4}) we can write:%
\begin{eqnarray}
\mathcal{J}_{\delta }\left( t\right) &=&a_{-\frac{2}{3}}\left( t\right)
\varphi \left( 0,0,t\right) \int_{\partial \mathcal{R}_{\delta }}\left[
\left( n_{x}v\right) G_{-\frac{2}{3}}-n_{v}D_{v}G_{-\frac{2}{3}}\right] ds
\notag \\
&&+\varphi \left( 0,0,t\right) \int_{\partial \mathcal{R}_{\delta }}\left[
n_{x}vQ-D_{v}Qn_{v}\right] ds+  \notag \\
&&+\int_{\partial \mathcal{R}_{\delta }}\left[ P\left(
n_{v}D_{v}W+n_{x}vW\right) -WD_{v}Pn_{v}\right] ds\   \label{S9E3}
\end{eqnarray}%
where $Q=P-a_{-\frac{2}{3}}\left( t\right) G_{-\frac{2}{3}}$ and $W=\varphi
-\varphi \left( 0,0,t\right) .$ We define functions $Q_{R}\left( x,v\right) =%
\frac{Q\left( Rx,R^{\frac{1}{3}}v\right) }{R^{-\frac{2}{3}}},$ $W_{R}\left(
x,v\right) =\frac{W\left( Rx,R^{\frac{1}{3}}v\right) }{R^{\beta }}.$ Using (%
\ref{Forigin}) and (\ref{phi_decom}) we obtain that $\left\vert Q_{R}\left(
x,v\right) \right\vert $ and $\left\vert W_{R}\left( x,v\right) \right\vert $
are bounded in $\mathcal{R}_{2}\diagdown \mathcal{R}_{\frac{1}{2}}$.
Hypoellipticity property in \cite{RS}, \cite{HVJ2} combined with (\ref{Forigin}) then imply that $\sup_{%
\mathcal{R}_{\frac{3}{2}}\diagdown \mathcal{R}_{\frac{2}{3}}}\left\vert
D_{v}Q_{R}\right\vert \rightarrow 0$ as $R\rightarrow 0.$ On the other hand, Hypoellipticity property
in \cite{RS}, \cite{HVJ2} and (\ref{phi_decom}) yield $\sup_{\mathcal{R}_{\frac{3}{%
2}}\diagdown \mathcal{R}_{\frac{2}{3}}}\left\vert D_{v}W_{R}\right\vert \leq
h\left( t\right) ,$ with $\int_{\left[ 0,T\right] }h\left( t\right)
dt<\infty .$ The definition of $Q_{R},\ W_{R}$ then yields:%
\begin{equation}
\lim_{R\rightarrow 0}\sup_{\frac{R}{2}\leq x+\left\vert v\right\vert
^{3}\leq 2R}\frac{\left\vert D_{v}Q\right\vert }{R^{-\frac{2}{3}-\frac{1}{3}}%
}=0\ \ ,\ \ \ \int_{\left[ 0,T\right] }\left[ \sup_{\frac{R}{2}\leq
x+\left\vert v\right\vert ^{3}\leq 2R}\left( \left\vert W\right\vert +R^{%
\frac{1}{3}}\left\vert D_{v}W\right\vert \right) \right] dt\leq CR^{\beta }
\label{DerVest}
\end{equation}

Therefore, the last two integral terms in (\ref{S9E3}) converge to zero
(after integrating in time) and we have:%
\begin{equation*}
\lim_{\delta \rightarrow 0}\int_{\left[ 0,T\right] }\mathcal{J}_{\delta
}\left( t\right) dt=\left( \int_{\left[ 0,T\right] }a_{-\frac{2}{3}}\left(
t\right) \varphi \left( 0,0,t\right) dt\right) \lim_{\delta \rightarrow
0}\int_{\partial \mathcal{R}_{\delta }}\left[ \left( n_{x}v\right) G_{-\frac{%
2}{3}}-n_{v}D_{v}G_{-\frac{2}{3}}\right] ds
\end{equation*}

The last limit can be computed using Proposition \ref{LimitFluxes}. Then:\
\begin{equation*}
\lim_{\delta \rightarrow 0}\int_{\left[ 0,T\right] }\mathcal{J}_{\delta
}\left( t\right) dt=-9^{\frac{2}{3}}\left[ \log \left( r\right) +\frac{\pi }{%
\sqrt{3}}\right] \left( \int_{\left[ 0,T\right] }a_{-\frac{2}{3}}\left(
t\right) \varphi \left( 0,0,t\right) dt\right)
\end{equation*}

We then obtain, using (\ref{S9E4}):%
\begin{eqnarray*}
\mathcal{I} &=&-9^{\frac{2}{3}}\left[ \log \left( r\right) +\frac{\pi }{%
\sqrt{3}}\right] \left( \int_{\left[ 0,T\right] }a_{-\frac{2}{3}}\left(
t\right) \varphi \left( 0,0,t\right) dt\right) \\
&&+\int_{\mathcal{U}}\left( P\varphi \right) \left( \cdot ,T\right)
dxdv-\int_{\mathcal{U}}\left( P\varphi \right) \left( \cdot ,0\right) dxdv
\end{eqnarray*}

Using now (\ref{S8E1}) and (\ref{meas2}) we obtain:%
\begin{equation*}
\mathcal{I}=\left( \int_{\left[ 0,T\right] }\frac{dm\left( t\right) }{dt}%
\varphi \left( 0,0,t\right) dt\right) +\int_{\mathcal{U}}\left( P\varphi
\right) \left( \cdot ,T\right) dxdv-\int_{\mathcal{U}}\left( P\varphi
\right) \left( \cdot ,0\right) dxdv
\end{equation*}%
whence, using (\ref{meas1}) and (\ref{S8E9}) we obtain (\ref{Ad1}). Notice
that we use also that $\varphi _{t}\left( 0,0,t\right) =-\mathcal{L}\varphi
\left( 0,0,t\right) =0$ due to (\ref{S7E1}). This concludes the Proof of
Proposition \ref{AdjAbs}.

We now consider the case of nontrapping boundary conditions (cf. Proposition %
\ref{AdjRef}). Arguing similarly we obtain formula (\ref{S9E4}) with:%
\begin{equation*}
\mathcal{I}=\int_{\left[ 0,T\right] }\int_{\mathcal{U}}P\left( \partial
_{t}\varphi +\Omega _{nt,sub}\varphi \right) dxdvdt
\end{equation*}

In order to compute $\lim_{\delta \rightarrow 0}\int_{\left[ 0,T\right] }%
\mathcal{J}_{\delta }\left( t\right) dt$ in this case we use (\ref{Forigin}%
), (\ref{phi_decom}) and (\ref{S5E3}). We then obtain, instead of (\ref{S9E3}%
):%
\begin{eqnarray}
\mathcal{J}_{\delta }\left( t\right) &=&a_{\alpha }\left( t\right) \varphi
\left( 0,0,t\right) \int_{\partial \mathcal{R}_{\delta }}\left[ \left(
n_{x}v\right) G_{\alpha }-n_{v}D_{v}G_{\alpha }\right] ds  \notag \\
&&+\varphi \left( 0,0,t\right) \int_{\partial \mathcal{R}_{\delta }}\left[
n_{x}vQ-D_{v}Qn_{v}\right] ds  \notag \\
&&+\int_{\partial \mathcal{R}_{\delta }}\left[ P\left(
n_{v}D_{v}W+n_{x}vW\right) -WD_{v}Pn_{v}\right] ds  \label{S9E5}
\end{eqnarray}%
where now $Q=P-a_{\alpha }\left( t\right) G_{\alpha }$ and $W=\varphi
-\varphi \left( 0,0,t\right) .$ Notice that the boundary condition $\mathcal{%
A}\left( \varphi \right) =0$ in (\ref{S7E2}) implies, arguing as in the
Proof of (\ref{DerVest}):%
\begin{equation}
\lim_{R\rightarrow 0}\sup_{\frac{R}{2}\leq x+\left\vert v\right\vert
^{3}\leq 2R}\left( \frac{\left\vert Q\right\vert }{R^{-\frac{2}{3}}}+\frac{%
\left\vert D_{v}Q\right\vert }{R^{-\frac{2}{3}-\frac{1}{3}}}\right) =0\
,\lim_{R\rightarrow 0}\int_{\left[ 0,T\right] }\left[ \sup_{\frac{R}{2}\leq
x+\left\vert v\right\vert ^{3}\leq 2R}\left( \frac{W}{R^{\beta }}+\frac{%
\left\vert D_{v}W\right\vert }{R^{\beta -\frac{1}{3}}}\right) \right] dt=0\
\   \label{S9E6}
\end{equation}

Then, the two last integrals in (\ref{S9E5}) tend to zero. On the other
hand, the remaining one vanishes due to Proposition \ref{FluxesGalpha}.
Therefore $\lim_{\delta \rightarrow 0}\int_{\left[ 0,T\right] }\mathcal{J}%
_{\delta }\left( t\right) dt=0.$ Using then (\ref{S5E3}), (\ref{S8E1}) and (%
\ref{meas2}) it follows that $m\left( t\right) =0,$ whence (\ref{Ad1})
follows. This shows Proposition \ref{AdjRef}.

Finally we consider the case of Partially Trapping Boundary Conditions (cf.
Proposition \ref{AdjMix}). In this case we obtain (\ref{S9E4}) with:%
\begin{equation}
\mathcal{I}=\int_{\left[ 0,T\right] }\int_{\mathcal{U}}P\left( \partial
_{t}\varphi +\Omega _{pt,sub}\varphi \right) dxdvdt \label{I}
\end{equation}

We now have:%
\begin{eqnarray}
\mathcal{J}_{\delta }\left( t\right) &=&a_{\alpha }\left( t\right) \varphi
\left( 0,0,t\right) \int_{\partial \mathcal{R}_{\delta }}\left[ \left(
n_{x}v\right) G_{\alpha }-nD_{v}G_{\alpha }\right] ds  \notag \\
&&+a_{-\frac{2}{3}}\left( t\right) \varphi \left( 0,0,t\right)
\int_{\partial \mathcal{R}_{\delta }}\left[ \left( n_{x}v\right) G_{-\frac{2%
}{3}}-nD_{v}G_{-\frac{2}{3}}\right] ds\   \notag \\
&&+\varphi \left( 0,0,t\right) \int_{\partial \mathcal{R}_{\delta }}\left[
n_{x}vQ-D_{v}Qn_{v}\right] ds+  \notag \\
&&+\int_{\partial \mathcal{R}_{\delta }}\left[ P\left(
n_{v}D_{v}W+n_{x}vW\right) -WD_{v}Pn_{v}\right] ds\   \label{S9E7}
\end{eqnarray}%
with $Q=P-a_{\alpha }\left( t\right) G_{\alpha }-a_{-\frac{2}{3}}\left(
t\right) G_{-\frac{2}{3}},$ $W=\varphi -\varphi \left( 0,0,t\right) .$
Arguing as in the previous cases, by means of a rescaling argument, we
obtain:%
\begin{eqnarray}
\lim_{R\rightarrow 0}\sup_{\frac{R}{2}\leq x+\left\vert v\right\vert
^{3}\leq 2R}\left( \frac{\left\vert Q\right\vert }{R^{-\frac{2}{3}}}+\frac{%
\left\vert D_{v}Q\right\vert }{R^{-\frac{2}{3}-\frac{1}{3}}}\right) &=&0\ \
\label{S9E8} \\
\ \lim_{R\rightarrow 0}\sup_{\frac{R}{2}\leq x+\left\vert v\right\vert
^{3}\leq 2R}\int_{\left[ 0,T\right] }\left( \frac{W}{R^{\beta }}+\frac{%
\left\vert D_{v}W\right\vert }{R^{\beta -\frac{1}{3}}}\right) &\leq &C\
\notag
\end{eqnarray}

The first integral on the right-hand side of (\ref{S9E7}) vanishes due to
Proposition \ref{FluxesGalpha}. The third integral tends to zero as $\delta
\rightarrow 0$ due to (\ref{S9E8}). Using Proposition \ref{LimitFluxes} we
obtain:%
\begin{eqnarray}
\lim_{\delta \rightarrow 0}\mathcal{J}_{\delta }\left( t\right) &=&-9^{\frac{%
2}{3}}\left[ \log \left( r\right) +\frac{\pi }{\sqrt{3}}\right] a_{-\frac{2}{%
3}}\left( t\right) \varphi \left( 0,0,t\right)  \notag \\
&&+\lim_{\delta \rightarrow 0}\int_{\partial \mathcal{R}_{\delta }}\left[
P\left( n_{v}D_{v}W+n_{x}vW\right) -WD_{v}Pn_{v}\right] ds  \label{S9E9}
\end{eqnarray}

We now notice that:%
\begin{equation}
\lim_{\delta \rightarrow 0}\int_{\partial \mathcal{R}_{\delta }}\left[
P\left( n_{v}D_{v}W+n_{x}vW\right) -WD_{v}Pn_{v}\right] ds=C_{\ast
}a_{\alpha }\left( t\right) \mathcal{A}\left( \varphi \right) \
\label{S9E9a}
\end{equation}%
where $C_{\ast }$ is as in (\ref{W1E5}). We recall that $C_{\ast }$ has been
computed in Proposition \ref{Cstar}. Therefore, using (\ref{S9E9a}) in (\ref%
{S9E9}) we obtain:%
\begin{equation*}
\lim_{\delta \rightarrow 0}\mathcal{J}_{\delta }\left( t\right) =-9^{\frac{2%
}{3}}\left[ \log \left( r\right) +\frac{\pi }{\sqrt{3}}\right] a_{-\frac{2}{3%
}}\left( t\right) \varphi \left( 0,0,t\right) +C_{\ast }a_{\alpha }\left(
t\right) \mathcal{A}\left( \varphi \right)
\end{equation*}

Combining (\ref{S8E1}), (\ref{meas2}), (\ref{S9E4}), (\ref{S9E9}) we obtain:%
\begin{eqnarray*}
\mathcal{I} &=&\int_{\left[ 0,T\right] }\frac{dm\left( t\right) }{dt}\varphi
\left( 0,0,t\right) dt+C_{\ast }\int_{\left[ 0,T\right] }a_{\alpha }\left(
t\right) \mathcal{A}\left( \varphi \right) dt \\
&&+\int_{\mathcal{U}}\left( P\varphi \right) \left( \cdot ,T\right)
dxdv-\int_{\mathcal{U}}\left( P\varphi \right) \left( \cdot ,0\right) dxdv
\end{eqnarray*}

Using then the definition of the measure $f$ in (\ref{meas1}) as well as the
fact that $\partial _{t}\varphi \left( 0,0,t\right) =-\mathcal{L}\varphi
\left( 0,0,t\right) $ we obtain:%
\begin{equation*}
\mathcal{I}=\int_{\left[ 0,T\right] }m\left( t\right) \mathcal{L}\varphi
\left( 0,0,t\right) dt+C_{\ast }\int_{\left[ 0,T\right] }a_{\alpha }\left(
t\right) \mathcal{A}\left( \varphi \right) dt+\int_{\mathcal{V}}f\left(
dxdv,T\right) -\int_{\mathcal{V}}f\left( dxdv,0\right)
\end{equation*}

Using now (\ref{S5E6}) and (\ref{S7E3}) we obtain:%
\begin{eqnarray*}
\mathcal{I} &=&-\mu _{\ast }C_{\ast }\int_{\left[ 0,T\right] }m\left(
t\right) \mathcal{A}\left( \varphi \right) dt+\mu _{\ast }C_{\ast }\int_{%
\left[ 0,T\right] }m\left( t\right) \mathcal{A}\left( \varphi \right) dt \\
&&+\int_{\mathcal{V}}f\left( dxdv,T\right) -\int_{\mathcal{V}}f\left(
dxdv,0\right) \\
&=&\int_{\mathcal{V}}f\left( dxdv,T\right) -\int_{\mathcal{V}}f\left(
dxdv,0\right)
\end{eqnarray*}%
where we used the fact that $C_{\ast }<0$ (cf. Lemma \ref{CstarSign}). Then,
using also the definition of $\mathcal{I}$ in (\ref{I}) we obtain :%
\begin{equation*}
\int_{\left[ 0,T\right] }\int_{\mathcal{U}}P\left( \partial _{t}\varphi
+\Omega _{a,m}\varphi \right) dxdvdt=\int_{\mathcal{V}}f\left( dxdv,T\right)
-\int_{\mathcal{V}}f\left( dxdv,0\right)
\end{equation*}%
whence (\ref{Ad1}) follows. This concludes the Proof of Proposition \ref%
{AdjMix}.
\end{proof}

\subsection{Formulation of the adjoint problem if $r>r_{c}$}

\begin{proposition}
\label{AdjSupCrt}Let $r>r_{c}.$ The operator $\Omega _{\sup }$ defined in
Subsection \ref{Supercrit} is an adjoint operator $\mathbb{A}$ for the problem (%
\ref{S0E1}), (\ref{S0E2}) with boundary conditions (\ref{T1E1}) in the sense
of Definition \ref{FormAdj}.
\end{proposition}

\begin{proof}
It is similar to the proof of Propositions \ref{AdjAbs}, \ref{AdjRef}, \ref%
{AdjMix}. We define
\begin{equation}
\mathcal{I}=\int_{\left[ 0,T\right] }\int_{\mathcal{U}}P\left( \partial
_{t}\varphi +\Omega _{sup}\varphi \right) dxdvdt  \label{T8E3}
\end{equation}

By assumption $P$ satisfies (\ref{Forigin}) with (\ref{T1E1}). We define $%
Q=P-a_{\alpha }\left( t\right) G_{\alpha }$ and $W=\varphi -\varphi \left(
0,0,t\right) .$ Then, arguing as in the previous proof we obtain (\ref{S9E4}%
) with:%
\begin{eqnarray}
\mathcal{J}_{\delta }\left( t\right) &=&a_{\alpha }\left( t\right) \varphi
\left( 0,0,t\right) \int_{\partial \mathcal{R}_{\delta }}\left[ \left(
n_{x}v\right) G_{\alpha }-n_{v}D_{v}G_{\alpha }\right] ds+\   \notag \\
&&+\varphi \left( 0,0,t\right) \int_{\partial \mathcal{R}_{\delta }}\left[
n_{x}vQ-D_{v}Qn_{v}\right] ds+  \notag \\
&&+\int_{\partial \mathcal{R}_{\delta }}\left[ P\left(
n_{v}D_{v}W+n_{x}vW\right) -WD_{v}Pn_{v}\right] ds\ \   \label{T8E1}
\end{eqnarray}

Notice that in this case we have $\alpha >-\frac{2}{3}$ and due to (\ref%
{T1E1}) we have, arguing as in the previous case:%
\begin{equation}
\lim_{R\rightarrow 0}\sup_{\frac{R}{2}\leq x+\left\vert v\right\vert
^{3}\leq 2R}\left( \frac{\left\vert Q\right\vert }{R^{-\frac{2}{3}}}+\frac{%
\left\vert D_{v}Q\right\vert }{R^{-\frac{2}{3}-\frac{1}{3}}}\right) =0\
,\lim_{R\rightarrow 0}\int_{\left[ 0,T\right] }\left[ \sup_{\frac{R}{2}\leq
x+\left\vert v\right\vert ^{3}\leq 2R}\left( \left\vert W\right\vert +\frac{%
\left\vert D_{v}W\right\vert }{R^{-\frac{1}{3}}}\right) \right] dt=0\
\label{T8E2}
\end{equation}

Notice that, since $\beta <0$ in this case, the only information that we
have about $\left\vert W\right\vert $ near the origin is that it converges
to zero, plus the estimates for the derivatives that can be obtained by
rescaling. The first integral on the right of (\ref{T8E1}) vanishes due to
Proposition \ref{FluxesGalpha}. The second converges to zero as $\delta
\rightarrow 0$ due to (\ref{T8E2}) and the third one can be estimated, using
the estimates for $P$ as well as (\ref{T8E2}) as $C\left( \delta \right)
^{\alpha +\frac{2}{3}}.$ Using that $\alpha >-\frac{2}{3}$ it then follows
that the this integral converges to zero as $\delta \rightarrow 0.$ Taking
the limit of (\ref{S9E4}) as $\delta \rightarrow 0$ we arrive at:
\begin{equation*}
\mathcal{I}=\int_{\mathcal{U}}\left( P\varphi \right) \left( \cdot ,T\right)
dxdv-\int_{\mathcal{U}}\left( P\varphi \right) \left( \cdot ,0\right) dxdv
\end{equation*}

Using then that $m\left( t\right) =0$ (cf. (\ref{S8E1}) and (\ref{meas1}))
we obtain (\ref{Ad1}) and the result follows.
\end{proof}


%
%
%

\section{Weak solutions for the original problem.\label{weakSolDef}}

We will define suitable measure valued solutions of the problem
(\ref{S0E1})-(\ref{S0E3}) by means of the corresponding adjoint problems. 
To this end, we argue by duality. We will use the index $%
\sigma $ to denote each of the four cases considered in Subsection \ref%
{Operators} of Section \ref{AdjProblems}, namely, in the case of subcritical
values of $r,$ we can use trapping, nontrapping or partially trapping
boundary conditions. We will consider also the supercritical case. The
following definition will be used to define measure valued solutions of the
problem (\ref{S0E1})-(\ref{S0E3}) with all the boundary conditions
considered above.

\begin{definition}
\label{weakSol}Given $P_{0}\in \mathcal{M}_{+}\left( X\right) $ we define a
measure valued function $P_{\sigma }\in C\left( \left[ 0,\infty \right) :%
\mathcal{M}_{+}\left( X\right) \right) $ by means of:%
\begin{equation}
\int \varphi \left( dP_{\sigma }\left( \cdot ,t\right) \right) =\int \left(
S_{\sigma }\left( t\right) \varphi \right) \left( dP_{0}\right) \ \ ,\ \
t\geq 0  \label{A3}
\end{equation}%
for any $\varphi \in C\left( \left[ 0,\infty \right) :C\left( X\right)
\right) $
\end{definition}

\begin{remark}
Notice that the notation in (\ref{A3}) must be understood as follows. Let $%
\psi \left( \cdot ,t\right) =S_{\sigma }\left( t\right) \varphi .$ Then, the
right-hand side of (\ref{A3}) is equivalent to $\int \psi \left( \cdot
,t\right) dP_{0}.$ The left-hand side of (\ref{A3}) is just $\int \varphi
\left( \cdot \right) dP_{\sigma }\left( t\right) .$
\end{remark}

It is convenient to write in detail the weak formulation of (\ref{S0E1})-(%
\ref{S0E3}) satisfied by each of the measures $P_{\sigma }.$

\begin{definition}
\label{WAC}Suppose that $0<r<r_{c}$ and $P_{0}\in \mathcal{M}_{+}\left(
X\right) .$ We will say that $P\in C\left( \left[ 0,\infty \right) :\mathcal{%
M}_{+}\left( X\right) \right) $ is a weak solution of (\ref{S0E1})-(\ref%
{S0E3}) with trapping boundary conditions if for any $T\in \left[ 0,\infty
\right) $ and $\varphi \ $such that, for any time $t\in \left[ 0,T\right] ,$
$\varphi \left( \cdot ,t\right) ,\mathcal{L}\varphi \left( \cdot ,t\right)
,\varphi _{t}\left( \cdot ,t\right) \in C\left( X\right) ,\ \varphi \left(
\cdot ,t\right) $ satisfies (\ref{phi_decom})\ and such that $\mathcal{L}%
\varphi \left( 0,0,t\right) =0$ for any $t\geq 0,$ the following identity
holds:%
\begin{equation}
\int_{0}^{T}\int_{X\diagdown \left\{ \left( 0,0\right) \right\} }\left[
\varphi _{t}\left( dP\left( \cdot ,t\right) ,t\right) +\mathcal{L}\varphi
\left( dP\left( \cdot ,t\right) ,t\right) \right] +\int_{X}\varphi \left(
dP_{0}\left( \cdot \right) ,0\right) -\int_{X}\varphi \left( dP\left( \cdot
,T\right) ,T\right) =0,  \label{S7E5}
\end{equation}%
\ where $\mathcal{L}$ is as in (\ref{diffOp}).
\end{definition}

\begin{definition}
\label{WRC}Suppose that $0<r<r_{c}$ and $P_{0}\in \mathcal{M}_{+}\left(
X\right) .$ We will say that $P\in C\left( \left[ 0,\infty \right) :\mathcal{%
M}_{+}\left( X\right) \right) $ is a weak solution of (\ref{S0E1})-(\ref%
{S0E3}) with nontrapping boundary conditions if for any $T\in \left[
0,\infty \right) $ and $\varphi $ such that, for any time $t\in \left[ 0,T%
\right] ,$ $\varphi \left( \cdot ,t\right) ,\mathcal{L}\varphi \left( \cdot
,t\right) ,\varphi _{t}\left( \cdot ,t\right) \in C\left( X\right) ,\
\varphi \left( \cdot ,t\right) $ satisfies (\ref{phi_decom})\ and\ $\mathcal{%
A}\left( \varphi \right) \left( \cdot ,t\right) =0$ the identity (\ref{S7E5})
holds,\ where $\mathcal{L}$ is as in (\ref{diffOp}).
\end{definition}

\begin{definition}
\label{WMC}Suppose that $0<r<r_{c}$ and $P_{0}\in \mathcal{M}_{+}\left(
X\right) .$ We will say that $P\in C\left( \left[ 0,\infty \right) :\mathcal{%
M}_{+}\left( X\right) \right) $ is a weak solution of (\ref{S0E1})-(\ref%
{S0E3}) with partially trapping boundary conditions if for any $T\in \left[
0,\infty \right) $ and and $\varphi $ such that, for any time $t\in \left[
0,T\right] ,$ $\varphi \left( \cdot ,t\right) ,\mathcal{L}\varphi \left(
\cdot ,t\right) ,\varphi _{t}\left( \cdot ,t\right) \in C\left( X\right) ,\
\varphi \left( \cdot ,t\right) $ satisfies (\ref{phi_decom})\ and\ $\mathcal{%
L}\varphi \left( 0,0,t\right) =\mu _{\ast }\left\vert C_{\ast }\right\vert
\mathcal{A}\left( \varphi \right) ,$ the identity (\ref{S7E5}) holds, where $%
\mathcal{L}$ is as in (\ref{diffOp}).
\end{definition}

\begin{definition}
\label{WSuper}Suppose that $r>r_{c}$ and $P_{0}\in \mathcal{M}_{+}\left(
X\right) .$ We will say that $P\in C\left( \left[ 0,\infty \right) :\mathcal{%
M}_{+}\left( X\right) \right) $ is a weak solution of (\ref{S0E1})-(\ref%
{S0E3}) for supercritical boundary conditions if for any $T\in \left[
0,\infty \right) $ and $\varphi \in C_{c}^{2}\left( \left[ 0,T\right) \times
X\right) $ the identity (\ref{S7E5}) holds, where $\mathcal{L}$ is as in (%
\ref{diffOp}).
\end{definition}

We state the following Theorem, which will be proved rigorously in \cite{HVJ2}.

\begin{theorem}
Suppose that we define measures $P_{t,sub},\ P_{nt,sub},\ P_{pt,sub},\
P_{sup}\in C\left( \left[ 0,\infty \right) :\mathcal{M}_{+}\left( X\right)
\right) $ as in Definition \ref{weakSol} and initial datum $P_{0}$. Then,
these measures can be decomposed as:%
\begin{equation}
dP_{\sigma }\left( \cdot ,t\right) =m_{\sigma }\left( t\right) \delta
_{\left( 0,0\right) }+p_{\sigma }\left( x,v,t\right) dxdv\ \ ,\ \ \ t\geq 0\
\label{A4}
\end{equation}%
where for each $\sigma $ we have $p_{\sigma }\left( \cdot ,t\right) \in
L^{1}\left( X\right) $ and $m_{\sigma }\left( t\right) \geq 0.$ If $%
m_{\sigma }\left( 0\right) =0$ we have also $m_{nt,sub}\left( t\right)
=m_{sup}\left( t\right) =0$ for any $t\geq 0.$ Moreover, if $P_{0}$ is not
identically zero, we have also $m_{t,sub}\left( t\right) >0,\
m_{pt,sub}\left( t\right) >0$ for any $t>0.$

The functions $p_{\sigma }\left( \cdot ,t\right) $ are infinitely
differentiable for $\left( x,v\right) \neq \left( 0,0\right) $ and they
satisfy (\ref{S0E2}).
\end{theorem}



\noindent{\bf Acknowledgements.} 
The authors would like to thank the Hausdorff Center for Mathematical
Sciences of the University of Bonn and Pohang Mathematics Institute where
part of this work was done. HJH is partly supported by the Basic Science
Research Program (2015R1A2A2A0100251) through the National
Research Foundation of Korea (NRF). JJ is supported in part by NSF grants
DMS-1608492 and DMS-1608494. The authors acknowledge support through the CRC
1060 The mathematics of emergent effects at the University of Bonn, that is
funded through the German Science Foundation (DFG). We thank Seongwon Lee
for helping with figures in the paper.




\end{document}